\documentstyle[12pt,suthesis]{report}
\newtheorem{theorem}{Theorem}
\newtheorem{lemma}{Lemma}

\newtheorem{definition}{Definition}
\newtheorem{conjecture}{Conjecture}
\newtheorem{condition}{Condition}
\newtheorem{proposition}{Proposition}
\newcommand{\real}{{\bf R}}
\newcommand{\ddt}{\frac{\partial}{\partial t}}
\newcommand{\qed}{$\Box$}
\begin{document}

\title{The Penrose
Inequality in General Relativity and Volume Comparison 
Theorems involving Scalar Curvature}

\author{Hubert L. Bray}
\principaladviser{Richard M. Schoen}
\firstreader{Ben Andrews}
\secondreader{Leon Simon}
\thirdreader{Brian White}

\dept{Mathematics}
\submitdate{August 1997}
\tablespagefalse
\figurespagetrue

\beforepreface

\prefacesection{Abstract}

In this thesis we describe how 
minimal surface techniques can be used 
to prove the Penrose inequality in general
relativity for two classes of 3-manifolds.
We also describe how a 
new volume comparison theorem involving
scalar curvature for 3-manifolds follows from these same 
techniques.

The Penrose inequality in
general relativity is closely related to the positive
mass theorem, first proved by Schoen and Yau in 1979.  
In physical terms, the 
positive mass theorem states that an isolated gravitational system with 
nonnegative local energy density must have nonnegative total energy.  The
idea is that nonnegative energy densities ``add up'' 
to something nonnegative.  The
Penrose inequality, 
on the other hand, states that if an isolated
gravitational system with nonnegative local energy density contains a 
black hole of mass $m$, then the total energy of the system must be at least
$m$.  

Given a $3$-manifold $M^3$, 
we consider the function $A(V)$ equal
to the minimum area required for a surface in $M^3$ to 
contain a volume $V$.  We find that lower bounds on the 
curvature of $M^3$ yield upper bounds on $A''(V)$.  
Furthermore, in the case of an asymptotically
flat manifold which has nonnegative scalar curvature
(which is the condition needed for nonnegative energy 
density), we find that the behavior of $A(V)$ for large $V$
describes the total mass of the manifold.  In this way
we are able to use the curvature bounds of the manifold to
achieve lower bounds on the total mass.
We can also use Ricci and scalar curvature bounds on a 
compact $3$-manifold $M^3$ 
to bound the total volume of $M^3$.  
Since $A(V)$ equals zero when $V$ is either equal 
to zero or the total volume of $M^3$, 
upper bounds on $A''(V)$
force the roots of $A(V)$ to be close together, giving
an upper bound on the volume of $M^3$.

\prefacesection{Acknowledgments}

I am deeply grateful to my adviser, Professor Richard Schoen,
for suggesting the topic of this thesis and for the remarkable
insight and ideas which he routinely provides.  Rick has had a 
tremendous positive influence not only on this thesis but also 
on my education as a mathematician, and I thank him.  

I would also like to thank Professors Ben Andrews, Leon Simon, 
and Brian White who have also always been enthusiastic and 
generous with their time and who have made important
contributions to this thesis.    
I also thank the 
Rice University mathematics department for nurturing my 
interest in mathematics as an undergraduate.
In addition, I wish to acknowledge
the Department of Defense and the ARCS Foundation 
for their financial support 
during my first four years of graduate school.

I am very appreciative of the help Kevin Iga has given me
with many of the technical aspects of this thesis.  I also 
thank Kevin for many interesting mathematical conversations.
  
Finally, I especially would like to 
thank my parents and my brother for their love and 
support.  
My parents have always encouraged our interest in 
mathematics, and without their appreciation for the beauty of 
mathematics  
this thesis never would have 
happened.

\afterpreface

\chapter{Introduction}

Einstein's theory of general relativity is a theory of gravity
which asserts that matter causes the four dimensional space-time 
in which we live to be curved, and that our perception of 
gravity is a consequence
of this curvature.  Let $(N^4,\bar{g})$ 
be the space-time manifold with
metric $\bar{g}$ of signature $(-+++)$.  Then the central formula of 
general relativity is Einstein's equation,
\begin{equation}\label{Einstein}
G = 8\pi T,
\end{equation}
where $T$ is the energy-momentum tensor,  
$G=Ric(\bar{g}) - \frac12 R(\bar{g})\cdot \bar{g}$ is the Einstein curvature
tensor, $Ric(\bar{g})$ is the Ricci curvature tensor, and 
$R(\bar{g})$ is the scalar curvature of $\bar{g}$.  
The beauty of general relativity is that this simple formula 
explains gravity more accurately than Newtonian physics and is 
entirely consistent with large scale experiments.

However, the nature of the behavior of mass in general relativity
is still not well understood.  It is not even well 
understood how to define how much energy and 
momentum exists in a given region, except in special cases.  There
does exist a well defined notion of local energy and momentum 
density which is simply given by the energy-momentum tensor which,
by equation \ref{Einstein}, can be computed in terms of the 
curvature of $N^4$.  Also, if we assume that the matter
of the space-time manifold $N^4$ is concentrated in 
some central region of the universe, then
$N^4$ becomes flatter as we get farther away from this
central region.  If the curvature of $N^4$ decays quickly enough,
then $N^4$ is said to be asymptotically flat, so that with these
assumptions it is then possible to define the total mass of the 
space-time $N^4$.  Interestingly enough, though, the definition
of local energy-momentum density, which involves curvature
terms of $N^4$, bears no obvious  
resemblance to the definition of the total mass of $N^4$, which
is defined in terms of how fast the metric becomes flat at infinity.

The Penrose inequality and the positive mass theorem can both be
thought of as basic attempts at understanding the 
relationship between the local energy density of a space-time $N^4$
and the total mass of $N^4$.
In physical terms, the 
positive mass theorem states 
that an isolated gravitational system with 
nonnegative local energy density must have nonnegative total energy.  The
idea is that nonnegative energy densities must ``add up'' 
to something nonnegative.  The
Penrose inequality, on the other hand, states that if an isolated
gravitational system with nonnegative local energy density contains a 
black hole of mass $m$, then the total energy of 
the system must be at least $m$. 

Important cases of the positive mass theorem 
and the Penrose inequality can be translated
into statements about complete, asymptotically flat 
$3$-manifolds $(M^3,g)$ 
with nonnegative scalar curvature. 
If we consider $(M^3,g)$ as a space-like hypersurface of 
$(N^4,\bar{g})$ with second fundamental form $h_{ij}$ in $N^4$,
then equation \ref{Einstein} implies that
\begin{equation}
\mu = \frac{1}{16\pi} 
      [R - \sum_{i,j} h^{ij}h_{ij} + (\sum_i h_i^i)^2],
\end{equation}
\begin{equation}
J^i = \frac{1}{8\pi}\sum_j \nabla_j[h^{ij} - (\sum_k h_k^k)g^{ij}],
\end{equation}
where $R$ is the scalar curvature of the metric $g$, $\mu$ is the
local energy density, and $J^i$ is the local current density.
These two equations are called the constraint equations for $M^3$
in $N^4$, and the assumption of nonnegative energy density everywhere
in $N^4$ implies that we must have
\begin{equation}
\mu \ge \left(\sum_i J^i J_i \right)^\frac12
\end{equation}
at all points on $M^3$ \cite{SY4}.  Thus we see that 
if we restrict our attention 
to $3$-manifolds which have zero mean curvature in $N^3$, 
the constraint equations and the assumption of nonnegative
energy density imply that $(M^3,g)$ has nonnegative scalar 
curvature everywhere.  We also assume that $(M^3,g)$ is asymptotically
flat, which is defined in section \ref{s6}, in which case 
we can define the total mass of $M^3$, also given
in section \ref{s6}.  

An ``end'' of an $n$-manifold 
is a region of the manifold diffeomorphic to $\real^n - B_1(0)$
where $B_1(0)$ is the ball of radius one in $\real^n$.  
In general, $M^3$ may have any number of disjoint ends, 
but for 
simplicity, let us assume that $M^3$ has only one disjoint end 
and that it is
asymptotically flat.
In section \ref{s6} we will show that without loss of 
generality (for stating the Penrose inequality and the positive mass theorem)
we may assume that $(M^3 - K,g)$ is isometric to 
$({\real}^3 - B,h)$ for some compact set $K$  
in $M^3$ and some ball $B$ 
in $\real^3$ centered around the origin, 
and for some constant $m$, where  
$h_{ij} = (1 + \frac{m}{2r})^4 \delta_{ij}$ and 
$r$ is the radial coordinate in $\real^3$.
This is a convenient 
assumption about $M^3$, because the total mass 
of the system is then just $m$.
The metric $(\real^3-\{0\},h)$ has zero scalar curvature,
is spherically symmetric, and 
is called the Schwarzschild metric of mass $m$, 
and we say that in the above case, $M^3$ is Schwarzschild with mass $m$
at infinity.  Using this simplified setup, we can make a statement which 
is equivalent to the positive mass theorem in this setting.

\vspace{.1in}\noindent
{\bf The Positive Mass Theorem (Schoen, Yau, 1979) }{\it
Suppose $(M^3,g)$ is complete, has nonnegative scalar curvature, 
and is Schwarzschild with mass $m$ at infinity.
Then $m \ge 0$, and $m=0$ if and only if
$(M^3,g)$ is isometric to $\real^3$ with the standard flat metric.}
\vspace{.1in}

Apparent horizons of black holes in $N^4$ 
correspond to outermost minimal spheres 
of $M^3$ if we assume $M^3$ has zero second fundamental form
in $N^4$.  An 
outermost minimal sphere is a sphere in $M^3$ which locally minimizes area 
(and hence has zero mean curvature) and which is not contained entirely inside
another minimal sphere.  We will also use the term horizon to 
mean an outermost minimal sphere in $M^3$.  
It is easy to show that two outermost horizons
never intersect.  Also, it follows from a stability argument that
these minimal surfaces are
always spheres \cite{HE}.  
However, there may be
more than one outermost minimal sphere, with each minimal sphere corresponding
to a different black hole.  As we will see in the next section, 
there is a 
strong motivation to define the mass of a black hole as 
$\sqrt{\frac{A}{16\pi}}$, where $A$ is the surface area of 
the horizon.  
Hence, the physical statement that a system
with nonnegative energy density containing 
a black hole of mass $m$ must 
have total mass
at least $m$ can be translated into the following geometric statement.

\vspace{.1in}\noindent
{\bf The Penrose Inequality (Huisken, Ilmanen, announced 1997) }{\it
Suppose $(M^3,g)$ is complete, 
has nonnegative scalar curvature, contains an outermost minimal sphere
with surface area $A$, and is Schwarzschild with mass $m$
at infinity.  Then $m \ge \sqrt{\frac{A}{16\pi}}$, with equality only in the case
that $(M^3,g)$ is isometric to the Schwarzschild metric of mass $m$ outside the horizon.
}\vspace{.1in}

When $m > 0$, the Schwarzschild metrics $(\real^3,h)$
have minimal spheres at $r = \frac{m}{2}$ with areas 
$16\pi m^2$ so that
these metrics give equality in the Penrose inequality, 
and in fact, according to the recent announcement of Huisken 
and Ilmanen, 
these are the only metrics (outside the horizon)
which give equality.  

The proof that Huisken and Ilmanen used to prove the Penrose
inequality is as interesting as the theorem itself.  We discuss
the main ideas of their proof in section \ref{Hawking}.  We also
give another proof of the Penrose inequality for two classes
of manifolds in chapter \ref{Penrose} using isoperimetric 
surface techniques.  Both approaches are also interesting
because they give hints about the nature of 
quasi-local mass in general relativity.

Also, using isoperimetric surface techniques we are able to 
prove a generalized Penrose inequality for a class of 
manifolds in the case that $(M^3,g)$ has more than one horizon.
The idea is that if $(M^3,g)$ has more than one horizon, then
it should be possible to bound the total mass from below by some
function of the areas of the horizons.  In this way we hope to
understand how masses ``add'' in general relativity.
We state the conditions under which we can prove a 
generalized Penrose inequality in the introduction
to chapter \ref{Penrose} and then conjecture that this generalized
Penrose inequality is always true.

\section{Motivation behind the Penrose Inequality}

In 1973, Roger Penrose proposed the Penrose inequality as a 
test of the cosmic censor hypothesis \cite{P}.  The cosmic 
censor hypothesis states that naked singularities do not 
develop starting with physically reasonable nonsingular 
generic initial
conditions for the Cauchy problem in general relativity.  (However,
it has been shown by Christodoulou \cite{C} 
that naked singularities
can develop from nongeneric initial conditions.)  
If naked singularities did typically develop from generic 
initial conditions, then this would be a serious problem for
general relativity since it would not be possible to solve 
the Einstein equations uniquely past these singularities.
Singularities such as black holes do develop but are shielded
from observers at infinity by their horizons so that the Einstein
equations can still be solved from the point of view of an 
observer at infinity.

A summary of Penrose's argument can be found in \cite{JW}.
The main idea is to consider a space-time $(N^4,\bar{g})$ 
with given initial
conditions for the Cauchy problem $(M^3,g)$ with zero second 
fundamental form in $N^4$.  We assume that $N^4$ has 
nonnegative energy density everywhere, so by the constraint
equations $M^3$ must have nonnegative scalar curvature.  
Suppose also $(M^3,g)$ has an
outermost 
apparent horizon of area $A$, and event horizon of area $A_i$, 
and total mass $m_i$ (see \cite{HE}, \cite{HP} for the definitions of 
these horizons).  As long as
a singularity does not form, then it is assumed that 
eventually the space-time 
should converge on some stationary final state.  From the 
theorems of Israel \cite{IJW}, Hawking \cite{H2}, and 
Robinson \cite{RJW}, the only stationary vacuum black holes
are the Kerr solutions which satisfy
\begin{equation}\label{Kerr}
A_{f} = 8\pi[m_f^2+ (m_f^4 - J^2)^\frac12] \le 16\pi m_f^2,
\end{equation}
where $A_f$ is the area of the horizon of the Kerr black hole,
$m_f$ is the mass at infinity, and $J$ is the angular momentum.

However, by the Hawing area theorem \cite{H2JW}, the area of 
the event horizon of the black hole is nondecreasing.  Thus,
$A_f \ge A_i$.  Also, presumably some energy radiates off 
to infinity, so we expect to have $m_i \ge m_f$.  

The apparent horizon is defined to be the outer boundary of the 
region in $M^3$ which contains trapped or marginally trapped 
surfaces \cite{HE}.  
The apparent horizon itself must then be a marginally
trapped surface, and hence satisfies
\begin{equation}
H + h^{ij}(g_{ij} - r_i r_j) = 0
\end{equation}
where $H$ is the mean curvature of the apparent horizon in 
$M^3$, $h$ is the second fundamental form of $(M^3,g)$ in 
$(N^4,\bar{g})$, and $r$ is the outward unit normal to the 
apparent horizon in $M^3$.  Hence, since we chose $M^3$ to 
have zero second fundamental form,  $h^{ij} = 0$, so that 
the apparent horizon is a 
zero mean curvature surface in $M^3$.  Furthermore, if we 
consider the surface of smallest area which encloses the apparent
horizon, it too must have zero mean curvature and hence is a 
marginally trapped surface in $M^3$.  Thus, the apparent horizon
is an outermost minimal surface of $M^3$, which by stability
arguments, must be a sphere \cite{HE}.  
Since the event horizon
always contains the apparent horizon, $A_i \ge A$, so putting
all the inequalities together we conclude that
\begin{equation}
m_i \ge m_f \ge \sqrt{\frac{A_f}{16\pi}} \ge 
\sqrt{\frac{A_i}{16\pi}} \ge \sqrt{\frac{A}{16\pi}} 
\end{equation}

Thus, Penrose argued, assuming the cosmic censor hypothesis 
and a few reasonable sounding 
assumptions as to the nature of gravitational
collapse, given a complete 
asymptotically flat $3$-manifold $M^3$ of total mass $m_i$ with
nonnegative scalar curvature which has an outermost minimal sphere
of total area $A$, then 
\begin{equation}
m_i \ge \sqrt{\frac{A}{16\pi}} 
\end{equation}
Conversely, he argued, if one could find an $M^3$ which was a 
counterexample to the above inequality, then it would be likely
that the counterexample, when used as initial conditions in
the Cauchy problem for Einstein's equation, would produce a 
naked singularity.  Since Huisken and Ilmanen have proved the 
above inequality, they
have ruled out one possible way of constructing counterexamples to
the cosmic censor hypothesis.

\section{The Schwarzschild Metric}

\begin{figure}
\vspace{3.5in}
\includegraphics{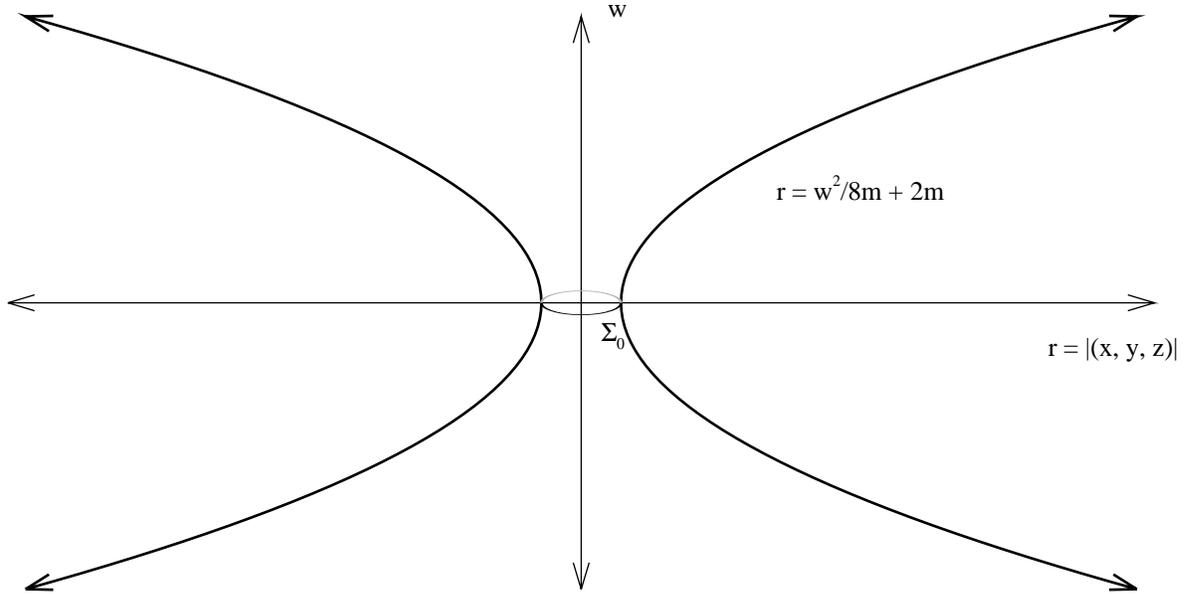}
\caption{The space-like Schwarzschild metric of mass $m$,
$(\real^3-\{0\}, h)$, 
viewed as a submanifold
of four-dimensional Euclidean space.\label{Schwarz}}
\end{figure}

The space-like 
Schwarzschild metric is a particularly important example
to consider when discussing the Penrose inequality.  First
of all, Huisken and Ilmanen prove that 
it is the only $3$-manifold which gives equality
in the Penrose inequality.  Also, if a $3$-manifold is assumed
to be complete, spherically symmetric, and have zero scalar
curvature, then it must be isometric to 
either a Schwarzschild metric of 
mass $m > 0$ or $\real^3$, which can be viewed as the 
Schwarzschild metric when $m=0$.   

In addition, understanding the Schwarzschild metric is 
particularly important for chapter \ref{Penrose} 
because we show in section \ref{s6} that without loss of
generality for proving the Penrose inequality for $M^3$ we may 
assume that outside a compact set $M^3$ is spherically 
symmetric with zero scalar curvature, which means that in 
this region it is isometric to the Schwarzschild metric of
some mass $m$.  
We also show in that same section that all asymptotically flat
metrics of nonnegative scalar curvature can be perturbed
pointwise less than $\epsilon$ in such a way that the total mass
is changed less than $\epsilon$ too and so that the new metric is 
isometric to the Schwarzschild metric outside a compact set.
Thus, the Schwarzschild metric is a useful picture to keep in
mind.

The space-like Schwarzschild metric, $(\real^3-\{0\},h)$, is a
time symmetric asymptotically flat 
three-dimensional maximal slice (chosen to have zero momentum at
infinity) of the four-dimensional Schwarzschild space-time metric.The
space-like Schwarzschild metric is conformal to $\real^3-\{0\}$ 
with
$h_{ij}=\left(1+\frac{m}{2r}\right)^4\delta_{ij}$.
The Schwarzschild metric of mass $m$, $(\real^3-\{0\},h)$, can also be
isometrically embedded into four-dimensional Euclidean space as the
three-dimensional set of points in $\real^4=\{(x,y,z,w)\}$ satisfying
$|(x,y,z)|=\frac{w^2}{8m}+2m$, seen in figure \ref{Schwarz}.  
Hence, $\Sigma_0$ is a
minimal sphere of area $16\pi m^2$, so we have equality in the Penrose
inequality.

\section{The Spherically Symmetric Case}\label{sscase}

In this section we sketch a proof of the Penrose inequality 
in the case that $M^3$ is spherically symmetric.  The proof 
is very easy conceptually, but what is more important is that
some of the ideas generalize.  In particular, we will see why
the minimal sphere in the Penrose inequality must be outermost.

Let $(M^3,g)$ be a complete 
asymptotically flat spherically symmetric 
$3$-manifold with nonnegative scalar curvature.  
For convenience, we also assume that $(M^3,g)$ is isometric to
the Schwarzschild metric of some mass $m$ outside a large
compact set.  Then the total mass of $(M^3,g)$ is $m$. 
Let $\Sigma(V)$ 
be the spherically symmetric sphere containing a volume $V$
in $M^3$.  Let $A(V)$ be the area of this sphere.  It turns out
that the function $A(V)$, $V \ge 0$, captures all the information
about $M^3$ since $M^3$ is spherically symmetric.

Let $R(V)$ be the scalar curvature of $M^3$ on $\Sigma(V)$.
From the calculations we will do in section \ref{s3}, it 
follows that 
\begin{equation}
R(V) = \frac{8\pi}{A} - 2A(V)A''(V) - \frac32 A'(V)^2
\end{equation}
Define 
\begin{equation}
m(V) = \left( \frac{A(V)}{16\pi}\right)^\frac12
\left(1 - \frac{1}{16\pi} A(V)A'(V)^2 \right)
\end{equation}
It turns out that $m'(V) \ge 0$ whenever $A'(V) \ge 0$ since
we find that 
\begin{equation}\label{mprime}
m'(V) = \frac{A'(V)}{16\pi} \left(\frac{A(V)}{16\pi}\right)
^\frac12 R(V)
\end{equation}
and $R(V) \ge 0$.

Let $\Sigma(V_0)$ be the outermost minimal sphere.  It 
follows that $A'(V) \ge 0$ for all $V \ge V_0$.  Hence, $m(V)$
is increasing in this range as well, so 
\begin{equation}\label{lim}
\lim_{V \rightarrow \infty} m(V) \ge m(V_0)
\end{equation}
Furthermore $m(V_0) = \sqrt{\frac{A(V_0)}{16\pi}}$ since 
$A'(V_0) = 0$.  Also, we assumed that $M^3$ was isometric to
the Schwarzschild metric outside a large compact set, and we 
claim that $m(V) = m$, the mass parameter of the Schwarzschild
metric, in this region, or equivalently, for $V > V_{LARGE}$ for
some $V_{LARGE}> 0$.  To see this, consider the mass function
$m(V)$ defined on the Schwarzschild metric, where now $V$ refers
to the volume contained by the spherically symmetric spheres 
of the Schwarzschild metric which is outside the horizon.  Then
by equation \ref{mprime}, $m(V)$ is constant for all $V$ on the
Schwarzschild metric since the Schwarzschild metric has zero 
scalar curvature.  Furthermore, setting $V=0$ and considering 
$m(V)$ at the horizon yields $m(0) = \sqrt{\frac{A(0)}
{16\pi}} = m$, the mass parameter of the Schwarzschild metric,
since the Schwarzschild metric gives equality in the Penrose 
inequality.  Thus, $m(V) = m$ for all $V$ in the Schwarzschild
metric, so going back to $(M^3,g)$, we see that $m(V) = m$, 
the mass parameter of the Schwarzschild metric, for 
$V > V_{LARGE}$.  Thus, it follows from inequality \ref{lim}
that   
\begin{equation}
m \ge \sqrt{\frac{A(V_0)}{16\pi}}
\end{equation}
which proves the Penrose inequality for spherically
symmetric manifolds.

\begin{figure}
\vspace{6in}
\includegraphics{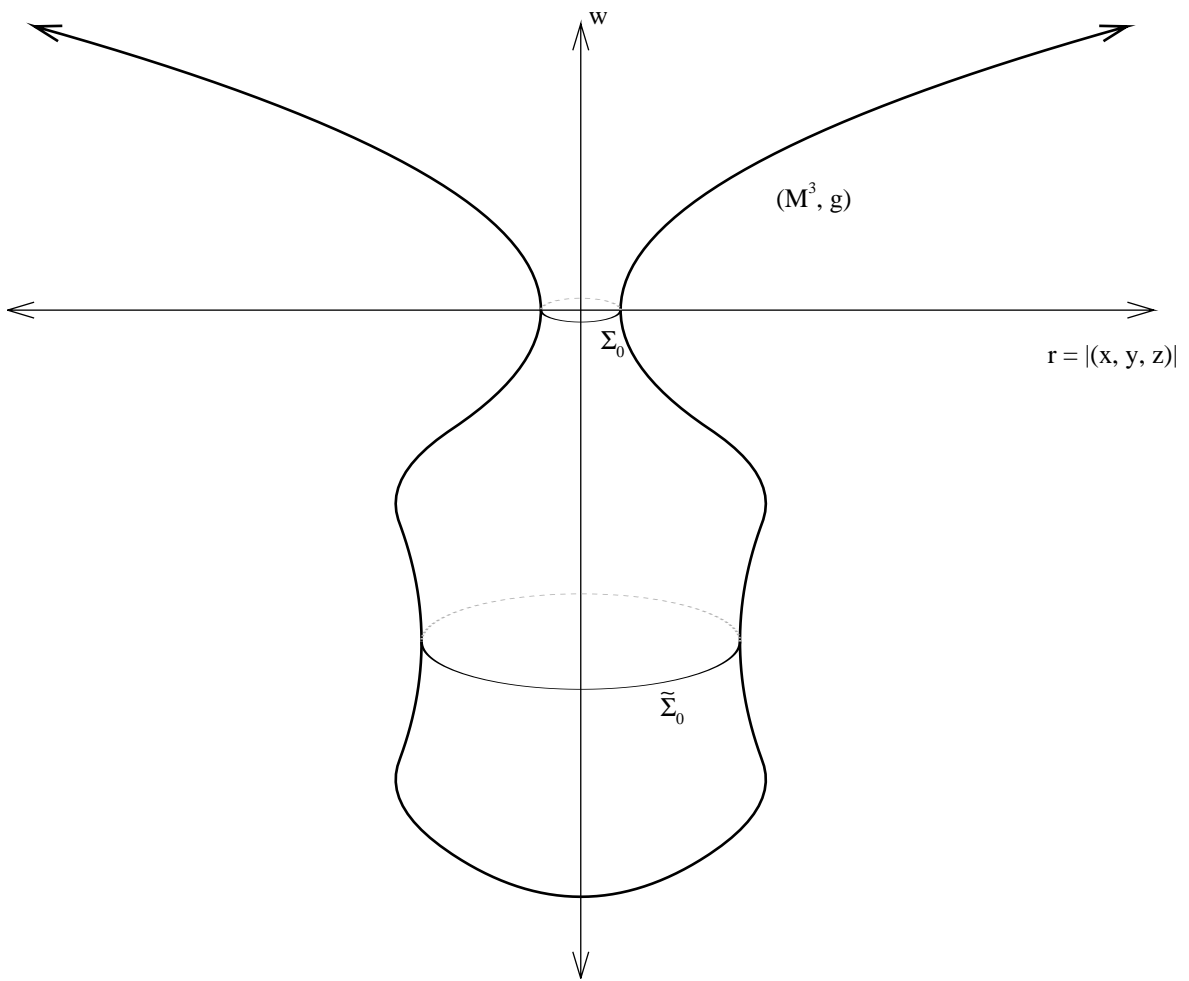}
\caption{Counterexample to Penrose inequality if the 
minimal sphere is not outermost.\label{outermost}}
\end{figure}

Conversely, we see that equation \ref{mprime} can be used to
construct spherically symmetric manifolds which do not satisfy
the Penrose inequality if we do not require the minimal sphere
to be outermost.  
In figure \ref{outermost}, 
we are viewing $(M^3,g)$ as an isometrically embedded 
submanifold of $\real^4$
with the standard Euclidean metric.  $(M^3,g)$ is spherically symmetric and
is constructed by rotating the curve shown above around the $w$-axis in
$\real^4$.  Hence, $\Sigma_0$ and $\tilde{\Sigma}_0$ are both
$2$-spheres, and we can choose the curve shown above so that the scalar
curvature of $(M^3,g)$, $R(g)$, is non-negative.

The Penrose inequality, $m\ge \sqrt{\frac{|\Sigma_0|}{16\pi}}$, is
true for $\Sigma_0$, but is not true for $\tilde{\Sigma}_0$.  However,
$\tilde{\Sigma}_0$ is not an outermost minimal sphere since it is
contained by another minimal sphere, namely, $\Sigma_0$.  In fact,
since $\tilde{\Sigma}_0$ is not outermost, we can construct a
spherically symmetric manifold like the one shown above so that the
area of $\tilde{\Sigma}_0$ is as large as we like and the total
mass of
$(M^3,g)$ is still one.

\section{The Hawking Mass and Inverse Mean Curvature Flows}
\label{Hawking}

One goal in general relativity is to  
understand how to 
define the amount of mass inside a given region.   
In section \ref{sscase}, we defined a function $m(V)$ which was 
increasing as a function of $V$ outside the outermost minimal
sphere.  Furthermore, for large $V$, $m(V)$ equaled the total
mass of the manifold.  Hence, it seems reasonable to say that
the spherically symmetric sphere $\Sigma(V)$ defined in 
section \ref{sscase} contains a mass $m(V)$.  The function $m(V)$
is called a quasi-local mass function.

Naturally we would like to define a quasi-local mass
function which would measure the amount of mass inside any 
surface $\Sigma$ which is the boundary of a region in any 
$3$-manifold $M^3$. 
We refer the reader to \cite{CY}, \cite{Ba4} ,\cite{G}, 
and \cite {H1} for a complete
discussion of this topic.  

In this section we discuss a definition
of quasi-local mass proposed by Hawking called the Hawking mass.
Going back to section \ref{sscase}, we recall that in the 
spherically symmetric case,
\begin{equation}\label{ssmv}
m(V) = \left( \frac{A(V)}{16\pi}\right)^\frac12
\left(1 - \frac{1}{16\pi} A(V)A'(V)^2 \right)
\end{equation}
As can be seen from the calculations in the appendix  
and as will be shown later, it so happens that $A'(V) = H(V)$,
where $H(V)$ is the mean curvature of the spherically symmetric
sphere $\Sigma(V)$.  Hence, one way to generalize equation
\ref{ssmv} is to define
\begin{equation}
m(\Sigma) = \left(\frac{A}{16\pi} \right)^\frac12 
\left(1 - \frac{1}{16\pi}\int_\Sigma H^2  \right)
\end{equation}
where $A$ is the area of $\Sigma$ and $H$ is the mean curvature
of $\Sigma$ in $M^3$.  It turns out that this definition of
quasi-local mass has some very important properties.

In \cite{G}, Geroch showed that if $M^3$ has nonnegative
scalar curvature, then the Hawking mass is nondecreasing
when the 
surface $\Sigma$ is flowed out at a rate equal to the inverse
of its mean curvature.  This is straight forward to check using 
equations \ref{eqn:darea} and \ref{eqn:dcurve} from the 
appendix and the Gauss equation which is given
in equation \ref{Gauss} of section \ref{s3}.  
In view of this, Jang and Wald proposed using the
Hawking mass function to prove the Penrose inequality
\cite{JW}.  They
suggested that we should let $\Sigma(0)$ be an outermost horizon 
and then to flow out using an inverse mean curvature flow to
create a family of surfaces $\Sigma(t)$ flowing out to infinity.
Since the Hawking mass function $m(t) = m(\Sigma(t))$ 
is nondecreasing as a function of $t$, we have 
\begin{equation}\label{limt}
\lim_{t \rightarrow \infty} m(t) \ge m(0)
\end{equation}
Furthermore $m(0) = \sqrt{\frac{A}{16\pi}}$, where $A$ is the 
area of the outermost horizon $\Sigma(0)$, since horizons have
zero mean curvature.  Hence, Jang and Wald proposed a proof
of the Penrose inequality which is basically a generalization
of the proof which works in the spherically symmetric case.

The main problem for this type of proof is the existence of
an inverse mean curvature flow.  Naturally, if the mean curvature
of the surface ever went to zero or became negative, the flow
could not exist, at least in this form.  However, Huisken
and Ilmanen have recently announced that they have been 
able to generalize the idea of an inverse mean curvature flow
to a ``weak'' inverse mean curvature flow which always exists
and hence can be used to prove the Penrose inequality \cite{HI}.
  
They introduce the notion of a ``maximal minimal mean convex
hull'' of a surface $\Sigma$ 
which equals the outermost surface of minimum area
needed to enclose $\Sigma$.  Then their weak inverse mean 
curvature flow can be thought of as continuously replacing
$\Sigma(t)$ with the maximal minimal mean convex hull of 
$\Sigma(t)$
while flowing out
using the inverse mean curvature flow.  The replacement process
can then be shown to never decrease the mass and also to keep the
mean curvature of $\Sigma(t)$ nonnegative.  The resulting weak
flow is a family of surfaces which occasionally has ``jumps''
and for which the Hawking mass is nondecreasing.  They use this 
to prove
the Penrose inequality as stated in the introduction to this
chapter.  

Partial results on the Penrose inequality have also been
found by Herzlich \cite {He} 
using the Dirac operator which Witten \cite{Wi}
used to prove the positive mass theorem, by 
Gibbons \cite{Gi} in the special case of collapsing shells, by
Tod \cite{T}, and by Bartnik \cite{Ba3} for quasi-spherical 
metrics.

However, other versions of the Penrose inequality remain
open.  As in the positive mass theorem \cite{SY5}, 
the assumption of nonnegative scalar curvature for $M^3$ should
be able to be modified to include a more general local nonnegative
energy condition.
It is also natural to ask what kind of generalized
Penrose inequality we should expect for manifolds 
with multiple horizons.  
In chapter \ref{Penrose} 
we prove a generalized Penrose inequality for a certain class
of manifolds and in section \ref{conjectures} 
conjecture that this generalized
inequality is always true.

\section{Volume Comparison Theorems}

The isoperimetric surface techniques which we will develop to 
study the Penrose inequality in general relativity also can be
used to prove several volume comparison theorems, including a 
new proof of Bishop's volume comparison theorem for positive 
Ricci curvature. 

Let $(S^n,g_0)$ be the standard metric (with any scaling) on $S^n$
with constant Ricci curvature $Ric_0 \cdot g_0$.  Bishop's theorem says that if
$(M^n,g)$ is a complete Riemannian manifold ($n\ge 2$) with $Ric(g) \ge Ric_0 \cdot g$,
then $\mbox{Vol}(M^n) \le \mbox{Vol}(S^n)$.  It is then natural to ask
whether a similar type of volume comparison theorem could be true for
scalar curvature.  
As it happens, a lower bound on scalar curvature by itself is not 
sufficient to give an upper bound on the total volume.  We can scale
a cylinder, $S^2 \times \real$, to have any positive scalar curvature and still have
infinite volume.  

Since a lower bound on scalar curvature is not enough to 
realize an upper bound on the volume of a manifold, 
in chapter \ref{volume} we consider $3$-manifolds $(M^3,g)$
which satisfy $R(g) \ge R_0$ and 
$Ric(g) \ge \epsilon \cdot Ric_0 \cdot g$, where 
$(S^3,g_0)$ is the standard metric on $S^3$ with constant
scalar curvature $R_0$ and constant Ricci curvature $Ric_0
\cdot g_0$ (so that naturally $R_0 = 3 Ric_0$).   
It turns out that there do exist values of 
$\epsilon < 1$ for which 
these curvature conditions imply that   
$\mbox{Vol}(M^3) \le \mbox{Vol}(S^3)$, giving us a volume 
comparison theorem for scalar curvature. 
We also find the best value for $\epsilon$ for which this
theorem is true.


\chapter{The Penrose Inequality}\label{Penrose}

We will use isoperimetric surfaces to prove the Penrose inequality for two cases.
In the first case, we will use surfaces which globally minimize area among surfaces which
contain the same volume, and in the second case, we will need to look at collections of 
surfaces, each of which locally minimizes area among surfaces enclosing the same volume.
In both cases these surfaces will have constant mean curvature.

\begin{definition}\label{defav}

Suppose $M^3$ is asymptotically flat, complete, and 
has only one outermost minimal sphere $\Sigma_0$.  Let
$\tilde{M}^3$ be the closure of the
component of $M^3 - \Sigma_0$
that contains the asymptotically flat end.  We define
\[A(V) = \inf_\Sigma \{\mbox{Area}(\Sigma)\,\,| \,\, 
\Sigma \mbox{ contains a volume $V$ outside } \Sigma_0    \}\]
where $\Sigma$ is the boundary of some 3-dimensional region in $M^3$
and $\Sigma$ is a surface in $\tilde{M}^3$ in the same
homology class of $\tilde{M}^3$ as the horizon $\Sigma_0$.

\end{definition}
  
If $\Sigma$ contains a volume $V$ outside the horizon and 
$\mbox{Area}(\Sigma) = A(V)$, then we say that $\Sigma$
minimizes area with the given volume constraint.  Naturally, $\Sigma$ 
could have multiple components, as long as one of the components contains
the horizon.   

\begin{condition}\label{R3}
$(M^3,g)$ has only one horizon $\Sigma_0$, and 
for each $V>0$, 
if one or more area minimizers exist for $V$, then at least
one of these area minimizers for the volume $V$ has exactly one component.
\end{condition}
Note that condition \ref{R3} does not assume that an area minimizer
exists for each $V$, just that if one or more do exist for $V$,
at least one of these minimizers for $V$ has only one component.  However,
once condition \ref{R3} is assumed, the existence of an area minimizer
$\Sigma(V)$ for each $V \ge 0$ follows from the behavior of $A(V)$ as will 
be shown in section \ref{s7}.  Also, assuming condition \ref{R3}, we can
prove the Penrose inequality.

\begin{theorem}\label{R4}
Suppose $(M^3,g)$ is complete, 
has nonnegative scalar curvature, contains an outermost
minimal sphere with surface area $A$, is Schwarzschild with mass $m$
at infinity, and satisfies condition \ref{R3}.  
Then $m \ge \sqrt{\frac{A}{16\pi}}$.
\end{theorem}

Naturally, we want to find a way to get around condition \ref{R3}.  This
can be partially accomplished if we 
consider a different minimization problem, minimizing 
the sum of the areas to the three halves power given a volume constraint. 
Posing this type of problem seems strange at first, but turns out to be 
surprisingly natural for manifolds with nonnegative scalar curvature.
 
\begin{definition}\label{deffv}
Suppose $M^3$ is asymptotically flat, complete, and 
has any number of horizons.  Let $\tilde{M}^3$ be
the closure of 
the component of $M^3 - \{\mbox{the horizons}\}$ that contains the 
asymptotically flat end.  Let
\[ F(V) = \inf_{\{\Sigma_i\}} \{ \sum_{i} \mbox{Area}(\Sigma_i)^\frac32 
\,\,|\,\,
\{\Sigma_i\} \mbox{ contain a volume $V$ outside the horizons} \} \]
where the $\{\Sigma_i\}$ are the boundaries of the components of some 
3-dimensional open region in $M^3$ and 
$\bigcup_i \Sigma_i$ is in $\tilde{M}^3$ and is in the
homology class of $\tilde{M}^3$ which contains both a large sphere at infinity
and the union of the horizons.  

\end{definition}

If the collection
$\{\Sigma_i\}$ contains a volume $V$ outside the horizons and 
$\sum_{i} \mbox{Area}(\Sigma_i)^\frac32 = F(V)$, then we say that 
$\{\Sigma_i\}$ minimizes $F$ for the volume $V$.  The only problem that 
occurs with this optimization problem is that two or more surfaces 
$\Sigma_i$ and $\Sigma_j$ can push up against each other.

\begin{condition}\label{R5}
For each $V>0$, if one or more sets of surfaces minimize $F$ for the 
volume $V$, then at least one of these sets of surfaces 
$\{\Sigma_i\}$ which minimize $F$
for the volume $V$ is pairwise disjoint, that is, $\Sigma_i \cap
\Sigma_j = \emptyset$ for all $i \ne j$.
\end{condition}
Note that condition \ref{R5} does not assume that an $F$ minimizer exists
for each $V$, just that if one or more do exist for $V$, at least one of 
these $F$ minimizers for $V$ does not have any of its surfaces pushing up 
against or touching each other.  However, once condition \ref{R5} is 
assumed, the existence of an $F$ minimizer for each $V \ge 0$ follows in a 
nice way from the
behavior of $F(V)$ as we will see in section \ref{s8}.  
Also, it is possible to verify experimentally using a computer
to construct axially symmetric, conformally flat metrics with multiple
horizons that there are examples of 3-manifolds which appear to 
satisfy condition \ref{R5} but not condition \ref{R3}. 

At first glance the issue of existence for this optimization problem
looks bleak for several reasons.  First, two components $\Sigma_i$ and 
$\Sigma_j$ can be joined into one component by a thread of area zero.  This
is always disadvantageous for minimizing $F$ and so is not a 
problem.  Also, one component of $\{\Sigma_i\}$ might run off to infinity.
In addition, ``bubbling'' might occur, where the optimal configuration is
an infinite number of tiny balls with a finite total volume.  Amazingly,
if we assume condition \ref{R5},
these last two problems do not occur and the Penrose inequality follows.

\begin{theorem}\label{R6}
Suppose $(M^3,g)$ is complete, 
has nonnegative scalar curvature, contains one or more outermost
minimal spheres with surface areas $\{A_i\}$, 
is Schwarzschild with mass $m$
at infinity, and satisfies condition \ref{R5}.  
Then $m \ge \left(\sum_{i=1}^n \left(\frac{A_i}{16\pi}\right)^\frac32 \right)^\frac13$.
\end{theorem}

Thus, condition 2 implies a stronger version of the Penrose
inequality since \\
$\left(\sum_{i=1}^n \left(\frac{A_i}
{16\pi}\right)^\frac32 \right)^\frac13 
\ge \sqrt{\frac{A_j}{16\pi}}$
for all $j$.
Based on this, it seems plausible to conjecture that for multiple
black holes,
this stronger Penrose inequality is always true.  We provide
additional motivation for this conjecture in section
\ref{conjectures}.

\begin{conjecture}\label{genpenrose}
Suppose $(M^3,g)$ is complete, 
has nonnegative scalar curvature, contains one of more outermost
minimal spheres with surface areas $\{A_i\}$, and 
is \\
Schwarzschild with mass $m$
at infinity.
Then $m \ge \left(\sum_{i=1}^n \left(\frac{A_i}{16\pi}\right)^\frac32 \right)
^\frac13$.
\end{conjecture}

\section{Isoperimetric Surface Techniques}   
\label{s3}

We will prove theorems \ref{R4} and \ref{R6} using constant mean curvature 
surfaces which minimize area given a volume constraint.  
First, let us assume the hypotheses of theorem \ref{R4} including 
condition \ref{R3} 
and recall the definition of $A(V)$ given in the previous section.  
Under these circumstances, we show in section \ref{s7} that for all $V \ge 0$,
there exists a smooth, constant mean curvature  
surface $\Sigma(V)$ which
minimizes area among surfaces which enclose a volume $V$ outside the horizon.
By condition \ref{R3}, we may choose the minimizer $\Sigma(V)$ 
to have only one component, 
and since $\Sigma(V)$
is a minimizer, the area of $\Sigma(V)$ is $A(V)$.

$A(V)$ contains important geometric information,
including the fact that $A(0)$ is the area of the horizon.  The fact 
that the horizon is outermost implies that $A(V)$ is nondecreasing.

We must use this last fact somewhere, because  
the Penrose inequality is definitely not true without the assumption
that the minimal sphere in the conjecture is outermost.
In fact, it is easy to
construct a complete, spherically symmetric 3-manifold with nonnegative scalar
curvature and total mass $1$ and an 
arbitrarily large (non-outermost) minimal sphere. 
It is worth noting that when considering other possible approaches to
the Penrose inequality, the hypothesis of the Penrose inequality
that the minimal sphere is outermost is often one of the more delicate 
and difficult points to handle.

Also, as we will prove in section \ref{s6.5}, the total mass $m$ of $(M^3,g)$,
is encoded in the asymptotic behavior of the function $A(V)$ for large $V$. 
Hence, the key to proving theorem \ref{R4} is understanding how the assumption
of nonnegative scalar curvature on $(M^3,g)$ bounds the behavior of $A(V)$.

\begin{theorem} \label{APP}
Suppose $(M^3,g)$ is complete, 
has nonnegative scalar curvature, is \\
Schwarzschild
at infinity, and satisfies condition \ref{R3}.  Then 
the function $A(V)$ defined in definition \ref{defav} satisfies
\[A''(V) \le \frac{4\pi}{A(V)^2} - \frac{3 A'(V)^2}{4 A(V)} \]
in the sense of comparison functions,
where this means that for all $V_0 \ge 0$ there exists a
smooth function $A_{V_0}(V)\ge A(V)$ with $A_{V_0}(V_0)=A(V_0)$
satisfying
\[A_{V_0}''(V_0) \le \frac{4\pi}{A_{V_0}(V_0)^2} - 
                     \frac{3 A_{V_0}'(V_0)^2 }{4 A_{V_0}(V_0)}\]
\end{theorem}
{\em Proof.} 
First we comment that all the inequalities which we state 
``in the sense of comparison function'' are also true 
distributionally.  We do not need this, so we do not prove it, 
but the proof is very similar to the proof of lemma 
\ref{AJ} in section \ref{s4}.  

To get an upper bound for $A''(V)$ at $V=V_0$, we will do
a unit normal variation on $\Sigma(V_0)$.  Let
$\Sigma_{V_0}(t)$ be the surface created by flowing $\Sigma(V_0)$
out at every point in the normal direction at unit speed for time
$t$.  Since $\Sigma(V_0)$ is smooth, we can do this variation for
$t \in (-\delta,\delta)$ for some $\delta > 0$.  Abusing notation
slightly, we can also parameterize these surfaces by their volumes as
$\Sigma_{V_0}(V)$ so that $V = V_0$ corresponds to
$t = 0$.  Let $A_{V_0}(V) = \mbox{Area}(\Sigma_{V_0}(V))$.  Then 
$A(V_0) = A_{V_0}(V_0)$ and $A(V) \le A_{V_0}(V)$ since
$\Sigma_{V_0}(V)$ is not necessarily minimizing for its volume. 
Hence,

\[A''(V_0) \le A_{V_0}''(V_0).\]

\begin{figure}
\vspace{3in}
\includegraphics{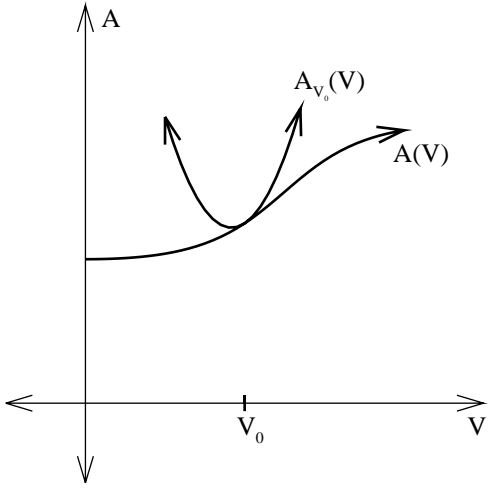}
\caption{Graphical demonstration that 
$A''(V_0)\le A''_{V_0}(V_0)$.\label{AofVgraph}}
\end{figure}

To compute $A''_{V_0}(V_0)$, we will need to compute the
first and second derivatives of the area of $\Sigma_{V_0}(t)$ and the volume
that it encloses.  We will use the formulas
\[\dot{d\mu}=H\,d\mu\,\,\,\,\mbox{ and }\,\,\,\,
\dot{H}=-||\Pi||^2-Ric(\nu,\nu)\]
where the dot represents differentiation with respect to $t$, $d\mu$ is
the surface area 2-form for $\Sigma_{V_0}(t)$, $\Pi$ is the second fundamental
form for $\Sigma_{V_0}(t)$, $H=\mbox{trace}(\Pi)$ is the mean curvature, and
$\nu$ is the outward pointing unit normal vector.  Since $A_{V_0}(t)=
\int_{\Sigma_{V_0}(t)}d\mu$,
\[A_{V_0}'(t)=\int_{\Sigma_{V_0}(t)}H\,d\mu\]
And since $V'(t)=\int_{\Sigma_{V_0}(t)}d\mu=A_{V_0}(t)$, 
we have that at $t = 0$,
\[A_{V_0}'(V)=A_{V_0}'(t)/V'(t)=H\]
By single variable calculus,
\[A_{V_0}''(V) = \frac{A_{V_0}''(t) - A_{V_0}'(V) V''(t)}{V'(t)^2}\]
so that at $t = 0$,
\begin{eqnarray*} 
A_{V_0}(V_0)^2 A_{V_0}''(V_0) &=& A_{V_0}''(t) - H V''(t) \\
                  &=& \frac{d}{dt} \int_{\Sigma_{V_0}(t)} H\, d\mu \,\, - \,\,
                      H \frac{d}{dt}\int_{\Sigma_{V_0}(t)} d\mu  \\
                  &=& \int_{\Sigma(V_0)} \dot{H}\, d\mu \\
                  &=& \int_{\Sigma(V_0)} -||\Pi||^2 - Ric(\nu,\nu) 
\end{eqnarray*}
By the Gauss equation,
\begin{equation}\label{Gauss}
Ric(\nu,\nu) = \frac12 R - K + \frac12 H^2 - \frac12 ||\Pi||^2 
\end{equation}
where $R$ is the scalar curvature of $M^3$ and $K$ is the Gauss curvature of 
$\Sigma(V_0)$.  Substituting we get,
\[A_{V_0}(V_0)^2 A_{V_0}''(V_0) = \int_{\Sigma(V_0)} 
-\frac12 R + K - \frac12 H^2 - \frac12 ||\Pi||^2 \]
Since $\Sigma(V_0)$ has only one component, 
$\int_{\Sigma(V_0)} K = 2\pi \chi (\Sigma(V_0)) \le 4\pi$
by the Gauss-Bonnet theorem.  
Since $R \ge 0$ and $||\Pi||^2 \ge \frac12 H^2$, we have
\begin{eqnarray*}
A_{V_0}(V_0)^2 A_{V_0}''(V_0) &\le& 4\pi - \int_{\Sigma(V_0)} 
\frac34 H^2 \\
     &=& 4\pi - \frac34 H^2 A_{V_0}(V_0)
\end{eqnarray*}
Hence,
\[A_{V_0}''(V_0) \le \frac{4\pi}{A_{V_0}(V_0)^2} - 
                     \frac{3 A_{V_0}'(V_0)^2 }{4 A_{V_0}(V_0)}\]
                   
Finally, since $A(V_0) = A_{V_0}(V_0)$ and 
$A(V) \le A_{V_0}(V)$ for every $V_0 \ge 0$, 

\begin{equation}\label{AD}
A''(V) \le \frac{4\pi}{A(V)^2} - \frac{3 A'(V)^2  }{4 A(V)}             
\end{equation}
in the sense of comparison functions.  \qed

It turns out that $\tilde{F}(V) = A(V)^\frac32$ is more  
convenient to work with
than $A(V)$.  Note that $\tilde{F}$ and $V$
have the same units.  Making this substitution, inequality \ref{AD} becomes 
\begin{equation}\label{R8}
\tilde{F}''(V) \le \frac{36\pi - \tilde{F}'(V)^2}{6\tilde{F}(V)} 
\end{equation}
in the sense of comparison functions.  This last inequality will be the key step
in proving theorem \ref{R4}.  

Now we turn to the other case in which we can prove the Penrose 
conjecture where there may be multiple horizons.  
We assume the hypotheses of theorem \ref{R6} including 
condition \ref{R5} and recall the definition of $F(V)$ 
given in the previous section.  
Under these circumstances, we show in section \ref{s8} that for all $V \ge 0$,
there exists a collection $\{\Sigma_i(V)\}$ 
of smooth surfaces which
minimizes $F$ among collections of surfaces 
which enclose a volume $V$ outside the horizons.  The mean curvature
is constant (but generally different) on each component.
By condition \ref{R5}, we may take $\{\Sigma_i(V)\}$ to be pairwise disjoint, 
and since $\{\Sigma_i(V)\}$
is a minimizer, 
$\sum_{i} \mbox{Area}(\Sigma_i(V))^\frac32 = F(V)$.

$F(V)$ also contains important geometric information,
including the fact that $F(0)$ is the sum of the areas of the horizons
to the three halves power.  The fact 
that the horizons are outermost implies that $F(V)$ is nondecreasing.  

Again, as is proved in section \ref{s6.5}, the total mass $m$ of $(M^3,g)$,
is encoded in the asymptotic behavior of the function $F(V)$ for large $V$,
since we will show
that for sufficiently large $V$, the minimizing 
collection of surfaces is a single large sphere.  
Hence, the key to proving theorem \ref{R6} is understanding how the assumption
of nonnegative scalar curvature on $(M^3,g)$ bounds the behavior of $F(V)$.

\begin{theorem}\label{FPP}
Suppose $(M^3,g)$ is complete, 
has nonnegative scalar curvature, is \\
Schwarzschild
at infinity, and satisfies condition \ref{R5}.  Then 
the function $F(V)$ defined in definition \ref{deffv} satisfies
\[ F''(V) \le \frac{36\pi - F'(V)^2}{6F(V)} \] 
in the sense of comparison functions,
where this means that for all $V_0 \ge 0$ there exists a
smooth function $F_{V_0}(V)\ge F(V)$ with $F_{V_0}(V_0)=F(V_0)$
satisfying
\[ F_{V_0}''(V_0) \le \frac{36\pi - F_{V_0}'(V_0)^2}{6F_{V_0}(V_0)} \] 
\end{theorem}
{\em Sketch of proof.}
The method of the proof here is exactly as in theorem \ref{APP}.
The reason condition \ref{R5} is needed is that if two
components of $\{\Sigma_i(V_0)\}$ push up against each other, then we can not
flow both of the surfaces out at the same time.  
We want to find a
flow on $\{\Sigma_i(V_0)\}$ 
which is constant (but different) on each component.  
Let $F_{V_0}(V)$ be the sum of the areas to the three halves power of these
new surfaces, parameterized as before by the total enclosed volume $V$.
First we consider a flow which is constant on $\Sigma_i(V_0)$ and
zero on all the other components.  As in inequality \ref{R8}, we get that
\[ F_{V_0}''(V_0) \le \frac{36\pi - F_{V_0}'(V_0)^2}{6A_i(V_0)^\frac32} \] 
The next observation to make is that the value we get for $F_{V_0}'(V_0)$ 
is independent of which component we flow out.  
Otherwise, we could find a volume preserving flow which flowed out on
one component and flowed in on the other component 
which decreased the value of $F$.  It follows that  $F_{V_0}'(V_0)$ 
is the same for any flow.
From these observations, it is possible to calculate $F_{V_0}''(V_0)$ for
any flow which is constant on each component.  We then choose the flow
which gives us the best estimate for $F_{V_0}''(V_0)$ which is 
\[ F_{V_0}''(V) \le \frac{36\pi - F_{V_0}'(V_0)^2}{6\sum_i A_i(V_0)^\frac32} 
                 =  \frac{36\pi - F_{V_0}'(V_0)^2}{6F_{V_0}(V_0)} \] 
Finally, since $F(V_0) = F_{V_0}(V_0)$ and 
$F(V) \le F_{V_0}(V)$ for every $V_0 \ge 0$, the theorem follows as before. \qed

\section{The Mass Function}                         
\label{s4}

The function $\tilde{F}(V) = A(V)^\frac32$ will be used to prove 
theorem \ref{R4} and the function $F(V)$ will be used to prove
theorem \ref{R6}.  We choose to abuse
notation slightly from this point on and call both functions $F(V)$
since both functions satisfy
\begin{equation}\label{FPP2}
F''(V) \le \frac{36\pi - F'(V)^2}{6F(V)} 
\end{equation} 
in the sense of comparison functions.    
It always will be clear from the context which function is
intended.
Given an inequality like the one above, it is natural to 
want to integrate it. 

\begin{definition}\label{def1}
For $V \ge 0$, let 
\[m(V) = F(V)^\frac13 \left(36\pi - F'(V)^2 \right) / c \]
be the mass function, where $c = 144\pi^\frac32$. 
\end{definition}

$F(V)$ is continuous, but $F'(V)$ does not necessarily exist
for all $V$, although it does exist almost everywhere
since $F(V)$ is monotone increasing.  
The left and right hand derivatives,
$F'_+(V)$ and $F'_-(V)$, do always exist though.  This follows
from the fact that $F(V)$ has a comparison function $F_{V_0}(V)$ 
(or $A_{V_0}(V)^\frac32$)
which touches $F$ at $V=V_0$ and is greater than $F$ in some
neighborhood of $V_0$ for all $V_0 \ge 0$.  Since the 
second derivatives of the comparison
functions are uniformly bounded from above in a bounded
interval we can add a quadratic to $F(V)$ to get a 
concave function, from which it follows that 
the left and right hand
derivatives always exist and are equal except at a countable 
number of points.

Furthermore, $F'_+(V) \le F'_-(V)$ using the comparison
function argument again since 
$F'_+(V_0) \le F'_{V_0}(V_0) \le F'_-(V_0)$.  If $F'(V)$ 
does not exist, then it is natural to define $F'(V)$ to be a 
multivalued function taking on every value in the interval
$(F'_+(V), F'_-(V))$.  This is consistent, since if $F'(V)$
does exist, then $F'_+(V) = F'_-(V)$.  Hence, $m(V)$ is 
multivalued for some $V$, which can be interpreted as the mass
``jumping up'' at these $V$, and the set of $V$ for which 
$m(V)$ and $F(V)$ are multivalued is a countable set.
Alternatively, one
could replace $F'(V)$ with $F'_+(V)$ (or $F'_-(V)$) in the 
formula for $m(V)$ so that $m(V)$ would always be single 
valued.

\begin{lemma} \label{AJ}
The quantity $m(V)$ is a nondecreasing 
function of $V$. 
\end{lemma}
{\em Proof.}  
The main idea is that if $F(V)$ were smooth,
\[m'(V) = 2 F^\frac13 F'(V) \left(-F''(V) + \frac{36\pi - F'(V)^2}{6F(V)}\right) /c \]
being nonnegative would follow from 
inequality \ref{FPP2} and the fact that $F(V)$ is nondecreasing.

More generally,
it is sufficient to prove that $m'(V)\ge 0$ distributionally.  Hence,
treating $m(V)$ as a distribution we may equivalently define
\[m(V)=F(V)^{1/3}\left(36\pi -F_+'(V)^2\right)\]
since $F_+'(V)=F'(V)$ except at a countable number of points.  It is
convenient to extend $F(V)$ and $m(V)$ to be defined for all real $V$,
so define $F(V)=F(0)$ for $V<0$.  Then since $F'(0)=0$ we still have
\[F''(V)\le\frac{36\pi-F'(V)^2}{6F(V)}\]
in the sense of comparison functions for all $V\in(-\infty,\infty)$,
where we recall that this means that for all $V_0$ there exists a
smooth function $F_{V_0}(V)\ge F(V)$ with $F_{V_0}(V_0)=F(V_0)$
satisfying
\begin{equation}
F_{V_0}''(V)\le\frac{36\pi-F_{V_0}'(V_0)^2}{6F_{V_0}(V_0)}.
\label{eqn:ineqcomp}
\end{equation}

To prove that $m'(V)\ge 0$ distributionally, we will prove that
\[-\int_{-\infty}^\infty m(V)\phi'(V)\,dV\ge 0\]
for all smooth positive test functions $\phi$ with compact support.
We will need the finite difference operator $\Delta_\delta$ which we
will define as
\[\Delta_\delta(g(V))=\frac{1}{\delta}(g(V+\delta)-g(V)).\]

Then
\begin{eqnarray*}
-\int_{-\infty}^\infty m(V)\phi'(V)\,dV &=&
-\int_{-\infty}^\infty F(V)^{1/3}\left(36\pi - F_+'(V)^2 \right)
\phi'(V)\, dV\\
&=&-\lim_{\delta\to 0} \int_{-\infty}^\infty F(V)^{1/3}\left(36\pi -
(\Delta_\delta F(V))^2 \right) (\Delta_\delta \phi(V))\, dV\\
&=&-\lim_{\delta\to 0} \int_{-\infty}^\infty \Delta_{-\delta}
\left\{F(V)^{1/3} \left(36\pi - (\Delta_\delta F(V))^2 \right)\right\}
\phi(V)\, dV\\
\end{eqnarray*}
where we have used the integration by parts formula for the finite
difference operator, $\int  f(x)(\Delta_\delta g(x))\,dx=\int
g(x)\Delta_{-\delta} f(x)\,dx$, which follows from a change of
variables.  Then since $F(V)$ has left-hand derivatives everywhere, in
the limit, we have
\[=\lim_{\delta\to 0} \int_{-\infty}^\infty F(V)^{1/3}\left\{\Delta_{-\delta}
[(\Delta_\delta F(V))^2] +  F_-'(V)
\frac{36\pi-F_+'(V)^2}{3F(V)}\right\} \phi(V)\, dV.\]

Using the comparison functions at each point, since
$F_{V_0}(V_0+\delta) \ge F(V_0 + \delta)$, 
$F_{V_0}(V_0-\delta)\ge F(V_0-\delta)$,
$F_{V_0}(V_0)=F(V_0)$, and $F(V)$ and $F_{V_0}(V)$ are increasing,
it follows that 
\[\Delta_{-\delta}\left[(\Delta_\delta F(V_0))^2\right] \ge \left.
\Delta_{-\delta} \left[  \left(\Delta_\delta F_{V_0}(V) \right)^2
\right] \right|_{V=V_0}.\]
Changing the integration variable to $V_0$, then, we have that
\begin{eqnarray*}
&&-\int_{-\infty}^\infty m(V)\phi'(V)\,dV\ge\\
&&\lim_{\delta\to 0}\int_{-\infty}^\infty F^{1/3}(V_0) \left\{
\left. \Delta_{-\delta} \left[  \left( \Delta_\delta
F_{V_0}(V)\right)^2\right]\right|_{V=V_0} + F_-'(V_0)\frac{36 -
F_+'(V_0)^2}{3F(V_0)} \right\} \phi(V_0)\,dV_0\\
\end{eqnarray*}
and since $F_+'(V_0)=F_-'(V_0)=F_{V_0}'(V_0)$ except at a countable
number of points,
\begin{eqnarray*}
&=&\lim_{\delta\to 0}\int_{-\infty}^\infty F^{1/3}(V_0) \left\{
\left. \Delta_{-\delta} \left[  \left( \Delta_\delta
F_{V_0}(V)\right)^2\right]\right|_{V=V_0} + F_{V_0}'(V_0)\frac{36 -
F_{V_0}'(V_0)^2}{3F_{V_0}(V_0)} \right\} \phi(V_0)\,dV_0\\
&=&\int_{-\infty}^\infty F^{1/3}(V_0) \left\{
-2F_{V_0}'(V_0)F_{V_0}''(V_0) + F_{V_0}'(V_0)\frac{36 -
F_{V_0}'(V_0)^2}{3F_{V_0}(V_0)} \right\} \phi(V_0)\,dV_0\\
&=&\int_{-\infty}^\infty 2F(V_0)^{1/3}F_{V_0}'(V_0) \left\{
-F_{V_0}''(V_0) + \frac{36 - F_{V_0}'(V_0)^2}{6F_{V_0}(V_0)} \right\}
\phi(V_0)\,dV_0\\
&\ge& 0
\end{eqnarray*}
since $F_{V_0}'(V_0)\ge 0$ and the comparison functions satisfy
inequality \ref{eqn:ineqcomp}.  Hence, \\
$m'(V)\ge 0$ distributionally,
so $m(V)$ is a nondecreasing function of $V$.\qed

\section{Proof of the Penrose Inequality 
         Assuming Condition \ref{R3} or \ref{R5}}  \label{s5}

The key to the proofs of theorems \ref{R4} and \ref{R6} is the mass
function $m(V)$.  First, we consider the context of theorem \ref{R4}
so that we have only one horizon and we are minimizing area with a 
volume constraint.  

In section \ref{s6.5}
we show that if $M^3$ is Schwarzschild with mass $m$ at infinity, 
then for large $V$ there is a unique area minimizer,
and that this minimizer is one of the spherically symmetric spheres of the 
Schwarzschild metric.  Hence, for large $V$, the functions 
$F(V) = A(V)^\frac32$ and hence
$m(V)$ are computable in terms of the parameter $m$, and in fact $m(V) = m$.
Also, since the horizon has zero mean curvature, $F'(0) = 0$, so 
$m(0) = \sqrt{\frac{A}{16\pi}}$, the mass of the black hole.
  
Since $m(0)$ equals the mass of the black hole and $m(\infty)$ equals
the total mass of the system, we now see why it is reasonable 
to call the function $m(V)$ mass.  (In fact, $m(V)$ is equal to the 
Hawking mass of $\Sigma(V)$, studied by Christodoulou and Yau in 
\cite{CY} and by
Huisken and Yau in \cite{HY}.)  
Since $m(V)$ is increasing, $m(\infty) \ge m(0)$, so
\[m \ge \sqrt{\frac{A}{16\pi}}\]
and we see that
theorem \ref{R4}, the Penrose inequality for manifolds 
which satisfy condition \ref{R3}, is true.

The proof of theorem \ref{R6} is exactly the same, but we get a stronger
result.  We are back in the context of multiple horizons, and we are 
minimizing the quantity $F$ (from definition \ref{deffv}) given a volume
constraint and assuming condition \ref{R5}.
Again, in section \ref{s6.5} we show that if $M^3$ is Schwarzschild
with mass $m$ at infinity, then for large $V$ there is a unique $F$ 
minimizer, and that this minimizer is a single spherically symmetric
sphere of the Schwarzschild metric.  Thus, once again, $m(V) = m$ for
sufficiently large $V$.  However, while $F'(0) = 0$ again since the horizons
still have zero mean curvature, 
$F(0) = \sum_i A_i^\frac32$, where the $\{A_i\}$
are the areas of the horizons.  Hence, 
\[m(0) = \left(\sum_{i=1}^n 
\left(\frac{A_i}{16\pi}\right)^\frac32 \right)^\frac13 ,\] 
so this time we get
\[ m \ge \left(\sum_{i=1}^n \left(\frac{A_i}{16\pi}\right)^\frac32 \right)
                                                             ^\frac13 \]
which proves that theorem \ref{R6}, the Penrose inequality for manifolds
which satisfy condition \ref{R5}, is true.

\section{Spherical Symmetry at Infinity}
\label{s6}

\begin{definition}
$(M^n,g)$ is said to be asymptotically flat if there is a compact set $K\subset
M$ and a diffeomorphism $\Phi:M-K\to \real^n-\{|x|<1\}$ such that, in
the coordinate chart defined by $\Phi$,
\[g=\sum_{i,j}g_{ij}(x)dx^i dx^j\]
where
\[g_{ij}(x)=\delta_{ij}+O(|x|^{-p})\]
\[|x||g_{ij,k}(x)| + |x|^2|g_{ij,kl(x)}| = O(|x|^{-p})\]
\[|R(g)|=O(|x|^{-q})\]
for some $p>\frac{n-2}{2}$ and some $q>n$, where we have used commas
to denote partial derivatives in the coordinate chart, and $R(g)$ is
the scalar curvature of $(M^n,g)$.
\end{definition}

These assumptions on the asymptotic behavior of $(M^n,g)$ at infinity
imply the existence of the limit
\[M_{ADM}(g)=(4\omega_{n-1})^{-1}\lim_{\sigma\to\infty}
\int_{S_\sigma}\sum_{i,j}(g_{ij,i}\nu_j-g_{ii,j}\nu_j)\,d\mu\]
where $\omega_{n-1}=Vol(S^{n-1}(1))$, $S_\sigma$ is the sphere
$\{|x|=\sigma\}$, $\nu$ is the unit normal to $S_\sigma$ in Euclidean
space, and $d\mu$ is the Euclidean area element of $S_\sigma$.  
The quantity $M_{ADM}$ is called the {\em total mass} of
$(M^n,g)$ (see \cite{ADM}, \cite{Ba2}, \cite{S}, and \cite{SY5}).

\begin{theorem} (Schoen-Yau \cite{SY5})
Let $(M^n,g)$, $n \ge 3$, 
be a complete asymptotically flat $n$-manifold with
$R(g)\ge 0$.  For any $\epsilon>0$, there is a metric $\bar{g}$ such
that $(M^n,\bar{g})$ is asymptotically flat, and outside a compact set
$(M,\bar{g})$ is conformally flat and has $R(\bar{g})=0$, and
$M_{ADM}(\bar{g}) < M_{ADM}(g) + \epsilon$.
\end{theorem}

Furthermore, although Schoen and Yau did not originally 
mention it, 
their proof of the above theorem also proves 
a stronger version of the theorem 
which we will use, namely that the theorem is still true
if we require $|M_{ADM}(\bar{g}) - M_{ADM}(g)| < \epsilon$ and
$\bar{g}$ and $g$ to be $\epsilon$-quasi isometric.
We say that two metrics on $M^n$ are $\epsilon$-quasi isometric if
for all $x\in M^n$
\[e^{-\epsilon} < \frac{\bar{g}(v,v)}{g(v,v)} < 
e^{\epsilon}\]
for all tangent vectors $v\in T_x(M^n)$.

Since $(M,g)$ is conformally flat and scalar flat outside a compact
set, we may choose $\real^n-B_{r_0}(0)$ as a coordinate chart for
$\bar{g}$ for some $r_0>0$, so that
\[\bar{g}_{ij}(x)=u(x)^{\frac{4}{n-2}}\delta_{ij},\quad \mbox{for }|x|>r_0.\]
Since 
\[R(\bar{g})=-\frac{4(n-1)}{n-2}u^{-\left(\frac{n+2}{n-2}\right)}\Delta
u\]
where $\Delta$ is the Euclidean Laplacian and $R(\bar{g})=0$, we see
that $\Delta u=0$ for $|x|>r_0$.  Since $(M^n,g)$ is asymptotically
flat, $u(x)$ tends to 1 as $x$ goes to infinity.

Thus, expanding $u(x)$ in terms of spherical harmonics of $\real^n$,
we find that
\[u(x)=1+\frac{M_{ADM}}{(n-1)|x|^{n-2}}+O(\frac{1}{|x|^{n-1}}).\]

Now we define a new metric $\tilde{g}$, $\epsilon$-quasi isometric to
$\bar{g}$
which will be spherically symmetric with zero scalar curvature outside
a compact set.  To do this, choose any $R>r_0$ and $\delta>0$ and let
\[v(x)=A+\frac{B}{|x|^{n-2}}\]
where $A$ and $B$ are chosen so that
\[A+\frac{B}{R^{n-2}} = \sup_{|x|=R}u(x)+\delta\]
\[A+\frac{B}{(2R)^{n-2}} = \inf_{|x|=2R}u(x)-\delta.\]
Define
\[ w (x)=\left\{\begin{array}{lll}
u(x)&,&|x|<R\\
\min(u(x),v(x))&,&R\le |x|\le 2R\\
v(x)&,&|x|>2R
\end{array}\right.\]
This function is continuous since $ w (x)=u(x)$ for $|x|=R$ and
$ w (x)=v(x)$ for $|x|=2R$.  Furthermore since $u$ and $v$ are harmonic
and the minimum value of two harmonic functions is weakly
superharmonic, $ w $ is weakly superharmonic.

Now define $\tilde{ w }(x)= w *b$ where $b$ is some smooth,
spherically symmetric, positive bump function of total integral 1 and
compact support in $B_\delta(0)$.  Then $\tilde{ w }$ is smooth
and superharmonic and $\tilde{ w }(x)=u(x)$ for $|x|<R-\delta$ and
$\tilde{ w }(x)=v(x)$ for $|x|>R+\delta$, since $u$ and $v$ are
harmonic and hence have the mean value property.

Let $\tilde{g}=\bar{g}$ everywhere except in the region that $\bar{g}$
is conformally flat and scalar flat.  In this region parameterized by
$\real^n-B_{r_0}(0)$, let
\[\tilde{g}_{ij}(x)=\tilde{ w }^{\frac{4}{n-2}}\delta_{ij},\quad\mbox{for
}|x|>r_0.\]
Since $\tilde{ w }$ is superharmonic,
\[R(\tilde{g}) = -\frac{4(n-1)}{n-2}
\tilde{ w }^{-\left(\frac{n+2}{n-2}\right)} \Delta \tilde{w} \ge 0.\] 
Furthermore, if we choose $R$ big enough and $\delta$ small enough we
can make $\tilde{g}$ $\epsilon$-quasi isometric to $g$ and 
\[|M_{ADM}(\tilde{g}) - M_{ADM}(\bar{g})| < \epsilon\]
for any $\epsilon > 0$.
The bound on the mass comes from the fact that
\[M_{ADM}(\tilde{g})=(n-1)AB\]
and choosing $R$ large and $\delta$ small gives $A$ close to 1 and $B$
close to $\frac{M_{ADM}(\tilde{g})}{n-1}$.  Hence we have the
following theorem:

\begin{theorem}
Let $(M^n,g)$, $n \ge 3$, be a complete asymptotically flat $n$-manifold with \\
$R(g)\ge 0$.  For any $\epsilon>0$ there is a metric $\tilde{g}$
$\epsilon$-quasi isometric to $g$ such that $(M^n,g)$ has
$R(\tilde{g})\ge 0$, is asymptotically flat, is spherically symmetric
with $R(\tilde{g})=0$ outside a compact set, and has
$|M_{ADM}(\tilde{g}) - M_{ADM}(g)| < \epsilon$.
\end{theorem}

The statement of this theorem can be simplified by introducing the
following terminology:  We define $(\real^n-\{0\},h)$ to be the
Schwarzschild metric of mass $m$ where
\[h_{ij}(x)=\left(1+\frac{m}{(n-1)|x|^{n-2}}\right)
^{\frac{4}{n-2}}\delta_{ij}.\]
This metric is spherically symmetric, asymptotically flat, has zero
scalar curvature, and has
total mass $m$.

\begin{definition}\label{definfinity}
We say that $(M^n,g)$ is Schwarzschild with mass $m$ at infinity if
$(M^n-K,g)$ is isometric to $(\real^n-B,h)$ for some compact set $K$
in $M^n$ and some ball $B$ in $\real^n$ centered 
around the origin.
\end{definition}
With this definition, the statement of the previous theorem can
be stated as follows.
\begin{theorem}\label{TH3}
Let $(M^n,g)$ be a complete asymptotically flat $n$-manifold with \\
$R(g)\ge 0$ and with total mass $M$.  
For any $\epsilon>0$, there exists a metric $\tilde{g}$
$\epsilon$-quasi isometric to $g$ with $R(\tilde{g}) \ge 0$,
$(M^n,\tilde{g})$ Schwarzschild with mass $m$ at infinity, and
$|m - M| < \epsilon$.
\end{theorem}

Hence, 
since the positive mass and Penrose inequalities are closed 
conditions, we see by theorem \ref{TH3} that without loss of generality, we may 
assume in the statements of the positive mass theorem and the Penrose
inequality that the manifolds are Schwarzschild at infinity.

\section{The Isoperimetric Surfaces of the \\
Schwarzschild Manifold} \label{isoperimetric}

In the previous section we justified the claim that without loss of
generality for proving the Penrose inequality, we could assume that
$(M^3,g)$ is isometric to the Schwarzschild metric of mass $m$ outside
a compact set.  By the positive mass theorem, $m\ge 0$, and since the
standard, flat $\real^3$ metric is the only metric with $m=0$, we
generally have $m>0$.  In the next two sections, we prove that the
isoperimetric surfaces of $M^3$, the surfaces $\Sigma(V)$ which
minimize area given a volume constraint $V$, are the
spherically-symmetric spheres of the Schwarzschild metric when $V$ is
large enough.

We recall that the Schwarzschild metric of mass $m$ is given by
$(\real^3-\{0\},h)$ where $h_{ij}=(1+\frac{m}{2r})^4\delta_{ij}$.  The
metric is spherically symmetric, asymptotically flat, has zero scalar
curvature, and has an outermost minimal sphere at $r=m/2$.  In fact,
the Schwarzschild metric is symmetric under the mapping $r\to
\frac{m^2}{4r}$ and so has two asymptotically flat ends.

\begin{theorem}
In the Schwarzschild metric of mass $m \ge 0$, 
$(\real^3-\{0\},h)$, described above,
the spherically symmetric spheres given by $r=\mbox{constant}$
minimize area among all other surfaces in their homology class
containing the same volume. 

\label{Aschwarz}
\end{theorem}
{\em Proof.}
Since there is an infinite amount of volume
inside the horizon of the \\
Schwarzschild metric, 
we first comment that ``containing the same volume''
is a well defined notion among surfaces in the same homology
class.  Equivalently, one could define the volume contained
by a surface in the horizon's homology class to be the volume
contained by the region outside the horizon keeping track of
signs if the region is not entirely outside the horizon.

Let $\Sigma$ be a spherically symmetric sphere $r=c>m/2$ of the
Schwarzschild metric $(\real^3-\{0\},h)$.  We will prove that $\Sigma$ is an
isoperimetric surface of  
$(\real^3-\{0\},h)$, and the case when $r<m/2$ will
follow from the symmetry of the Schwarzschild metric under
$r\to\frac{m^2}{4r}$.  We omit the case $r=m/2$, but in this case it
is easy to show that $\Sigma$ minimizes area among all surfaces even
without the volume constraint (see figure \ref{Schwarz}).

By direct calculation, it is easy to compute that the Hawking mass of
$\Sigma$ is always $m$, that is
\begin{eqnarray*}
m&=&\left(\frac{A}{16\pi}\right)^{1/2}\left(1-\frac{1}{16\pi}\int_\Sigma
H^2\right)\\
&=&\left(\frac{A}{16\pi}\right)^{1/2}\left(1-\frac{1}{16\pi}H^2 A\right)
\end{eqnarray*}
where $A$ is the area of $\Sigma$ and $H$ is the mean curvature of
$\Sigma$ which is constant on $\Sigma$ by symmetry.  Since $m>0$,
$H^2A<16\pi$ for $\Sigma$, and since $c>m/2$, it is easy to check that
$H>0$.  Notice that we have already used the positivity of the 
mass $m$.

Now we construct a new metric $(\real^3,k)$ which is isometric to
$(\real^3-\{0\}, h)$ outside $\Sigma$ but is isometric to a
spherically symmetric connected neighborhood of the tip of a 
spherically symmetric cone inside $\Sigma$, where
the proportions of the cone are chosen to give $\Sigma$ the same area
$A$ and mean curvature $H$ from the inside as the outside.  The form
of the metric $(\real^3,k)$ can then be written
\[ds^2_k=u(r)^{-2}dr^2+u(r)r^2d\sigma^2\]
in spherical coordinates $(r,\vec{\sigma})$ in $\real^3$, where
$d\sigma^2$ is the standard metric on the sphere of radius 1 in
$\real^3$.  Notice that $u(r)\equiv 1$ would represent the standard
flat metric of $\real^3$, and that $u(r)\equiv \mbox{constant}$ gives
a cone.  But the main point of this form 
for the metric is
that the volume element of $ds^2$ is the standard volume element in
$\real^3$ no matter what $u(r)$ is.

\begin{figure}
\vspace{3.5in}
\includegraphics{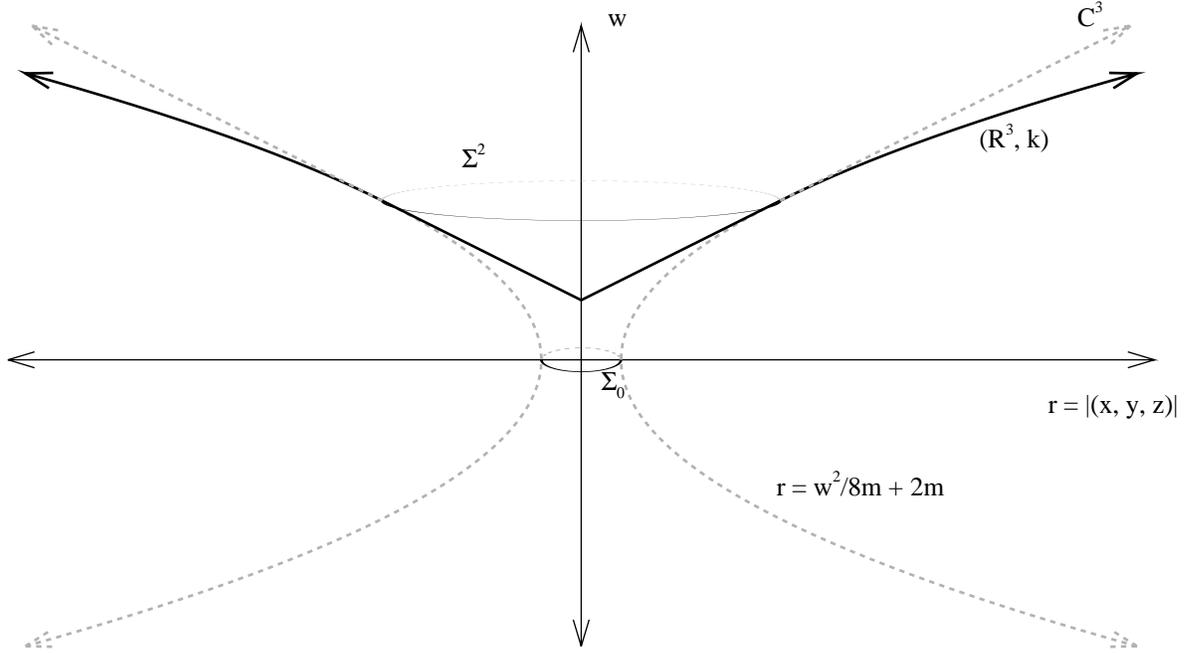}
\caption{Picture of $(\real^3,k)$ isometrically embedded in four-dimensional 
Euclidean space.\label{Schwarzcone}}
\end{figure}

It is helpful to view $(\real^3,k)$ as a submanifold of $\real^4$
(see figure \ref{Schwarzcone}).  
We recall from figure \ref{Schwarz} 
that the submanifold $|(x,y,z)|=\frac{w^2}{8m}+2m$
is the Schwarzschild metric of mass $m$.  Let $C^3$ be the
cone in $\real^4$ which is tangent to the Schwarzschild metric on $\Sigma^2$.
Then $(\real^3,k)$ is the spherically symmetric solid portions 
shown in figure \ref{Schwarzcone},
equal to the union of the Schwarzschild 
metric outside $\Sigma^2$ and 
the cone $C^3$ inside $\Sigma^2$.

Suppose $\Sigma$ is at $r=\bar{c}$ in this new metric $(\real^3,k)$.
Since $(\real^3,k)$ is a cone inside $\Sigma$, $u(r)=a$ for
$r<\bar{c}$ for some constant $a$.  Since $H^2A<16\pi$ for $\Sigma$,
it turns out that $0<a<1$.

In the Schwarzschild metric, if we parameterized the areas of the
spherically symmetric spheres by the enclosed volume (outside
the horizon) we get from the vanishing of the scalar curvature 
that $A(V)$ satisfies
\begin{equation}
A^2A''(V)=4\pi-\frac{3}{4}A'(V)^2A
\label{eqn:s1}
\end{equation}

We will now reparameterize $u(r)=\bar{u}(V)$ where $V=\frac{4}{3}\pi
r^3$ is the enclosed volume.  This is convenient since the coordinate
chart is volume preserving.
Hence, since $(\real^3,k)$ is Schwarzschild outside $\Sigma$, we have
\begin{eqnarray*}
\bar{u}(V)&=&\left\{\begin{array}{lll}
a&,&V<V_0\\
\frac{A(V)}{(36\pi)^{1/3}V^{2/3}}&,&V\ge V_0
\end{array}\right.
\end{eqnarray*}
where $V_0=\frac{4}{3}\pi\bar{c}^3$ and $A(V)$ satisfies the second
order differential equation (\ref{eqn:s1}) with initial conditions
$A(V_0)$ and $A'(V_0)$ such that
\[\bar{u}(V_0)=a\]
and
\[\bar{u}'(V_0)=0.\]

These last two initial conditions guarantee that $\Sigma$ (which is
the sphere at the boundary between the cone and the Schwarzschild
metric) has the same area and mean curvatures on the inside and the
outside, and that consequently the curvature of $(\real^3,k)$ is
bounded on $\Sigma$.
Finally,
\[\bar{u}(V)=\frac{A(V)}{(36\pi)^{1/3}V^{2/3}}\]
for $V\ge V_0$ since this factor guarantees that the sphere containing
a volume $V$ in $(\real^3,k)$ will have area $A(V)$, and hence
$(\real^3,k)$ will be Schwarzschild outside $\Sigma$.

\begin{lemma}  The following inequality holds.
\[a\le \bar{u}(V) \le 1.\]
\end{lemma}
{\em Proof.}
First we show that $\bar{u}(V)<1$.  To do this we note that
$(\real^3,k)$ has scalar curvature $R(k)\ge 0$ everywhere since $a<1$.
If we redefine $A(V)$ to be the area of the spherically symmetric sphere in
$(\real^3, k)$ containing a volume $V$, then by direct calculation we
have
\[A^2A''(V)\le 4\pi-\frac{3}{4}A'(V)^2A.\]
for all $V \ge 0$, with equality outside $\Sigma$.
Integrating this inequality implies that $m'(V)\ge 0$ where
\[m(V)=\left(\frac{A}{16\pi}\right)^{1/2}\left(1-\frac{1}{16\pi} 
A(V)A'(V)^2\right)\]
and hence that $m(V)\ge 0$ since $m(0)=0$.
Thus, $A'(V)^2 < \frac{16\pi}{A}$, from which it follows that
\[A(V)^3\le 36\pi V^2\]
which implies $\bar{u}(V)\le 1$.

The fact that $\bar{u}(V)\ge a$ follows from the fact that
$\bar{u}'(V)\ge 0$.  We show that $\bar{u}(V)$ is increasing for $V\ge
V_0$ by proving that $\bar{u}'(V)\le 0$ would imply $\bar{u}''(V)\ge 0$.
Hence, since $\bar{u}'(V_0)=0$, it follows that the minimum value of
$\bar{u}'(V)$ for $V\ge V_0$ is zero, so $\bar{u}'(V)\ge 0$.
We compute for $V\ge V_0$.

\[\bar{u}(V)=\frac{A(V)}{(36\pi)^{1/3}V^{2/3}}\]
\[(36\pi)^{1/3}\bar{u}'(V) = (A'(V)-\frac{2}{3}AV^{-1}) V^{-2/3}\]
so $\bar{u}'(V)\le 0$ implies
\[A'(V)\le \frac{2}{3}AV^{-1}\]
But
\[(V^{8/3}A^2(36\pi)^{1/3})\bar{u}''(V)=A^2 A''(V)V^2-\frac{4}{3}A'(V)
A^2 V + \frac{10}{9}A^3\]
Since $A^2A''(V)=4\pi-\frac{3}{4}A'(V)^2A$ in the Schwarzschild
metrics, it follows that for $V \ge V_0$,
\[(V^{8/3}A^2(36\pi)^{1/3})\bar{u}''(V)=(4\pi-\frac{3}{4}A'(V)^2A)V^2
- \frac{4}{3}A'(V)A^2 V + \frac{10}{9}A^3.\]
If $\bar{u}'(V)\le 0$, then $A'(V)\le \frac{2}{3}AV^{-1}$, so
\[(V^{8/3}A^2(36\pi)^{1/3})\bar{u}''(V)\ge \frac{1}{9}(36\pi
V^2-A^3)\ge 0\]
from before.  Thus, for $V\ge V_0$, $\bar{u}'(V)\ge 0$, so $\bar{u}\ge
a$.
This completes the proof that $a\le \bar{u}(V)<1$, so it follows that
$a\le u(r)\le 1$ for all $r$, too.  \qed

Now we prove that $\Sigma$ at $r=\bar{c}$ is an isoperimetric sphere
of $(\real^3,k)$.  Let $\tilde{\Sigma}$ be any other surface in
$(\real^3, k)$ containing the same volume $V_0$ as $\Sigma$ 
(or greater volume).
Let $A$ and $\tilde{A}$ be the $(\real^3,k)$ areas of $\Sigma$ and
$\tilde{\Sigma}$ respectively and let $A_0$ and $\tilde{A_0}$ be the
areas of $\Sigma$ and $\tilde{\Sigma}$ in the $\real^3$ coordinate
chart where $(\real^3,k)$ is represented by
\[ds^2_k=u(r)^{-2}dr^2+u(r)r^2d\sigma^2,\,\,\,\,\,\,
a\le u(r)\le 1.\]
Since this coordinate chart is volume preserving, $\Sigma$ and
$\tilde{\Sigma}$ both contain the same volume $V_0$ in the $\real^3$
coordinate chart.  Hence, by the isoperimetric inequality,
$\tilde{A}_0\ge A_0$.  Thus, since $u(r)\ge a$ and $u(r)^{-2}\ge a$,
and $u(\bar{c})=a$, we have
\[\tilde{A}\ge a\tilde{A}_0 \ge a A_0 = A.\]
Hence $\Sigma$ minimizes area among all surfaces which contain a
volume $V_0$ in $(\real^3,k)$.

Since $(\real^3-\{0\},k)$ and $(\real^3-\{0\},h)$, the Schwarzschild
metric, are both spherically symmetric, they are conformally
equivalent.  In fact
we can represent $(\real^3-\{0\},h)$ as
\[ds^2_h=w(r)^4\left(u(r)^{-2}dr^2+u(r)r^2 d\sigma^2\right)\]
where $w(r)\equiv 1$ for $r\ge \bar{c}$.  Furthermore, by the scalar
curvature formula for conformal metrics \cite{S}, 
$\Delta w \ge 0$, so 
$w(r)>1$ for $r<\bar{c}$
since Schwarzschild has zero scalar curvature and the cones with
$u(r)=a<1$ have positive scalar curvature.

Now we prove $\Sigma$ at $r=\bar{c}$ minimizes area among all surfaces
in the Schwarzschild metric in its homology class containing the same
(relative) volume.  Let $\tilde{\Sigma}$ be any other such sphere.
Since $\tilde{\Sigma}$ contains the same volume in $(\real^3-\{0\},h)$
as $\Sigma$, it must contain more volume than $\Sigma$ in the
$(\real^3,k)$ metric under the conformal identification since $w\ge
1$.  Hence it has more area in the $(\real^3,k)$ metric since $\Sigma$
is an isoperimetric sphere of $(\real^3,k)$.  Then again, since $w\ge
1$ but $w(\bar{c})=1$, $\tilde{\Sigma}$ must have more area than
$\Sigma$ in the Schwarzschild metric $(\real^3-\{0\},h)$.  Thus, we
have proved theorem \ref{Aschwarz}.  \qed

Now we consider minimizing $F$, ``the sum of the areas to the
three halves power'' from definition \ref{deffv}, with a volume
constraint in the Schwarzschild metric.  Using the same argument
as in the proof of theorem \ref{Aschwarz}, we find that the 
collection of surfaces $\{\Sigma_i(V)\}$ which minimizes $F$
among collections of surfaces containing the horizon and 
a volume $V$
outside the horizon is always a single spherically symmetric 
sphere of the Schwarzschild metric.

The only real difference in the proof is understanding 
minimizers of $F$ in $\real^3$.  Whereas a sphere minimizes 
area given a volume constraint in $\real^3$, any collection of
spheres minimizes $F$ in $\real^3$ given a volume constraint.
This follows from the isoperimetric inequality,
$A_i^\frac32 \ge \sqrt{36\pi} V_i$ with equality for spheres.
Hence, $\sum_i A_i^\frac32 \ge \sqrt{36\pi} V$ with equality
for collections of spheres.  However, a single sphere containing
a volume $V$ still minimizes $F$, and so the proof from theorem
\ref{Aschwarz} still applies.

\begin{theorem}
In the Schwarzschild metric of mass $m \ge 0$, 
$(\real^3-\{0\},h)$, 
the spherically symmetric spheres
minimize F among all other surfaces in their homology class 
containing the same volume.
\label{Fschwarz}
\end{theorem}

\section{Mass and Isoperimetric Spheres at Infinity}
\label{s6.5}

Now we consider manifolds $(M^3,g)$ which are Schwarzschild of
mass $m$ at 
infinity (see definition \ref{definfinity}), are complete, and 
have nonnegative scalar curvature.  By the positive mass theorem,
$m \ge 0$.  Since $(M^3,g)$ is isometric to the Schwarzschild
metric outside a compact set, we expect that when we minimize
area with a volume constraint $V$, the minimizers are still the
spherically symmetric spheres of the Schwarzschild metric when
$V$ is large enough.  In fact this is the case not only for
area minimization, but is also true when we minimize $F$ with
a volume constraint.  To prove this, we begin with three
definitions and a lemma.

\begin{definition}
Suppose $\Sigma^2 = \partial U^3 \subset M^3$ minimizes area among all
surfaces bounding a compact region of the same volume, $|U^3|$.  
Then we call $\Sigma^2$ an
isoperimetric surface of $M^3$.
\label{defn:isoperi}
\end{definition}

\begin{definition}
Suppose $\Sigma^2 = \partial U^3 \subset M^3$ minimizes area among all
surfaces bounding a compact region of volume greater than or equal to $|U^3|$.
Then we call $\Sigma^2$ an
outer isoperimetric surface of $M^3$.
\label{defn:outeriso}
\end{definition}

\begin{definition}\label{area_nonincreasing}
A mapping $\phi:A^3 \rightarrow B^3$ is area nonincreasing 
if and only if for all surfaces with boundary $\Sigma^2 \subset A^3$, the area
of $\phi(\Sigma^2)$ is less than or equal to the area of $\Sigma^2$.
\end{definition}

\begin{lemma}
Suppose $\Sigma^2 = \partial U^3$ is a smooth surface in $M^3$, 
$U^3$ is compact, and 
there exists a $C^1$, onto, area nonincreasing mapping 
$\phi:M^3\to N^3$, which is an isometry outside of the interior of $U^3$, 
such that
$\phi(\Sigma^2) = \partial(\phi(U^3))$ is an outer isoperimetric surface 
of $N^3$.  Then $\Sigma^2$ is an outer isoperimetric surface of $M^3$.  
\label{lem:outer}
\end{lemma}

{\em Proof.}  First we claim that $\phi$ is volume nonincreasing inside
$\Sigma^2$.  Let $\{e_i\}$ be an orthonormal basis at some point $p\in
U^3$, the region contained by $\Sigma^2$.  Let $G_{ij}=\langle
D\phi(e_i), D\phi(e_j) \rangle_{N^3}$, and
$\bar{G}_{ij}=G_{ij}^{-1}\cdot \det(G)$.  Then for $\vec{v}\in
T_p(M^3)$, $v_i^t G_{ij}v_j$ is the square of the length of
$D\phi(\vec{v})$ and $v_i^t\bar{G}_{ij}v_j$ is the factor by which
areas (more generally, $(n-1)$-volumes) get increased in the direction
orthogonal to the unit vector $\vec{v}$.  Since we are given that areas are not
increased, all the eigenvalues of $\bar{G}$ (which is symmetric and so
has all real eigenvalues) are less than or equal to one.  Thus
$\det\bar{G}\le 1$, which implies $\det(G)\le 1$ from the formula for
$\bar{G}$.  But $\det(G)$ is the factor by which volumes are changed 
at $p$, so $\phi$ is volume nonincreasing inside $\Sigma^2$.

Now we prove $\Sigma^2$ is outer isoperimetric.  Let
$\bar{\Sigma}^2$ be any competitor for $\Sigma^2$, that
is, suppose $\bar{\Sigma}^2$ 
contains at least as much volume in $M^3$ as $\Sigma^2$.  Then
upon reflection we see that $\phi(\bar{\Sigma}^2)$ must contain at least as
much volume as $\phi(\Sigma^2)$ in $M^3$ since $\phi$ is volume nonincreasing.
Hence, since $\phi(\Sigma^2)$ is outer
isoperimetric, $\mbox{Area}(\phi(\bar{\Sigma}^2))\ge \mbox{Area}(\phi(\Sigma^2))$.  
But $\phi$ is an isometry on $\Sigma^2$ and 
area nonincreasing everywhere, so $\mbox{Area}(\bar{\Sigma}^2)\ge
\mbox{Area}(\Sigma^2)$.  \qed

Using this lemma, we can prove theorem \ref{Aschwarz2}

\begin{theorem}\label{Aschwarz2}
Suppose $(M^3,g)$ is complete, has nonnegative scalar curvature,
and is Schwarzschild at infinity.  Then there exists a $V_0$
such that for all $V \ge V_0$, the spherically 
symmetric spheres of the Schwarzschild metric minimize area
among all other surfaces in their homology class 
containing the same volume $V$ (outside
the horizons, if any exist).
\end{theorem}

{\em Proof.}  Since $M^3$ is isometric to the Schwarzschild manifold
of mass $m$ outside a compact set, then for some $A_{min}\ge 16\pi
m^2$, $\Sigma^2(A)\subset M^3$ exists for $A\ge A_{min}$ and is the
spherically symmetric sphere of area $A$ of the Schwarzschild portion
of $M^3$.  We claim that the spheres of area
\begin{equation}
A\ge\frac{1}{\pi}\left(\frac{A_{min}}{m}\right)^2
\label{eqn:aineq}
\end{equation}
must be outer isoperimetric spheres of $M^3$.  Let $\Sigma^2\subset
M^3$ be one of these spheres with area satisfying inequality
\ref{eqn:aineq}.  As before in the proof of theorem 
\ref{Aschwarz}, we
construct a spherically symmetric manifold $(\real^3,k)$ which is
isometric to the Schwarzschild manifold of mass $m$ outside a
spherically symmetric sphere $\bar{\Sigma}^2$ and is isometric to a
cone inside $\bar{\Sigma}^2$, where the proportions of the cone are
chosen so that $\bar{\Sigma}^2$ has the same area and mean curvature
from the inside as the outside.  As we proved in theorem 
\ref{Aschwarz},
$\bar{\Sigma}^2$ is outer isoperimetric in $(\real^3,k)$.  We perform
this construction so that $(\real^3,k)$ has the same mass as $M^3$ and
$\bar{\Sigma}^2$ has the same area as $\Sigma^2$.  Hence $(\real^3,k)$
outside $\bar{\Sigma^2}$ is isometric to $M^3$ outside $\Sigma^2$.

We want to construct a map $\phi:M^3 \to (\real^3,k)$ which satisfies
the conditions of lemma \ref{lem:outer}.  Define $\phi$ to be the
identity isometry map outside $\Sigma^2$ so that
$\bar{\Sigma}^2=\phi(\Sigma^2)$.  Inside $\Sigma^2$ we will make
$\phi$ spherically symmetric where $M^3$ is spherically symmetric, so 
in the region where $\phi$ is injective we can
characterize $\phi$ with the function $A(\bar{A})$,
where $A$ is the area of the spherically symmetric pre-image in 
$M^3$ of
the spherically symmetric sphere in $(\real^3,k)$ of area $\bar{A}$.
Hence, if we let $A_0$ be the area of $\Sigma^2$ and of
$\bar{\Sigma}^2$, then $A(A_0)=A_0$.  We define $A(\bar{A})$ for
$\bar{A}<A_0$ so that $A'(\bar{A})$ is as small as possible for each
$\bar{A}<A_0$ such that $\phi$ is area nonincreasing.

In fact, since $A\ge \bar{A}$, lengths in the spherically symmetric
directions get decreased by a factor of
$\left(\frac{\bar{A}}{A}\right)^{1/2}$ by $\phi$, so that if we define
$\phi$ to increase lengths by a factor of
$\left(\frac{\bar{A}}{A}\right)^{-1/2}$ in the radial direction, $\phi$
will be area nonincreasing.  Hence, volumes will be decreased by
$\phi$ locally by a factor of 
$\left(\frac{\bar{A}}{A}\right)^{1/2}
\left(\frac{\bar{A}}{A}\right)^{1/2}
\left(\frac{\bar{A}}{A}\right)^{-1/2}
=
\left(\frac{\bar{A}}{A}\right)^{1/2}$.

It is convenient to consider the spherically symmetric functions $V$
on $M^3$ and $\bar{V}$ on $(\real^3,k)$, where $V$ and $\bar{V}$ are
volumes enclosed by the corresponding spherically symmetric spheres.
We note that $\bar{V}$ is defined everywhere on $(\real^3,k)$ but $V$
is defined only where $M^3$ is Schwarzschild and hence spherically
symmetric.  Then it is easy to compute that on the cone,
$\bar{A}(\bar{V}) = a(36\pi)^{1/3}\bar{V}^{2/3}$, so that for some
constant $a$ 
\begin{equation}
\bar{A}'(\bar{V})=a^{3/2}\cdot\frac{2}{3}\cdot(36\pi)^{1/2}\bar{A}^{-1/2}.
\label{eqn:aprime}
\end{equation}
Furthermore, 
in the Schwarzschild portion of $M^3$, since from the Hawking mass
we have $m= \left(\frac{A}{16\pi}\right)^{1/2} \left( 1 -
\frac{A}{16\pi} A'(V)^2 \right)$,
\begin{equation}
A'(V)=\sqrt{\frac{16\pi}{A}\left(1-m\left(\frac{16\pi}{A}\right)^{1/2}\right)}
\label{eqn:asqrt}
\end{equation}
Also, since we already noted that $\phi$ decreases volumes locally by
a factor of $\left(\frac{\bar{A}}{A}\right)^{1/2}$,
\begin{equation}
\frac{dV}{d\bar{V}}=\left(\frac{A}{\bar{A}}\right)^{1/2}.
\label{eqn:dvdv}
\end{equation}
Then, since 
\[\frac{dA}{dV}\cdot\frac{dV}{d\bar{V}}=\frac{dA}{d\bar{A}}
\cdot \frac{d\bar{A}}{d\bar{V}},\] 
we find that
\begin{equation}
A'(\bar{A})=a^{-3/2}\sqrt{1-m\left(\frac{16\pi}{A}\right)^{1/2}}
\label{eqn:aprime2}
\end{equation}
for $\bar{A}\le A_0$ with initial condition $A(A_0)=A_0$.  Hence, this
differential equation determines $A(\bar{A})$.

We claim that $A(0)\ge A_{min}$ if $A_0\ge\frac{1}{\pi} \left(
\frac{A_{min}}{m} \right)^2$ as in inequality (\ref{eqn:aineq}).  This
will prove that an area nonincreasing map $\phi$ exists from
$M^3$ to $(\real^3,k)$, where $\phi$ will be defined to map 
everything in $M^3$ inside the spherically symmetric sphere of area $A(0)$ 
to the tip of the cone $(\real^3,k)$.  Actually, this mapping is not $C^1$ 
as required in the lemma, but the mapping can be perturbed to be $C^1$
and still stay area nonincreasing. 

We need the inequality
\[\sqrt{b-x}\le\sqrt{b}(1-\frac{1}{2}x)\]
for $0\le x\le b\le 1$.  Hence,
\begin{eqnarray*}
A'(\bar{A})&=&a^{-3/2}\sqrt{1-m\left(\frac{16\pi}{A_0}\right)^{1/2} -
\left[m(16\pi)^{1/2} \left(A^{-1/2}-A_0^{-1/2}\right)\right]}\\
&\le&a^{-3/2}\sqrt{1-m\left(\frac{16\pi}{A_0}\right)^{1/2}}\left[ 1  -
\frac{m}{2}(16\pi)^{1/2} \left(A^{-1/2}-A_0^{-1/2}\right)\right]
\end{eqnarray*}
where we will verify later that $m\left(\frac{16\pi}{A}\right)^{1/2}\le
1$ for all $\bar{A}\ge 0$.  Since the dimensions of the cone
(including $a$) were chosen so that $\bar{\Sigma^2}$ had the same mean
curvature in $(\real^3,k)$ on the inside as on the outside,
$A'(A_0)=1$ since to first order $A(\bar{A})\cong \bar{A}$ for
$\bar{A}$ near $A_0$.  Thus, $a$ must satisfy
\[1=A'(A_0)=a^{-3/2}\sqrt{1-m\left(\frac{16\pi}{A_0}\right)^{1/2}},\]
so
\[A'(\bar{A})\le 1- \frac{m}{2} (16\pi)^{1/2} \left(A^{-1/2} -
A_0^{-1/2} \right).\]
Hence, if we let $D(\bar{A})=A(\bar{A})-\bar{A}$, then
\[D'(\bar{A})\le -\frac{1}{2}m(16\pi)^{1/2}\left((D(\bar{A})+\bar{A})^{-1/2}
-A_0^{-1/2}\right).\]
Since $D(A_0)=0$, it follows that $D(\bar{A})\ge 0$ for $\bar{A}\le
A_0$ and that $D'(\bar{A})\le 0$.  Thus, $D(\bar{A})$ attains its
maximum value at zero.  Hence
\[D'(\bar{A})\le -\frac{1}{2}m(16\pi)^{1/2}\left((D(0)+\bar{A})^{-1/2}
-A_0^{-1/2}\right)\]
so that if we integrate both sides from $\bar{A}=0$ to $\bar{A}=A_0$
we get 
\begin{eqnarray*}
D(A_0)-D(0)&\le&\int_0^{A_0} -\frac{1}{2}m(16\pi)^{1/2}\left((D(0) +
\bar{A})^{-1/2}-A_0^{-1/2}\right)\,d\bar{A}\\
&=& -\frac{1}{2}m(16\pi)^{1/2}\left[2(D(0) +
\bar{A})^{1/2}-A_0^{-1/2}\bar{A}\right]_{\bar{A}=0}^{\bar{A}=A_0}\\
&=& -\frac{1}{2}m(16\pi)^{1/2}\left[2(D(0) + A_0)^{1/2}-
2D(0)^{1/2}-A_0^{1/2}\right]
\end{eqnarray*}
so that since $D(A_0)=0$ and $(D(0)+A_0)^{1/2}\ge A_0^{1/2}$, 
\begin{equation}
D(0)+m(16\pi)^{1/2}D(0)^{1/2}\ge \frac{1}{2}m(16\pi)^{1/2}A_0^{1/2}.
\label{eqn:d0}
\end{equation}
Since $A_0\ge\frac{1}{\pi}\left(\frac{A_{min}}{m}\right)^2$ from
inequality (\ref{eqn:aineq}) and $A_{min}\ge 16\pi m^2$ since the
minimal sphere in the Schwarzschild manifold has area $16\pi m^2$,
\[D(0)+m(16\pi)^{1/2}D(0)^{1/2}\ge 32\pi m^2.\]
Hence, $D(0)^{1/2} \ge m(16\pi)^{1/2}$, 
since the left side of the above
inequality is an increasing function of $D(0)$.  Thus, plugging this
into inequality (\ref{eqn:d0}) we get
\[2D(0)\ge\frac{1}{2}m(16\pi)^{1/2}A_0^{1/2}\]
so that from inequality (\ref{eqn:aineq}) we have
\[D(0)\ge A_{min}.\]
But $D(0)=A(0)-0=A(0)$, so
\[A(0)\ge A_{min}\]
which means that we have stayed in the spherically symmetric portion
of $M^3$ for $0\le\bar{A}\le A_0$.  We notice that the spherically 
symmetric sphere in $M^3$ of area $A(0)$ gets mapped to the tip of the
cone $(\real^3,k)$, so we might as well define $\phi$ to send everything
inside the sphere of area $A(0)$ in $M^3$ to the tip of the cone.
Certainly this is an area nonincreasing map.
Thus, we have defined a 
mapping $\phi:M^3 \to (\real^3.k)$ which is an isometry outside of
$\Sigma$ in $M^3$ and which is area decreasing inside
of $\Sigma$ in $M^3$.  The mapping is not $C^1$ on the sphere of 
area $A(0)$ in $M^3$, but $\phi$ can be perturbed slightly
around the sphere of area $A(0)$ so that
it is $C^1$ and still area nonincreasing.    
Since $\bar{\Sigma}$ is outer
isoperimetric in $(\real^3,k)$ and $\phi(\Sigma) = \bar{\Sigma}$,
it follows from lemma \ref{lem:outer}
that $\Sigma$ is outer isoperimetric in $M^3$ and
hence minimizes area among surfaces containing the same volume in
$M^3$.  Since $\Sigma$ was any of the spherically symmetric spheres 
of area $A \ge \frac{1}{\pi}
\left( \frac{A_{min}}{m}\right)^2$, this proves theorem 
\ref{Aschwarz2}.
\qed

Since spheres are $F$-minimizers given a volume constraint,
we mentioned in the previous section that $\bar{\Sigma}^2$
minimizes $F$ as well as area given a volume constraint in 
$(\real^3,k)$.  Thus, lemma \ref{lem:outer} implies
theorem \ref{Fschwarz2} as well as theorem \ref{Aschwarz2}.

\begin{theorem}\label{Fschwarz2}
Suppose $(M^3,g)$ is complete, has nonnegative scalar curvature,
and is Schwarzschild at infinity.  Then there exists a $V_0$
such that for all $V \ge V_0$, the spherically 
symmetric spheres of the Schwarzschild metric minimize F
among all other surfaces in their homology class 
containing the same volume $V$ (outside
the horizons, if any exist).
\end{theorem}

Theorems \ref{Aschwarz}, \ref{Fschwarz}, \ref{Aschwarz2}, and
\ref{Fschwarz2} are also true in higher dimensions except the 
exponent in the definition of $F$ is more generally 
$\frac{n}{n-1}$ instead of $\frac32$.  Also, note that theorems
\ref{Aschwarz2} and \ref{Fschwarz2} are not true if we merely
require $(M^3,g)$ to be Schwarzschild at infinity and drop
the conditions that $(M^3,g)$ is complete and has nonnegative
scalar curvature.  We need $R(g) \ge 0$ and completeness to use
the positive mass theorem to conclude that $m \ge 0$, which is 
essential for theorems \ref{Aschwarz} and \ref{Fschwarz}.  In 
fact, for $m < 0$, the spherically symmetric spheres of the 
Schwarzschild metric are unstable and hence do not minimize
area among surfaces enclosing the same volume.

Theorems \ref{Aschwarz2} and \ref{Fschwarz2} are important
because they allow us to evaluate $\lim_{V \rightarrow \infty}
m(V)$.  In fact, since the minimizing surface (when minimizing
area) or collection of surfaces (when minimizing F) enclosing
a volume V outside the horizons is always a spherically
symmetric sphere of the Schwarzschild metric for $V \ge V_0$,
$m(V) = m$, the mass parameter of the Schwarzschild metric,
for $V > V_0$.

\begin{theorem}\label{masslimit}
Suppose $(M^3,g)$ is complete, has nonnegative scalar curvature,
and is Schwarzschild with mass $m$ at infinity.  Then for
both definitions of $m(V)$ in section \ref{def1}
(whether we are minimizing area
or $F$ with a volume constraint), we have
\[ \lim_{V \rightarrow \infty} m(V) = m. \]
\end{theorem}

This theorem is true in higher dimensions as well.
We also conjecture that theorem \ref{masslimit} is true for 
asymptotically flat manifolds with positive total mass, 
where $m$ is replaced by the 
total mass $M_{ADM}$.  With additional decay conditions on the
asymptotic flatness of $(M^3,g)$, Huisken and Yau show that
the region at infinity is foliated by constant mean curvature
spheres which are stable and hence locally minimize area with
a volume constraint \cite{HY}.
We conjecture that these same spheres also
globally minimize area among surfaces in the same homology class
containing the same volume.

\section{Existence of Surfaces which Minimize Area \\
Given a Volume Constraint}
\label{s7}

In our proof of theorem \ref{R4}, we used the fact that the mass function
$m(V)$ is nondecreasing which relies on equation \ref{R8}, which in turn
followed from doing a unit normal flow on $\Sigma(V)$, the surface which
minimizes area among surfaces containing 
a volume $V$ outside the horizon.  Thus, it
is essential that the surface $\Sigma(V)$ actually exists. 

In this section, we assume the hypotheses of theorem \ref{R4} again, 
including condition \ref{R3}, and recall the definition of $A(V)$
given in the introduction to this chapter.  
We will prove that for all $V \ge 0$, 
there exists a surface $\Sigma(V)$ which encloses a volume $V$ outside
the horizon and which minimizes area, so that 
$\mbox{Area}(\Sigma(V)) = A(V)$.

Existence theory for compact manifolds is well understood using geometric
measure theory, since the space of rectifiable currents 
of bounded mass on compact manifolds
is compact.
The main problem with this type of existence question is that $(M^3,g)$
is not compact.  However, we will be able to use the mass function $m(V)$
to combat this problem and prove that the minimizers always exist and lie
inside a bounded domain for each $V$.

First we will prove existence of $\Sigma(V)$ for $0 \le V \le V_{MAX}$,
but the 
approach will work for all nonnegative $V_{MAX}$.  Now consider 
$M^3 \cup S^3$, where the union is a disjoint union and $S^3$ is a 
constant curvature 3-sphere with total volume $V_S \gg V_{MAX}$.
The approach will be to prove existence of an area minimizer on this 
manifold, $M^3 \cup S^3$, for volumes less than or equal to $V_{MAX}$,
and then to use the mass function $m(V)$ to 
prove that the minimizers actually contain zero volume in the $S^3$ if
we choose $V_S$ to be large enough.

For the moment, let us redefine $A(V)$ to be exactly as before in 
definition \ref{defav}, 
except that we replace $M^3$ with $M^3 \cup S^3$ 
and $\tilde{M}^3$ with $\tilde{M}^3 \cup S^3$.
Since $A(V)$ is the infimum of the areas of surfaces which contain
a volume $V$ outside the horizon, there exists 
a sequence of surfaces $\{ \Sigma_i \}$ in 
$M^3 \cup S^3$, each containing a volume $V$, and whose areas approach $A(V)$
from above.  Again, since $M^3 \cup S^3$ is not compact, we can not
conclude that the sequence converges to a limit surface with area $A(V)$.
However, using the two propositions below, we will be able to modify this 
sequence of surfaces so that the areas still converge to $A(V)$ and 
each surface stays inside a compact region.  Then since the space of 
rectifiable currents with bounded mass 
in a compact region is compact, we will get a limit
surface in $M^3 \cup S^3$ with area $A(V)$.

The first proposition uses the fact that the Schwarzschild metric becomes very flat 
as we move out to infinity.  
Recall that the Schwarzschild metric of mass $m$ is
$(\real^3-\{0\}, h)$, where $h_{ij} = (1+\frac{m}{2r})^4 \delta_{ij}$ and 
$r$ is the radial coordinate in $\real^3$.  In the next 
proposition, we allow
$m$ to be positive, zero, or negative.  If $m$ is negative, then the
Schwarzschild metric has a singularity at $r = - \frac{m}{2}$.  If $m$
is positive, then the metric has a horizon at $r = \frac{m}{2}$.  In these
cases, for the purposes of the two propositions below, we will say that a surface
contains a volume $V$ when it contains a volume $V$ outside the horizon
or singularity.

\begin{proposition} \label{proposition1}

Consider the Schwarzschild metric $(\real^3-\{0\},h)$ of mass $m$
disjoint union a constant curvature 3-sphere with volume $V_S$.  Then there 
exists an $r$ such that if we choose any $r_1 > r$ and let $r_2 = 2r_1$, then 
if $\Sigma$ is any connected surface containing a volume $V \le V_S$ 
intersecting both the coordinate
sphere of radius $r_1$ and the coordinate sphere of radius $r_2$
(using $\real^3$ coordinates here), then we can modify $\Sigma$ outside of
the coordinate ball of radius $r_1$ to be three surfaces
$\Sigma_1$, $\Sigma_2$, and $\Sigma_3$, with $\Sigma_1$ and $\Sigma_2$ 
in the closed 3-dimensional region contained by $\Sigma$ and with
$\Sigma_3$ in the constant curvature sphere, such that $\Sigma_1$ intersects
the coordinate sphere of radius $r_1$ 
but not the coordinate sphere of radius $r_2$, $\Sigma_2$
intersects the coordinate sphere of radius $r_2$ 
but not the coordinate sphere of radius $r_1$, 
and $\Sigma_1 \cup \Sigma_2 \cup \Sigma_3$ has less area than $\Sigma$
but still contains the same volume $V$. 

\end{proposition}

The main idea of this proposition is that if $r_1$ and $r_2$ are large enough,
then any connected surface with a finite volume $V$ intersecting both spheres
must have at least one very long tentacle.  Since $\Sigma$ has finite total
volume, these tentacles must get very thin.  Then we can snip the tentacles
somewhere in the region between the two spheres so that we get two surfaces,
$\Sigma_1$ and $\Sigma_2$, with $\Sigma_1$ entirely inside the 
coordinate ball of 
radius $r_2$ and $\Sigma_2$ entirely outside the coordinate 
ball of radius $r_1$.
The simplest snipping process would simply be to remove a section of the 
tentacle.  By doing this, we've decreased the volume by $\Delta V$, so we
define $\Sigma_3$ to be a constant curvature 2-sphere of volume $\Delta V$
in the constant curvature 3-sphere.  Thus the total enclosed volume stays
the same, and if we snip the tentacle correctly where we remove a sufficiently
long and skinny section, the total area will decrease.
Proposition \ref{proposition1} follows as a generalization of theorem
\ref{cutting} which is proved in section \ref{another}.
We leave the
details to the reader.

\begin{proposition} \label{proposition2}

Consider the Schwarzschild metric $(\real^3-\{0\},h)$ of mass $m \ge 0$
disjoint union a constant curvature 3-sphere with volume $V_S$.  Then there 
exists an $\tilde{r}$ such that 
if $\Sigma$ is any surface bounding a region of volume $V \le V_S$
entirely outside the 
coordinate ball of radius $\tilde{r}$ , 
then the area of $\Sigma$ is
greater than the area of a constant curvature 2-sphere containing a volume
$V$ in the constant curvature 3-sphere.

\end{proposition}

This proposition follows from the fact that Schwarzschild is very nearly flat
outside a coordinate ball of large radius.  Thus
we get that surfaces nearly satisfy the 
isoperimetric inequality for surfaces in $\real^3$, that 
$A^\frac32 \ge \sqrt{36\pi} V$, whereas constant curvature 2-spheres
in any constant curvature 3-sphere always have $A^\frac32 < \sqrt{36\pi} V$.
Making this idea rigorous is delicate, particularly for small volumes.
We prove proposition \ref{proposition2} later in section 
\ref{another}.
                                     
Now we are ready to prove existence of $\Sigma(V)$ on $M^3 \cup S^3$.
Again, since $A(V)$ is the infimum of the areas of surfaces which contain
a volume $V$ outside the horizon, there exists 
a sequence of surfaces $\{ \Sigma_i \}$ in 
$M^3 \cup S^3$, each containing a volume $V$, and whose areas approach $A(V)$
from above.  Note that if we modify the surfaces in the sequence in a way 
which preserves their enclosed volumes but decreases their areas, then the
areas of the surfaces still approaches $A(V)$.

The first modification we will make to each surface in the sequence is to take
whatever volume is in $S^3$ and to change that part of the surface to be a
single constant curvature 2-sphere in the $S^3$ enclosing that volume.  This
always decreases area since it is known that 2-spheres minimize area with
a volume constraint in $S^3$.  We will repeat this step whenever more 
volume is sent to $S^3$ from $M^3$ in subsequent modifications of the surfaces.

Next we use propositions \ref{proposition1} and \ref{proposition2} to modify 
each surface in the sequence.  By assumption, $M^3$ is isometric to the
Schwarzschild metric outside a compact set.  Again, we use the standard
coordinate chart for the region of $M^3$ which is Schwarzschild, just as
we did in propositions \ref{proposition1} and \ref{proposition2}.  Now we choose $r_1$ to be 
greater than the $r$ of proposition \ref{proposition1} and the $\tilde{r}$ of 
proposition \ref{proposition2} and large enough that $M^3$ is Schwarzschild outside the
coordinate sphere of radius $r_1$.  As in proposition \ref{proposition1}, $r_2 = 2r_1$.
By proposition \ref{proposition1}, we can modify each surface in the sequence so that 
each component of each surface is either entirely inside the coordinate
ball of radius $r_2$ or entirely outside the coordinate ball or radius $r_1$.
A portion of the volume gets sent to $S^3$, but the total volume of the 
surfaces stays the same and the total area decreases.

Next, using proposition \ref{proposition2}, we take all of the components of the surfaces
outside the coordinate ball of radius $r_1$ and send them to spheres of the
same volume in $S^3$, one at a time.  
By proposition \ref{proposition2}, this also decreases
the areas and preserves the volumes of the surfaces in the sequence.  We
send these components to $S^3$ one at a time in the sense that if at any
point there are two spheres (or any other surface that is not a single 
constant curvature sphere) in $S^3$, we immediately turn this portion of the
surface into one constant curvature sphere in $S^3$ with the same volume.  
This always decreases the area, preserves volume, and guarantees that there
will be room for more spheres to be sent to $S^3$.

But now every surface in the sequence is contained in
the coordinate ball of radius $r_2$ union $S^3$.  Since the sequence
of surfaces in now contained in a compact set and the areas of the surfaces
still converge to $A(V)$ from above, it follows from the compactness of the
space of rectifiable currents in a compact manifold that a limit surface
$\Sigma(V)$ exists and that $\mbox{Area}(\Sigma(V)) = A(V)$.

\begin{theorem}

Suppose $(M^3,g)$ is complete, 
has nonnegative scalar curvature, contains a single outermost
minimal sphere $\Sigma_0$, is Schwarzschild
at infinity, and satisfies condition \ref{R3}.
Let $\tilde{M}^3$ be the closure of the component of $M^3 - \Sigma_0$
that contains the asymptotically flat end, and 
let $S^3$ be a constant curvature sphere of volume $V_S$.  
Define
\[A(V) = \inf_\Sigma \{\mbox{Area}(\Sigma)\,\,| \,\, 
\Sigma \mbox{ contains a volume $V$ outside } \Sigma_0    \}\]
where $\Sigma$ is the boundary of some 3-dimensional region in 
$M^3 \cup S^3$ and $\Sigma$ is a surface in $\tilde{M}^3 \cup S^3$ in the same
homology class of $\tilde{M}^3 \cup S^3$ as the horizon $\Sigma_0$.  

Then for all $V \in [0, V_S]$, there exists a surface $\Sigma(V)$
containing a volume $V$ outside $\Sigma_0$
in the same class of surfaces just described such that 
$\mbox{Area}(\Sigma(V)) = A(V)$.

\end{theorem}

Now we will prove existence of $\Sigma(V)$ on $M^3$ for $0 \le V \le V_{MAX}$,
for any nonnegative $V_{MAX}$.  Again, consider 
$M^3 \cup S^3$, where the union is a disjoint union and $S^3$ is a 
constant curvature 3-sphere with total volume $V_S$ much bigger 
than $V_{MAX}$.
We will describe how much bigger in a moment.  Since we have already proven
existence of a minimizer $\Sigma(V)$ on $M^3 \cup S^3$ for volumes up to 
$V_S$, we certainly have existence on the same manifold 
up to the volume $V_{MAX}$.  
Furthermore, by condition \ref{R3}, 
we can choose $\Sigma(V)$ to have at most two
components, with only one component in $M^3$ and possibly one component
in $S^3$.
 
If we choose $V_S$ to be large enough, we can use the mass function $m(V)$ to 
prove that the minimizers actually contain zero volume in the $S^3$,
and hence are entirely contained in $M^3$.
Let $\tilde{V}$ be the supremum of all volumes $\bar{V} \le V_{MAX}$
with the property that $\Sigma(V)$ has zero volume in the $S^3$ for 
$0 \le V \le \bar{V}$.  Since we are assuming $M^3$ satisfies condition  
\ref{R3}, $M^3$ has exactly one horizon, and $\Sigma(0)$ is this horizon,
which of course is contained entirely in $M^3$.
Hence, $\tilde{V} \ge 0$.  

Furthermore, for $0 \le V \le \tilde{V}$, $\Sigma(V)$ has zero volume in
the $S^3$ and hence is in $M^3$ and has only one component.  Thus, by
lemma \ref{AJ}, $m(V)$ is a nondecreasing function of $V$ in this range, 
and since 
$m(0) = \sqrt{\frac{A}{16\pi}}$, where A is the area of the horizon, 
$m(V)$ is positive for $0 \le V \le \tilde{V}$.

For $V \ge \tilde{V}$, $m(V)$ is no longer necessarily nondecreasing.
However, if we reexamine the proof of lemma \ref{AJ} and the derivation
of inequality \ref{R8}, it turns out that there exists a uniform 
$\epsilon > 0$ which is only a function of $V_{MAX}$ and the area of the
horizon such that $m(V) \ge \epsilon$ for $0 \le V \le \tilde{V}+\epsilon$. 

The reason for this is that in this range, $m'(V)$ can be bounded below
uniformly in terms of $V_{MAX}$ and the area of the horizon.  
Inequality \ref{AD} is changed where the $4\pi$ is replaced by an $8\pi$
because the minimizers on which we do a unit normal variation in 
section \ref{s3} may now have up to two components, so the Euler 
characteristic may be as large as 4.  The function $A(V)$ is bounded on
both sides since it is larger than the area of the horizon and smaller than
the area of the horizon plus
$(36\pi)^\frac13 V_{MAX}^\frac23$.  The upper bound on $A(V)$ comes from
comparing $\Sigma(V)$ with a surface which is the horizon union a 
roughly spherical surface containing a volume $V$ very far out on the 
asymptotically flat end of $M^3$.
Finally, since the horizon is 
outermost, $A'(V)$ is bounded below by zero and bounded above for 
$\epsilon$ small enough since we have an upper bound on $A''(V)$ from 
inequality \ref{AD}.  We leave the details of this to the interested
reader.

Since the mass function $m(V) \ge \epsilon$, it follows from 
definition \ref{def1} that $F'(V) = \frac32 A(V)^\frac12 A'(V) 
\le \sqrt{36\pi} - \epsilon'$ for some $\epsilon' > 0$, which is equivalent
to $A'(V) \le \sqrt{\frac{16\pi}{A(V)}} - \epsilon''$ for some uniform
$\epsilon'' > 0$.  On the other hand,  
consider $\Sigma(V)$ where $0 \le V \le \tilde{V} + \epsilon$.
The surface $\Sigma(V)$ has constant mean curvature $H(V)$ on all the
components, and by looking
at a unit normal variation of $\Sigma(V)$ and comparing the 
areas of the variation surfaces
with $A(V)$, we get that 
\begin{equation}\label{hbound}
H(V) \le \sqrt{\frac{16\pi}{A(V)}} - \epsilon''
\end{equation}
where $\epsilon'' > 0$ is a function of $V_{MAX}$ and the area of the 
horizon only.

The mean curvature H of a constant curvature sphere of area A in $\real^3$ 
is $\sqrt{\frac{16\pi}{A}}$.  Furthermore, the mean curvature H of a 
constant curvature sphere of area A in a constant curvature 3-sphere $S^3$
of volume $V_S$ is as close to the $\real^3$ value as we like if we choose
$V_S$ to be large enough.  Now suppose $\Sigma(V)$, with 
$0 \le V \le \tilde{V} + \epsilon$, had a component in $S^3$.
This component is a constant curvature sphere, and we can define $V_S$
in terms of $V_{MAX}$ and the area of the horizon to be large enough so
that the mean curvature of the sphere cannot satisfy inequality \ref{hbound}.
Hence, we have a contradiction, so $\Sigma(V)$ is in $M^3$ for
$0 \le V \le \tilde{V} + \epsilon$.

But $\tilde{V}$ is the supremum of all volumes $\bar{V} \le V_{MAX}$
with the property that $\Sigma(V)$ is entirely contained in $M^3$ for 
$0 \le V \le \bar{V}$.  Hence, since $\epsilon$ is a function
of $V_{MAX}$ and the area of the horizon only, $\tilde{V} = V_{MAX}$, 
proving that the minimizer $\Sigma(V)$ exists in $M^3$ for all $V \ge 0$
since $V_{MAX}$ was arbitrary.

\begin{theorem}

Suppose $(M^3,g)$ is complete, 
has nonnegative scalar curvature, contains a single outermost
minimal sphere $\Sigma_0$, is Schwarzschild
at infinity, and satisfies condition \ref{R3}.
Let $\tilde{M}^3$ be the closure of the component of $M^3 - \Sigma_0$
that contains the asymptotically flat end.
Define
\[A(V) = \inf_\Sigma \{\mbox{Area}(\Sigma)\,\,| \,\, 
\Sigma \mbox{ contains a volume $V$ outside } \Sigma_0    \}\]
where $\Sigma$ is the boundary of some 3-dimensional region in 
$M^3$ and $\Sigma$ is a surface in $\tilde{M}^3$ in the same
homology class of $\tilde{M}^3$ as the horizon $\Sigma_0$.  

Then for all $V \ge 0$ there exists a surface $\Sigma(V)$
containing a volume $V$ outside $\Sigma_0$
in the same class of surfaces just described such that 
$\mbox{Area}(\Sigma(V)) = A(V)$.

\end{theorem}

\section{Existence of Surfaces which Minimize $F$ Given a Volume Constraint}
\label{s8}

We go back to the case 
that $M^3$ has any number of horizons
and assume the hypotheses of theorem \ref{R6}, 
including condition \ref{R5}.  Let $\tilde{M}^3$ be
the closure of 
the component of $M^3 - \{\mbox{the horizons}\}$ that contains the 
asymptotically flat end.  Let
\[ F(V) = \inf_{\{\Sigma_i\}} \{ \sum_{i} \mbox{Area}(\Sigma_i)^\frac32 
\,\,|\,\,
\{\Sigma_i\} \mbox{ contain a volume $V$ outside the horizons} \} \]
where the $\{\Sigma_i\}$ are the boundaries of the components of some 
3-dimensional open region in $M^3$ and 
$\bigcup_i \Sigma_i$ is in $\tilde{M}^3$ and is in the
homology class of $\tilde{M}^3$ which contains both a large sphere at infinity
and the union of the horizons.  If the collection
$\{\Sigma_i\}$ contains a volume $V$ outside the horizons and 
$\sum_{i} \mbox{Area}(\Sigma_i)^\frac32 = F(V)$, then we say that 
$\{\Sigma_i\}$ minimizes $F$ for the volume $V$. 

In this section we will prove that if $M^3$ satisfies condition \ref{R5}, 
then an $F$-minimizer always exists.  Existence of an $F$-minimizer for
all volumes $V \ge 0$ is necessary to prove theorem \ref{R6} since the 
theorem relied on the fact that we had an increasing mass function $m(V)$.
The proof that the mass function was increasing though relied on doing a 
variation of the $F$-minimizers for each $V \ge 0$.  Thus, it is essential
that their exists a collection of surfaces $\Phi(V) = \{\Sigma_i(V)\}$
which minimize $F$ among collections of surfaces in the correct homology
class containing a volume $V$ outside the horizon.  

We will abuse notation slightly again and define 
\[ F(\Phi) = \sum_i \mbox{Area}(\Sigma_i)^\frac32 \]
where $\Phi = \{ \Sigma_i \}$ is any collection of surfaces in $M^3$
which are the boundaries of the components of some 3-dimensional open
region in $M^3$.

First we will prove existence of $\Phi(V)$ for $0 \le V \le V_{MAX}$,
but the 
approach will work for all nonnegative $V_{MAX}$.  Now consider 
$M^3 \cup S^3$, where the union is a disjoint union and $S^3$ is a 
constant curvature 3-sphere with total volume $V_S \gg V_{MAX}$.
The approach will be to prove existence of an $F$-minimizer on this 
manifold, $M^3 \cup S^3$, for volumes less than or equal to $V_{MAX}$,
and then to use the mass function $m(V)$ to 
prove that the minimizers actually contain zero volume in the $S^3$ if
we choose $V_S$ to be large enough.

For the moment, let us redefine $F(V)$ to be exactly as above 
except that we replace $M^3$ with $M^3 \cup S^3$ 
and $\tilde{M}^3$ with $\tilde{M}^3 \cup S^3$.
Since $F(V)$ is the infimum of the $F$-values of collections of surfaces 
which contain
a volume $V$ outside the horizon, there exists 
a sequence of collections of surfaces $\{ \Phi_i \}$ in 
$M^3 \cup S^3$, each containing a volume $V$, and whose $F$-values 
approach $F(V)$
from above.  

In the previous section, we used propositions \ref{proposition1} and \ref{proposition2}
to show that we could modify any  
sequence of surfaces in 
$M^3 \cup S^3$ to lie inside a compact region of $M^3 \cup S^3$ without
increasing the areas of any of the surfaces.  The technique, though, did
increase the number of the components of the surfaces.  However, since the
total area went down and the number of components went up, it follows
from the fact that $(a+b)^\frac32 > a^\frac32 + b^\frac32$ 
for $a$ and $b$ positive that these
same techniques can be used to modify a sequence of collections of surfaces 
$\{\Phi_i\}$ so that the sequence lies inside a compact region of 
$M^3 \cup S^3$ without increasing the $F$-values of any of the collections
of surfaces.  Also, since it can be checked by direct calculation that the
collection of surfaces which minimizes $F$ inside $S^3$ is a single spherically
symmetric sphere, each of the new modified collections of surfaces will have at
most one component in the $S^3$ which will always be a spherically symmetric
sphere as before.
Since the $F$-values are not increased, the $F$-values of 
the new modified sequence $\{\Phi_i\}$ still converge to $F(V)$ from above.

In the introduction to this chapter 
we commented that there were two problem to look out
for in the existence of $F$-minimizers.  The first is that a component of 
the $F$-minimizer could run off to infinity.  This problem is taken care 
of since we are able to require our minimizing sequence to stay inside a 
compact set.  The other problem with $F$-minimization, though, is that 
``bubbling'' might occur, where the optimal configuration is an infinite 
number of tiny balls with a finite total volume.  To combat this, we 
modify the sequence $\{\Phi_i\}$ one last time.  We know bubbling cannot
happen in the $S^3$, since, by direct calculation, the $F$-minimizers in 
$S^3$ are single spherically symmetric spheres.  Hence, we modify a given
collection of surfaces $\{\Sigma_i\}$ using the following rule which we
will call the ``sphere replacement rule''.  If any 
subcollection of the surfaces in $M^3$ would have smaller $F$-value by 
replacing them with a single sphere in $S^3$, then we make the replacement.
This rule puts an upper bound on how much volume can be used for tiny balls 
in $M^3$ since at some point $F$ can be reduced by replacing a large 
number of tiny balls in $M^3$ by a single sphere with the same volume in 
$S^3$.  The reason for this is that since $M^3$ is smooth, it has bounded
curvature and hence on the small scale is approximately flat.  Since $F$ 
scales like volume, it follows that a bunch of tiny balls containing a 
volume $V$ will have $F$-value close to $\sqrt{36\pi}V$, which by direct
calculation is larger than the $F$-value of a single sphere 
in $S^3$ containing the same volume.  And as before, if at any time there 
are two or more spheres in $S^3$, then we combine them into one sphere 
containing the same volume and this also always decreases the $F$-value.

Now we are ready to take a limit of a subsequence of $\{\Phi_i\}$.  This
is a little tricky since each $\Phi_i$ is not a surface, but a collection
of connected surfaces.  
For each $i$, order the surfaces of each collection $\Phi_i$
by the volume outside the outermost horizons enclosed by 
each surface, with largest volumes first.  If two surfaces
enclose the same volume, then choose either ordering.  By the Federer-Fleming
compactness theorem, there exists a subsequence $\{\Phi_{1,i}\}$
of $\{\Phi_i\}$ in which
the largest surfaces of each $\Phi$ converge to a limit.
Similarly, there exists a subsequence $\{\Phi_{2,i}\}$ of $\{\Phi_{1,i}\}$
in which the second largest surfaces of each $\Phi$ converge to a limit. 
Repeating this process we define the sequence $\{\Phi_{n,i}\}$ for $n\ge1$.
Finally, we define $\{\tilde{\Phi}_i\} = \{\Phi_{i,i}\}$, which has the 
property that the largest surfaces converge to a limit, the second largest
surfaces converge to a limit, and so on, and we define the collection of 
the limit surfaces to be $\tilde{\Phi}$.

While we do get a collection of limit surfaces $\tilde{\Phi}$, we still need
to show that they enclose the correct volume $V$ that each collection of
surfaces in the original sequence enclosed.  Suppose $\tilde{\Phi}$ did not
enclose the volume $V$ but instead only enclosed a volume $V-v$ for some
$v > 0$.  Then for any $\epsilon>0$ there must exist an $i$ such that in 
the collection $\Phi_i$ there is a large subcollection of tiny surfaces each 
containing less than $\epsilon$ volume each but containing a total volume of 
$v$.  In other words, bubbling has occurred.  But by the sphere replacement
rule, this can not happen, since the $F$-value of $\Phi_i$ would have been
reduced by replacing the large subcollection of tiny surfaces by a single
sphere in $S^3$ containing a volume $v$, for some value of $\epsilon>0$.  
Hence, $\tilde{\Phi}$ encloses a volume $V$, and $F(\tilde{\Phi}) = F(V)$.

\begin{theorem}

Suppose $(M^3,g)$ is complete, 
has nonnegative scalar curvature, contains any number of outermost
minimal spheres $\{\tilde{\Sigma}_i\}$, is Schwarzschild
at infinity, and satisfies condition \ref{R5}.
Let $\tilde{M}^3$ be the closure of the component of $M^3 - 
\{\tilde{\Sigma_i}\}$
that contains the asymptotically flat end, and 
let $S^3$ be a constant curvature sphere of volume $V_S$.  
Define
\[ F(V) = \inf_{\{\Sigma_i\}} \{ \sum_{i} \mbox{Area}(\Sigma_i)^\frac32 
\,\,|\,\,
\{\Sigma_i\} \mbox{ contain a volume $V$ outside the horizons} \} \]
where the $\{\Sigma_i\}$ are the boundaries of the components of some 
3-dimensional open region in $M^3 \cup S^3$ and 
$\bigcup_i \Sigma_i$ is in $\tilde{M}^3 \cup S^3$ and is in the
homology class of $\tilde{M}^3 \cup S^3$ 
which contains both a large sphere at infinity
and the union of the horizons.

Then for all $V \in [0, V_S]$, there exists a collection of surfaces 
$\Phi(V) = \{\Sigma_i(V)\}$
containing a volume $V$ outside the horizons
in the same class of surfaces just described such that 
$F(\Phi(V)) = F(V)$.

\end{theorem}

Now we will prove existence of $\Phi(V)$ on $M^3$ for $0 \le V \le V_{MAX}$,
for any nonnegative $V_{MAX}$.  Again, consider 
$M^3 \cup S^3$, where the union is a disjoint union and $S^3$ is a 
constant curvature 3-sphere with total volume $V_S$ much bigger 
than $V_{MAX}$.
We will describe how much bigger in a moment.  Since we have already proven
existence of a minimizer $\Phi(V)$ on $M^3 \cup S^3$ for volumes up to 
$V_S$, we certainly have existence on the same manifold 
up to the volume $V_{MAX}$.  
Furthermore, by condition \ref{R5}, 
we can choose the collection of surfaces $\Phi(V)$ so that no two of its 
surfaces touch.  Hence, we can perform unit normal variations on each surface
of the collection, so the mass function $m(V)$ is nondecreasing as long as 
we require $V_S \ge 2V_{MAX}$ (since we need the mean curvature of the 
minimizers to be positive to get $F'(V)$ nonnegative which is required 
for nondecreasing mass).   

Since the mass function is initially positive since 
$m(0) = \left(\sum_{i=1}^n\left(\frac{A_i}{16\pi}\right)^\frac32\right)
^\frac13$, 
and since $m$ is nondecreasing, $m(V)$ is always positive for 
$0 \le V \le V_{MAX}$.  Since 
$m(V) = F(V)^\frac13 (36\pi - F'(V)^2)/c$
where $c=144\pi^\frac32$, it follows that $F'(V) < \sqrt{36\pi} - \epsilon$
for $0 \le V \le V_{MAX}$ for some $\epsilon > 0$.

On the other hand, suppose the surfaces $\{\Sigma_i\}$ minimize
$F$ while enclosing a volume $V$.  It follows from the first variation
of area on each surface that each surface has constant (generally distinct)
mean curvature $H_i$.  Furthermore, from this same first variational 
computation it follows that if we consider any smooth variation on these 
surfaces, the rate of change of $F$ with respect to $V$ will be
\[ \frac{dF}{dV} = \frac32 A_i^\frac12 H_i \] 
for all $i$.  By comparing this variation with other minimizers, it follows
that 
\begin{equation}\label{I1}
\frac32 A_i^\frac12 H_i \le \mbox{the left sided derivative of } F(V)
\le \sqrt{36\pi} - \epsilon
\end{equation}
for some fixed $\epsilon > 0$.

But if we choose $V_S$ to be large enough, then the local 
geometry of the sphere 
$S^3$ can be made as close to that of $\real^3$ as we like.  Hence, for
a sphere containing a volume less than $V_{MAX}$ in $S^3$,  
$\frac32 A^\frac12 H$ can be made as close to $\frac32 \sqrt{4\pi r^2}
\frac{2}{r} = \sqrt{36\pi}$ as we like if we choose $V_S$ large enough, 
violating inequality \ref{I1}.  Hence, if we choose $V_S$ large enough, 
then the minimizer $\Phi(V) = \{\Sigma_i(V)\}$ will not have
any components in the $S^3$, which proves that $\Phi(V)$ minimizes $F$
in $M^3$ among all other collections of surfaces in $M^3$ in the correct
homology class containing the same volume $V$, for $0\le V \le V_{MAX}$.
But since $V_{MAX}$ was arbitrary, we have a $F$-minimizer for all $V \ge 0$.

\begin{theorem}

Suppose $(M^3,g)$ is complete, 
has nonnegative scalar curvature, contains any number of outermost
minimal spheres $\{\tilde{\Sigma}_{i}\}$, is Schwarzschild
at infinity, and satisfies condition \ref{R5}.
Let $\tilde{M}^3$ be the closure of the component of $M^3 - 
\{\tilde{\Sigma}_{i}\}$
that contains the asymptotically flat end.
Define
\[ F(V) = \inf_{\{\Sigma_i\}} \{ \sum_{i} \mbox{Area}(\Sigma_i)^\frac32 
\,\,|\,\,
\{\Sigma_i\} \mbox{ contain a volume $V$ outside the horizons} \} \]
where the $\{\Sigma_i\}$ are the boundaries of the components of some 
3-dimensional open region in $M^3$ and 
$\bigcup_i \Sigma_i$ is in $\tilde{M}^3$ and is in the
homology class of $\tilde{M}^3$ 
which contains both a large sphere at infinity
and the union of the horizons.

Then for all $V \ge 0$, there exists a collection of surfaces 
$\Phi(V) = \{\Sigma_i(V)\}$
containing a volume $V$ outside the horizons
in the same class of surfaces just described such that 
$F(\Phi(V)) = F(V)$.

\end{theorem}

\section{Another Isoperimetric Inequality for the \\
         Schwarzschild Metric}                   \label{another}

In section \ref{s6.5}, we proved that the spherically symmetric spheres of
the Schwarzschild metric minimize area among all surfaces in their homology
class containing the same volume outside the horizon.  This gives a lower
bound for the area of any surface containing the horizon in terms of the 
volume outside the horizon that the surface encloses, and so is an 
isoperimetric inequality.

In this section, we lead up to proving proposition \ref{proposition2} 
of section 
\ref{s7} which is an isoperimetric inequality for the asymptotically flat
portion of the Schwarzschild metric since it gives lower bounds for the 
areas of surfaces in terms of their enclosed volumes.  We also prove
theorem \ref{cutting} below, which is necessary in the proof of proposition
\ref{proposition2}, and which, when generalized sufficiently, proves
proposition \ref{proposition1} of section \ref{s7} as well.  We begin 
with two definitions.

\begin{definition}
Let $D^3$ be a region in $\real^3$ and let $\Sigma^2=\partial D^3$ which we
assume is smooth.
Let
\begin{eqnarray*}
D^3(x_1,x_2)&=&\{(x,y,z)\in D^3 | x_1<x<x_2\}\\
\Sigma^2(x_1,x_2)&=&\{(x,y,z)\in \Sigma^2 | x_1<x<x_2\}\\
C^2(x_1)&=&\{(x,y,z)\in D^3 | x=x_1\}
\end{eqnarray*}
\end{definition}

\begin{definition}  Define $\Delta A(x_1,x_2)$ by
\[\Delta A(x_1,x_2)=|C^2(x_1)|+|C^2(x_2)|+(36\pi)^{1/3}|D^3(x_1,x_2)|^{2/3}
-|\Sigma^2(x_1,x_2)|\]
\end{definition}
Hence, $\Delta A$ is the change in the surface area of $D^3$ if we cut
a section of $D^3$ out from $x=x_1$ to $x=x_2$ and replace this
section with a ball of equal volume.

\begin{theorem}\label{cutting}
There exist $\alpha, \beta>0$ such that if for some $d>0$
\begin{enumerate}
\item $\inf_{0\le x\le d} |C^2(x)| > 0$

\item $\frac{|\Sigma^2(0,d)|}{d^2} < \alpha$,
\end{enumerate}

then there exists $x_1, x_2\in [0,d], x_1 < x_2,,$
such that
\[\frac{\Delta A(x_1,x_2)}{|\Sigma^2(x_1,x_2)|}<-\beta.\]
\end{theorem}

{\em Proof.}
Since the theorem is scale-invariant, we may as well assume
$d=4$.  Now we break the interval $[0, 4]$ into six intervals, the two at
the ends having length $2\epsilon$ and the four in the middle having
length $1-\epsilon$, for some positive $\epsilon \ll 1$.  Let
\begin{eqnarray*}
I_1&=&[0,2\epsilon]\\
I_2&=&[2\epsilon,1+\epsilon]\\
I_3&=&[1+\epsilon,2]\\
I_4&=&[2,3-\epsilon]\\
I_5&=&[3-\epsilon,4-2\epsilon]\\
I_6&=&[4-2\epsilon,4]
\end{eqnarray*}
be these six intervals.  Abusing notation slightly, let
\[p_k=\frac{|\Sigma^2(I_k)|}{|\Sigma^2(0,4)|}, \,\,\,\,\,\, 1\le k\le 6,\]
so that $\sum_{k=1}^6 p_k = 1$.  Hence, $p_k$ is the fraction of the area of
$\Sigma^2(0,4)$ which is in the region $I_k\times \real^2$.

Note that, in general, if
\begin{equation}\label{star1}
 \frac{|\Sigma^2(x_1,x_2)|}{|\Sigma^2(0,4)|} <
 \left(\frac{x_2-x_1}{4}\right)^2,
\end{equation}
then we can iterate this proof by substituting the interval $(x_1, x_2)$ for
$(0, 4)$ then rescaling $\real^3$ (in all directions)
so that $(x_1, x_2)$ becomes $(0,4)$.
Hypothesis 1 of the theorem is still satisfied, and hypothesis 2 is still
satisfied too since areas scale as the square of distances.  Hence, since
the conclusion of the new rescaled theorem is stronger, the original theorem
follows from the rescaled theorem.

We choose to rescale if inequality \ref{star1} is satisfied for $[x_1, x_2] =
I_2, I_3, I_4, I_5, \cup_{k=1}^5 I_k, $ or $\cup_{k=2}^6 I_k$.
We claim this iteration process can only happen a finite number of times.
First we note that
\[|\Sigma^2(0,4)|\ge \int_0^4 [\mbox{length of } \partial C^2(x)]\,dx
\ge \int_0^4\sqrt{4\pi|C^2(x)|}\, dx\ge 4\sqrt{4\pi a}\]
where we let $a=\inf_{0\le x\le 4} |C^2(x)|$.  Hence
\begin{equation}\label{star2}
 a \le \frac{|\Sigma^2(0,4)|^2}{64\pi}.
\end{equation}
Each time we rescale, it follows from inequality \ref{star1} 
that the area of $\Sigma^2(0,4)$
does not increase.  However, $a$ goes up by at least a factor of
$\frac{4}{4-2\epsilon}$.  Hence, since in the theorem we assumed
$a>0$, we must only rescale a finite number of times or 
inequality \ref{star2} would
be violated.

In the final rescaled interval, we must therefore have
\[\frac{|\Sigma^2(x_1,x_2)|}{|\Sigma^2(0,4)|}\ge \left(\frac{x_2-x_1}{4}
\right)^2\]
for $[x_1,x_2] =I_2, I_3, I_4, I_5, \cup_{k=1}^5 I_k, $and $\cup_{k=2}^6 I_k.$
Thus,
\[p_2, p_3, p_4, p_5 \ge \left(\frac{1-\epsilon}{4}\right)^2\]
and
\[\sum_{k=1}^5 p_k, \sum_{k=2}^6 p_k\ge\left(\frac{4-2\epsilon}{4}\right)^2.\]
Since $\sum_{k=1}^6 p_k = 1$, it follows (but it is not equivalent to) that

\begin{eqnarray}\label{star3}
p_1, p_6 & \le & \epsilon\\ \label{star4}
p_2, p_3, p_4, p_5 & \ge & \frac{1}{16}(1-2\epsilon).
\end{eqnarray}

To get an upper bound on $\frac{\Delta A(x_1,x_2)}{|\Sigma^2(x_1,x_2)|}$ in the
conclusion of the theorem, we choose $x_1=0$, $x_2=4$ and determine the
region $D^3$ which maximizes
\begin{equation}\label{maximize}
\frac{\Delta A(0,4)}{|\Sigma^2(0,4)|}
\end{equation}
while still satisfying \ref{star3} and \ref{star4}.

This optimal region $D^3$ must be
axially symmetric around the $x$-axis.  This follows from the
following symmetrization argument.  Given a region $D^3$, 
symmetrize it about the $x$-axis by defining another region
$D_{SYM}^3$ to be axially symmetric around the $x$-axis but 
having
the same cross sectional area as $D^3$ when intersected by planes
given by $x$ equal to a constant.  $D_{SYM}^3$ and $D^3$ have
the same volume, and it is known that this symmetrization 
process decreases surface area.  In fact, $D_{SYM}^3$ will have
less surface area than $D^3$ in each region $I_k \times \real^2$.
We want to preserve inequalities \ref{star3} and \ref{star4},
so define $\bar{D}^3$ to be $D_{SYM}^3$ union any regions in
$\real^3$ so that $\bar{D}^3$ has the same area as $D^3$ in each
$I_k \times \real^2$.  Then since $\bar{D}^3$ has more volume
than $D^3$, we see that $D^3$ can only maximize
\ref{maximize} if it is axially symmetric. 

Furthermore, from the first variation
formula, $\Sigma^2=\partial D^3$ must have constant mean curvature in each of
the six intervals.  Hence, in each interval $\Sigma^2$ is either a collection
of spheres or a Delaunay surface.  If $\alpha$ (from the statement of the
theorem) is small enough, we can rule out Delaunay surfaces 
since they are unstable.  
We can also rule out more than one sphere completely contained 
in the interiors of each of the six
intervals using stability since decreasing the area
of one of the spheres while increasing the area of one of the
other spheres at the same rate always increases volume
to second order.

For convenience, let's rescale again so that $|\Sigma^2(0,d)|=1$.  Then
checking all the possibilities we find that one of the optimal regions
$D^3$ which maximizes $\Delta A(0,d)$ is the right portion (with outside
surface area $\epsilon$) of a ball in $I_1\times \real^2$ union
the left portion (with outside surface area $\epsilon$) of a ball in $I_6\times
\real^2$ union a ball with surface area $\frac{1}{8}(1-2\epsilon)$
centered in $(I_2\cup I_3)\times \real^2$ union a ball with surface area
$\frac{7}{8}(1-2\epsilon)$ centered in $(I_4\cup I_5)\times \real^2$.  For
this region $D^3$, we can then estimate that
\[|C^2(0)|, |C^2(d)| \le \epsilon\]
and
\[|D^3(0,d)|\le (36\pi)^{-1/2}\left[(2\epsilon)^{3/2}+(\frac{1}{8}(1-2\epsilon)
)^{3/2}+(\frac{7}{8}(1-2\epsilon))^{3/2}\right]\]
so that
\[\Delta A(0,d)\le 2\epsilon+\left[(2\epsilon)^{3/2}+((\frac{1}{8})^{3/2}
+(\frac{7}{8})^{3/2})(1-2\epsilon)^{3/2}\right]^{2/3}-1\]
where again we recall that we have rescaled so that $|\Sigma^2(0,d)|=1$.
Note that when $\epsilon=0$, the right hand side of the above equation
equals
\[\left[(\frac{1}{8})^{3/2}+(\frac{7}{8})^{3/2}\right]^{2/3}-1<0.\]
Hence, by choosing $\epsilon$ small enough, we have
\[\frac{\Delta A}{|\Sigma^2(0,d)|} < -\beta\]
for some $\beta>0$.  Since this was for the maximal configuration for $D^3$,
the theorem follows.  \qed

We call theorem \ref{cutting} the ``cutting theorem'' since it tells us
that if a region is long and skinny enough, then we can cut out a portion
of it and replace that portion with a ball of equal volume and decrease
the total surface area in the process.  $\Delta A(x_1,x_2)$ is the amount
the area changes when we cut out the section $D^3(x_1,x_2)$, and 
\[ \frac{|\Sigma^2(0,d)|}{d^2} < \alpha \]
is the condition we need to know that $D^3$ is long and skinny enough.

Intuitively, this theorem is clear, but we see that the proof was nontrivial.
We claim, but neglect to prove here, two generalizations of theorem
\ref{cutting}.  First, we will need a cutting theorem like theorem
\ref{cutting} for the Schwarzschild metric outside a coordinate ball
of radius $R$ to prove proposition \ref{proposition2}, for some large $R>0$.
Since the Schwarzschild metric $(\real^3-\{0\},h)$ has conformal factor
$(1+\frac{m}{2r})^4$ which approaches $1$ as $r$ approaches infinity, we can
view the Schwarzschild metric as a perturbation of $\real^3$, with the 
perturbation being as small as we like if we choose $R$ large enough.
Hence, it is reasonable to use rotated and translated versions of the standard
$\real^3-\{0\}$ coordinate chart for the Schwarzschild metric of mass $m$
to define $D^3(x_1,x_2)$, $\Sigma^2(x_1,x_2)$, and $C^2(x_1)$, and then to
state a generalized version of theorem \ref{cutting} for regions $D^3$ in
the Schwarzschild metric entirely outside the coordinate ball of radius $R$, 
for some $R > 0$.  Secondly, we claim that proposition \ref{proposition1}
follows as a further generalization of theorem \ref{cutting}, where not only
are we now in the Schwarzschild metric of mass $m$, but the cuts are being made
along planes parallel to the sides of a large polyhedron
contained inside the coordinate ball of radius $r_2$ minus
the coordinate ball of radius $r_1$.  
In this way proposition \ref{proposition1} follows, 
after sufficient adaptation, from the proof of 
theorem \ref{cutting}.

Now we prove proposition \ref{proposition2} from section \ref{s7}.
Since the Schwarzschild metric is conformal to $\real^3-\{0\}$, with conformal
factor $(1+\frac{m}{2r})^4$, then for surfaces outside the coordinate ball
of radius $\tilde{r}$ we can use the isoperimetric inequality for $\real^3$
to conclude that
\[A^{3/2}\ge\sqrt{36\pi}V(1+\frac{m}{2\tilde{r}})^{-6}.\]

Let $A(V_S,V)$ be the area of a constant curvature 2-sphere containing a
volume $V$ in the constant curvature 3-sphere of volume $V_S$.
Then since $A(V_S,V)^{3/2}<\sqrt{36\pi}V$ for $V>0$, with the inequality being
by a uniform amount for $V\ge \epsilon$ given an $\epsilon>0$, we see that
we can simply choose $\tilde{r}$ large enough to prove 
proposition \ref{proposition2} for $V\ge 
\epsilon$.

To prove proposition \ref{proposition2}
for small $V$, we observe that
\begin{equation}\label{tri}
A(V_S,V)^{3/2} = \sqrt{36\pi} V\left[1-
k\left(\frac{V}{V_S}\right)^{2/3}+O_2[\left(\frac{V}{V_S}\right)^{2/3}]\right]
\end{equation}
where $k=\frac{3}{10}\left(\frac{3\pi}{2}\right)^{2/3}$.  We will show that for
small volumes $V<\epsilon$, if we choose $\tilde{r}$ large enough, then in
$(\real^3-\{0\},h)$ outside the coordinate ball of radius $\tilde{r}$ that
all surfaces of area $A$ containing a volume $V$ satisfy
\begin{equation}\label{one}
A^{3/2}\ge \sqrt{36\pi}V\left[1-\frac{1}{2}k\left(\frac{V}{V_S}
\right)^{2/3}\right]
\end{equation}
which will prove proposition \ref{proposition2} 
for $V<\epsilon$ if we choose $\epsilon$ small enough.
Thus, all that remains is to establish inequality \ref{one} 
for $V<\epsilon$, for some $\epsilon>0$.

Suppose $\Sigma^2=\partial D^3$ contains a volume $V=|D^3|<\epsilon$ and
is entirely outside the coordinate ball of radius $\tilde{r}$ in
Schwarzschild.  We assume $\Sigma^2$ is smooth, but $\Sigma^2$
could have tentacles extending long distances, for example, which is
troublesome.
We find it necessary to regularize $\Sigma^2$ first before proving
inequality \ref{one}.

Let $U^3$ be any open subset of the bounded open set $D^3$.  Define
\[f(U^3)=
\mbox{Area}(\partial(D^3-U^3))+(36\pi)^{1/3}\mbox{Volume($U^3$)}^{2/3}.\]
Since we have uniform bounds on $|U^3|$ and $|\partial U^3|$ 
when $f$
is being minimized since $|\partial U^3|\le|\partial(D^3-U^3)|+|\partial D^3|$,
and since $D^3$ is bounded, there exists a region $U_0\subset D^3$ which
minimizes $f$.  Note that since $f(\emptyset)=|\partial D^3|$,
$f(U_0)\le |\partial D^3|$.

Finally, we regularize $D^3$ by removing the region $U_0$ from $D^3$.  This,
of course, decreases the total volume, so to keep the total volume constant
we add a ball of volume $|U_0|$ to a copy of $\real^3$.  Thus, we've modified
$D^3$ and replaced it with $\bar{D}^3=(D^3-U_0)\cup B^3\subset
\mbox{Schwarzschild}\cup\real^3$, where $B^3$ is the ball of volume
$|U_0|$ in $\real^3$.  Note that $|\bar{D}^3|=|D^3|$ and that the
area has decreased since
\[|\partial \bar{D}^3|=|\partial(D^3-U_0)|+(36\pi)^{1/3}|U_0|^{2/3}
=f(U_0)\le f(\emptyset) = |\partial D^3|.\]
Thus, it is sufficient to prove inequality \ref{one} for the regularized region
$\bar{D}^3$ in \\
Schwarzschild (disjoint) union $\real^3$.

It is also sufficient to prove inequality \ref{one} for each component of 
$\bar{D}^3$ individually.  The ball $B^3$ in $\real^3$ satisfies inequality 
\ref{one} since $A^{3/2}=\sqrt{36\pi} V$ for balls.  Now consider one of the
components $\Sigma^2_i = \partial \bar{D}^3_i$ in Schwarzschild outside
the coordinate ball of radius $\tilde{r}$.  Note that for $\Sigma^2_i$,
\begin{equation}\label{two} 
A\le (36\pi)^{1/3}V^{2/3}
\end{equation}
because otherwise $U_0$ would have included the region $\bar{D}^3_i$.
Furthermore,
\begin{equation}\label{three}
\mbox{diam}(\Sigma^2_i)\le\alpha^{-1/2}(36\pi)^{1/6}V^{1/3},
\end{equation}
where $\alpha$ is the constant from the cutting theorem and $\mbox{diam}
(S)$ is the diameter of $S$.
Otherwise, we would have $\mbox{diam}(\Sigma^2_i)>\alpha^{-1/2}((36\pi)^{1/3}
V^{2/3})^{1/2}\ge\alpha^{-1/2}A^{1/2}$ which means we could use the cutting
theorem to remove a section of $\Sigma^2_i$, form a ball 
in $\real^3$ with it, and
decrease the boundary area while preserving the total volume.  This cannot
happen, since $U_0$ would have included this section of $\Sigma^2_i$ if
removing it and forming a ball with it decreased the total area.
Hence, we must have inequality \ref{three}.  
This diameter bound is central to the rest
of the proof and is the reason we needed to regularize $D^3$.

Pick any point $p_0$ in $\bar{D}^3_i$, where again $\Sigma^2_i
=\partial \bar{D}^3_i$.  In coordinates, Schwarzschild can be represented
as $(\real^3-\{0\},h)$, where $h_{ij}=\left(1+\frac{m}{2r}\right)^4\delta_{ij}$
is the metric and $r$ is the radial coordinate in $\real^3$.  Suppose
$p_0$ has radial coordinate $r_0$.  Since $\Sigma^2_i$ is outside the
coordinate ball of radius $\tilde{r}$, $r_0\ge\tilde{r}$.

We construct a spherically-symmetric mapping $\phi$ from a
spherically-symmetric connected neighborhood of Schwarzschild containing $p_0$ to
a spherically-symmetric connected annular neighborhood of a large 3-sphere $S^3$ of
radius $R_0$ (when embedded in $\real^4$).  We want $\phi:(\real^3-\{0\},h)
\to(S^3,g_0)$ to be spherically-symmetric, locally volume preserving, and
``tangent'' (to be defined in a moment) at $p_0$.

Let
\[u(r)^{-1}=\frac{\|D\phi(\partial_r)\|_{g_0}}{\|\partial_r\|_h}\]
where $\partial_r=\frac{\partial}{\partial_r}$ is a radial tangent vector in
Schwarzschild, $\|\cdot\|_h$ is the length in the Schwarzschild metric, and
$\|\cdot\|_{g_0}$ is the length in the sphere of radius $R_0$ metric.  Thus
$\phi$ increases lengths in the radial direction by a factor of $u(r)^{-1}$.
Since $\phi$ preserves volumes locally, lengths in the two other mutually
orthogonal directions must be increased by a factor of $u(r)^{1/2}$,
so that the areas of the spherically symmetric spheres of 
the Schwarzschild metric 
get increased by a factor of $u(r)$ by $\phi$.  We choose the radius
$R_0$ of $S^3$ and define $\phi$ such that $u(r_0)=1$ and
$\frac{du}{dr}(r_0)=0$, in which case we say $\phi:(\real^3-\{0\},h)
\to (S^3,g_0)$ is tangent at $r=r_0$, and in particular at $p_0$.

Since volume is preserved by $\phi$ locally, it is most convenient to
parameterize the spherically symmetric functions by the enclosed volume
of the corresponding spherically symmetric spheres, or
at least relative enclosed volume. On Schwarzschild, define $v(r)$ to be
the volume enclosed by the coordinate ball of radius $r$ outside the coordinate
ball of radius $r_0$.  When $r<r_0$, $v(r)<0$, and $v(r_0)=0$.  Let $U(v)$
be $u(r)$ changed into $v$-coordinates, and let $A_0(v)$ be the area of
the spherically symmetric sphere in Schwarzschild containing a volume $v$
outside the coordinate ball of radius $r_0$.  Use $\phi$ to define $v$
on $(S^3,g_0)$, and let $A_1(v)$ be the area of the spherically symmetric
spheres in $(S^3,g_0)$.  Since $\phi$ is locally volume preserving,
$v$ is relative enclosed volume on $(S^3,g_0)$ as well as $(\real^3-\{0\},h)$.

In the Schwarzschild metric of mass $m$,
\[\left(\frac{A_0(v)}{16\pi}\right)^{1/2}\left(1-\frac{1}{16\pi}
A_0(v)A_0'(v)^2\right)=m\]
for all $v$.  This follows from the fact that the mean curvature of the
spheres is given by $H=A'(v)$ and the formula for the Hawking mass.
In a 3-sphere $(S^3,g_0)$ of radius $R_0$ (when embedded in $\real^4$),
we compute directly that
\[\frac{4\pi}{A_1(v)}\left(1-\frac{1}{16\pi}A_1(v)A_1'(v)^2\right)=R_0^{-2}.\]
At the point of tangency ($v=0$), $A_0(0)=A_1(0)$ and $A_0'(0)=A_1'(0)$.
Hence, dividing the two previous formulas at $v=0$ gives us that
\[mR_0^2=\frac{1}{2}\left(\frac{A_0(0)}{4\pi}\right)^{3/2}.\]

Hence, the further out $p_0$ is in Schwarzschild, the larger $A_0(0)$ is
and the larger $R_0$ is.  Thus, we may guarantee $R_0$ to be as large as we
like if we choose $\tilde{r}$ large enough.  (Also, we see that this
construction only works when $m>0$, which, by the positive mass theorem,
is all we need.  Proposition \ref{proposition2} 
is true for $m\le 0$, but the proof requires
constructing tangent hyperbolic spaces instead of tangent spheres.)

Furthermore, since $A_1(v)=U(v)A_0(v)$, we can differentiate this twice,
and use $A_0(0)=A_1(0)$ and $A_0'(0)=A_1'(0)$ to get
\[U''(0)=\frac{A_1''(0)-A_0''(0)}{A_0(0)}.\]
Changing the $U(v)$ back to $r$ coordinates, we get
\[u''(r_0)=A_0(0)(A_1''(0)-A_0''(0))\]
where the derivatives on $A_0$ and $A_1$ are still with respect to $v$.
Working out the behavior of $A_0''(0)$ and $A_1''(0)$ for large $A_0(0)$,
we find that $u''(r_0)$ goes down as $k/A_0(0)$, for some constant $k$
for large $A_0(0)$.  Thus, we can guarantee $u''(r_0)$ to be as small as 
we like if we choose $\tilde{r}$ large enough.

Since $V<\epsilon$, $\mbox{diam}(\Sigma^2_i)\le \alpha^{-1/2}(36\pi)^{1/6}
\epsilon^{1/3},$ so we need only to extend $\phi$ this distance both ways
from $p_0$.  Hence, on $\Sigma^2_i$,
\[u(r)\ge 1-\delta(r-r_0)^2\]
for some $\delta>0$ since $u(r_0)=1$ and $u'(r_0)=0$, and
we can choose $\delta$ as small as we want if we choose $\tilde{r}$ large
enough since this will make $u''(r_0)$ small.  Thus, since direct calculation
shows that $u(r)\le 1$, $\phi$ increases areas pointwise by a factor less
than or equal to $u(r)^{-1/2}$, and
\begin{eqnarray*}
u(r)&\ge& 1-\delta\left[\mbox{diam}(\Sigma^2_i)\right]^2\\
&\ge& 1-\delta\alpha^{-1}(36\pi)^{1/3}V^{2/3}.
\end{eqnarray*}

We will use the isoperimetric inequality for $(S^3,g_0)$, that the
spherically symmetric spheres minimize area among surfaces enclosing the
same volume, to prove inequality \ref{one} for $V<\epsilon$.  If we choose
$\tilde{r}$ large enough, the radius $R_0$ (and total volume) of
$(S^3,g_0)$ will be as large as we want, so by inequality \ref{tri}
\[A^{3/2}\ge\sqrt{36\pi}V[1-\delta V^{2/3}] 
\,\,\,\,\mbox{in}\,\,\,\, (S^3,g_0)\]
for $V < \epsilon$.  Since
$\phi:(\real^3-\{0\},h)\to(S^3,g_0)$ increases areas less than
$u(r)^{-1/2}$, then in Schwarzschild we have
\begin{eqnarray*}
A^{3/2}&\ge& u(r)^{3/4}\sqrt{36\pi}V[1-\delta V^{2/3}]\\
&\ge&\sqrt{36\pi}V[1-\delta V^{2/3}][1-\delta\alpha^{-1}
(36\pi)^{1/3}V^{2/3}]^{3/4}
\end{eqnarray*}
for $V<\epsilon$ if we choose $\tilde{r}$ large enough.  Since $\delta>0$
could be chosen as small as we like provided $\tilde{r}$ was chosen large
enough, inequality \ref{one} follows, proving the theorem.  \qed

\section{Conjectures}
\label{conjectures}

We have seen that isoperimetric surfaces can be used to prove the Penrose
inequality for two classes of manifolds.  Naturally we want to generalize 
these results.  First we minimized area with a volume constraint and found
that this approach worked as long as the minimizing surfaces always had only
one component.  Then we realized that if we minimized the sum of the areas to
the three halves power with a volume constraint, then this approach worked 
even when the minimizer had multiple components.  However, this second approach
has a new problem, that two or more surfaces in the minimizing configuration 
can push up against each other.  

This suggests that we are still not optimizing the correct quantity.  
Minimizing the sum of the areas to the three halves power is a 
generalization of minimizing area (with a volume constraint) in the sense 
that these two optimization problems give the same answer when the minimizers
have only one component.  Hence, it is natural to consider how we can 
generalize the quantity ``sum of the areas to the three halves power'' in 
such a way that the new quantity equals the sum of the areas to the three
halves power in certain cases.  We recall the definition
of area nonincreasing maps given in definition 
\ref{area_nonincreasing} of section \ref{s6.5}
and propose the following functional.

\begin{definition}
Suppose $(M^3,g)$ is complete, has nonnegative scalar curvature, contains
one or more outermost minimal spheres $\{\Sigma_i\}$, and is asymptotically
flat at infinity.  Let
\[ f(D^3) = \sup_{\phi} 
\{\mbox{Vol}_{\real^3}(\phi(D^3)) \,\,|\,\, \phi:D^3 \rightarrow \real^3 
\mbox{is area nonincreasing}\} \]
and
\[ f(V) = \inf_{D^3} \{f(D^3) \,\,|\,\, D^3 \mbox{contains a volume V 
outside the horizons} \{\Sigma_i\} \}  \]
where $D^3$ is any open region in $M^3$ containing everything inside the 
horizons. 
\end{definition}

In general, there will not be a unique map $\phi$ which maximizes the 
$\real^3$ volume of $\phi(D^3)$ among area nonincreasing maps.  Instead, we
generally expect there to be a two dimensional ``critical set'' of $D^3$
which determines how large the $\real^3$ volume of $\phi(D^3)$ can be.  
For example, if the boundary of $D^3$ has area $A$, then by the isoperimetric
inequality in $\real^3$, the volume of $\phi(D^3)$ can be at most the volume
of a sphere in $\real^3$ with surface area $A$, which is $A^\frac32 / 
\sqrt{36\pi}$, and sometimes this upper bound is realized.  In fact, if
$D^3$ has several components each with boundary area $A_i$, then the volume 
of $\phi(D^3)$ can be at most the volume of disjoint balls with surface areas
$A_i$, which is $\sum A_i^\frac32 / \sqrt{36\pi}$, and sometimes this upper
bound is realized too.  Hence, we see that optimizing the functional $f$ 
is sometimes equivalent to optimizing the sum of the areas to the three halves
power.

\begin{conjecture}\label{anm}
Under the assumptions stated in the definition of $f(V)$, 
\[f''(V) \le \frac{1-f'(V)^2}{6f(V)} \]
\end{conjecture}

\begin{definition}
For $V \ge 0$, let
\[m(V) = f(V)^\frac13 (1 - f'(V)^2)/k \]
be the new mass function, where $k = (32 \pi /3)^\frac13$. 
\end{definition}   
 
From conjecture \ref{anm}, it follows that $m(V)$ is nondecreasing for 
$V \ge 0$.  Also, the original manifold can always be modified so that 
$m(0) =  \left(\sum_{i=1}^n \left(\frac{A_i}{16\pi}\right)^\frac32 \right)
^\frac13$ where $\{A_i\}$ are the areas of the horizons.  
Hence, it is likely that conjecture \ref{anm} would
imply the following generalized
Penrose inequality.

\begin{conjecture}\label{genpenrose2}
Suppose $(M^3,g)$ is complete, 
has nonnegative scalar curvature, contains one or more outermost
minimal spheres $\{\Sigma_i\}$ with surface areas $\{A_i\}$, and 
is Schwarzschild with mass $m$
at infinity.
Then $m \ge \left(\sum_{i=1}^n \left(\frac{A_i}{16\pi}\right)^\frac32 \right)
^\frac13$.
\end{conjecture}


\chapter{Volume Comparison Theorems}\label{volume}

The isoperimetric surface techniques which we developed to 
study the Penrose inequality in general relativity also can be
used to prove several volume comparison theorems, including a 
new proof of Bishop's volume comparison theorem for positive 
Ricci curvature. 
Let $(S^n,g_0)$ be the standard metric (with any scaling) on $S^n$
with constant Ricci curvature $Ric_0 \cdot g_0$.  Bishop's theorem says that if
$(M^n,g)$ is a complete Riemannian manifold ($n\ge 2$) with $Ric(g) \ge Ric_0 \cdot g$,
then $\mbox{Vol}(M^n) \le \mbox{Vol}(S^n)$.  It is then natural to ask
whether a similar type of volume comparison theorem could be true for
scalar curvature.  
We prove the following theorems for $3$-manifolds.

\begin{theorem} \label{T1}
Let $(S^3,g_0)$ be the constant curvature metric on $S^3$ with scalar 
curvature $R_0$, Ricci curvature $Ric_0 \cdot g_0$, and volume $V_0$.
Then there exists a positive $\epsilon_0 < 1$ such that 
if $(M^3,g)$ is any complete smooth Riemannian manifold of volume
$V$ satisfying
\[R(g) \ge R_0 \] 
\[Ric(g) \ge \epsilon_0 \cdot Ric_0 \cdot g\] 
then 
\[ V \le V_0. \]
\end{theorem}

As it happens, a lower bound on scalar curvature by itself is not 
sufficient to give an upper bound on the total volume.  We can scale
a cylinder, $S^2 \times \real$, to have any positive scalar curvature and still have
infinite volume.  
However, if we consider a neighborhood of metrics
around the standard metric $g_0$ on $S^3$ which satisfy
$Ric(g) \ge \epsilon_0 \cdot Ric_0 \cdot g$, then from the above theorem  we see
that $R(g) \ge R_0$ implies
that $V \le V_0$ for these metrics.  Hence, we see that a volume
comparison theorem 
for scalar curvature is true for metrics close to the standard
metric on $S^3$.  Moreover,
 
\begin{theorem}\label{T2} 
Let $(S^3,g_0)$ be the constant curvature metric on $S^3$ with scalar 
curvature $R_0$, Ricci curvature $Ric_0 \cdot g_0$, and volume $V_0$.
If $\epsilon \in (0,1]$   
and $(M^3,g)$ is any complete smooth Riemannian 3-manifold of volume
$V$ satisfying
\[R(g) \ge R_0 \] 
\[Ric(g) \ge \epsilon \cdot Ric_0 \cdot g\] 
then 
\[ V \le \alpha(\epsilon) V_0 \]
where 
\[\alpha(\epsilon)=\sup_{\frac{4\pi}{3-2\epsilon} \le z \le 4\pi}
\frac{1}{\pi^2}\left(\begin{array}{c}  
\int_0^{y(z)} 
\left(36\pi - 27 (1-\epsilon) y(z)^\frac23 - 9\epsilon\cdot
x^\frac23  \right)^{-\frac12} dx \\
+ \int_{y(z)}^{z^\frac32} 
(36\pi - 18(1-\epsilon)y(z) x^{-\frac13} - 9 x^{\frac23})^{-\frac12} dx
\end{array}  \right)\]
where
\[y(z) = \frac{z^\frac12 (4\pi - z)}{2(1-\epsilon)}. \]
Furthermore, this expression for $\alpha(\epsilon)$ is sharp.
\end{theorem}

Interestingly enough, $\alpha(\epsilon) = 1$ for many values of $\epsilon$.  
And since the above expression for $\alpha(\epsilon)$ is sharp, this allows
us to define the best value for $\epsilon_0$ which works in theorem \ref{T1},
namely 
\[\epsilon_0 = \inf \{ \epsilon \in (0,1] \,\,|\,\, \alpha(\epsilon) = 1 \}\]
Naturally it is desirable to estimate the actual value of 
$\epsilon_0$.  It is not too hard to show that 
$\epsilon_0 < 1$.  However, getting an accurate estimate for $\epsilon_0$
definitely seems to be a job for a computer, and it seems reasonable to conjecture
that $\epsilon_0$ is transcendental.  From preliminary computer calculations,
it looks like $0.134 < \epsilon_0 < 0.135$, although these bounds are not 
currently rigorous.


\section{Isoperimetric Surface Techniques}                              \label{S1}

As before in chapter \ref{Penrose},
isoperimetric surfaces will be used to prove these theorems.
The main difference is that we will be minimizing area with a 
volume constraint on compact manifolds in this chapter, so
existence of area minimizers is already known.  Also, the 
manifolds we will be considering all have positive Ricci 
curvature, from which it will follow from a stability argument 
that the area minimizers
always have exactly one component.  Hence, condition 
\ref{R3} from chapter \ref{Penrose} will always apply, so it 
will not be necessary to consider minimizing $F$ with a volume 
constraint.

\begin{definition}
Let $(M^n,g)$ be a complete Riemannian n-manifold.  Define
\[A(V) = \inf_R \{\mbox{Area}(\partial R)\,\,| \,\, 
                  \mbox{Vol}(R) = V\}\]
where $R$ is any region in $M^n$, $\mbox{Vol}(R)$ is the $n$
dimensional volume of $R$, and $\mbox{Area}(\partial R)$ is 
the $n-1$ dimensional volume of $\partial R$. 
If there exists a region $R$ with $\mbox{Vol}(R) = V$ such that
$\mbox{Area}(\partial R) = A(V)$, 
then we say that $\Sigma = \partial R$
minimizes area with the given volume constraint.  
\end{definition}

The manifolds we will be dealing with in this chapter all have
$Ric(g) \ge \delta > 0$.  Hence, these manifolds are compact, so there
will always exist a minimizer $\Sigma(V)$ (not necessarily
unique) for all $V$.  These minimal surfaces have constant mean
curvature and are smooth.    

We will use the function $A(V)$ to achieve the volume bounds on
$M^n$.  We will use the curvature bounds on $M^n$ to get an upper bound on
$A''(V)$.  Intuitively, this will force the two roots of $A(V)$ to
be close together.  Since the two roots of $A(V)$ are $0$ and
$\mbox{Vol}(M^n)$, we will get an upper bound for $\mbox{Vol}(M^n)$.

To get an upper bound for $A''(V)$ at $V=V_0$, we will do
a unit normal variation on $\Sigma(V_0)$.  That is, let
$\Sigma_{V_0}(t)$ be the surface created by flowing $\Sigma(V_0)$
out at every point in the normal direction at unit speed for time
$t$.  Since $\Sigma(V_0)$ is smooth, we can do this variation for
$t \in (-\delta,\delta)$ for some $\delta > 0$.  Abusing notation
slightly, we can also parameterize these surfaces by their volumes as
$\Sigma_{V_0}(V)$ so that $V = V_0$ will correspond to
$t = 0$.  Let $A_{V_0}(V) = \mbox{Area}(\Sigma_{V_0}(V))$.  Then 
$A(V_0) = A_{V_0}(V_0)$ and $A(V) \le A_{V_0}(V)$ since
$\Sigma_{V_0}(V)$ is not necessarily minimizing for its volume. 
Hence,

\[A''(V_0) \le A_{V_0}''(V_0).\] 
 
Let us 
suppose that $(M^3,g)$ satisfies $R(g) \ge R_0$ and $Ric(g) \ge 
\epsilon \cdot Ric_0 \cdot g$ as in theorem \ref{T2} and 
compute $A''_{V_0}(V_0)$.  To do this, we will need to compute the
first and second derivatives of the area of $\Sigma_{V_0}(t)$ and the volume
that it encloses.  We will use the formulas
\begin{equation}\label{formulas}
\dot{d\mu}=H\,d\mu\,\,\,\,\mbox{ and }\,\,\,\,\dot{H}=-||\Pi||^2-Ric(\nu,\nu)
\end{equation}
where the dot represents differentiation with respect to $t$, $d\mu$ is
the surface area 2-form for $\Sigma_{V_0}(t)$, $\Pi$ is the second fundamental
form for $\Sigma_{V_0}(t)$, $H=\mbox{trace}(\Pi)$ is the mean curvature, and
$\nu$ is the outward pointing unit normal vector.  Since $A_{V_0}(t)=
\int_{\Sigma_{V_0}(t)}d\mu$,
\[A_{V_0}'(t)=\int_{\Sigma_{V_0}(t)}H\,d\mu\]
And since $V'(t)=\int_{\Sigma_{V_0}(t)}d\mu=A_{V_0}(t)$,
\[A_{V_0}'(V_0)=A_{V_0}'(0)/V'(0)=H\]
By single variable calculus,
\[A_{V_0}''(V) = \frac{A_{V_0}''(t) - A_{V_0}'(V) V''(t)}{V'(t)^2}\]
so that at $t = 0$,
\begin{eqnarray*} 
A_{V_0}(V_0)^2 A_{V_0}''(V_0) &=& A_{V_0}''(t) - H V''(t) \\
                  &=& \frac{d}{dt} \int_{\Sigma_{V_0}(t)} H\, d\mu \,\, - \,\,
                      H \frac{d}{dt}\int_{\Sigma_{V_0}(t)} d\mu  \\
                  &=& \int_{\Sigma(V_0)} \dot{H}\, d\mu \\
                  &=& \int_{\Sigma(V_0)} -||\Pi||^2 - Ric(\nu,\nu) 
\end{eqnarray*}
Finally, since $||\Pi||^2 \ge \frac12 \mbox{trace}(\Pi)^2 = \frac12 H^2$ and 
$Ric(\nu,\nu) \ge \epsilon \cdot Ric_0$,

\begin{eqnarray*} 
A_{V_0}(V_0)^2 A_{V_0}''(V_0)&\le&\int_{\Sigma(V_0)} -\frac12 H^2 
- \epsilon \cdot Ric_0\\
  &=& - A_{V_0}(V_0) \left( \frac12 A_{V_0}'(V)^2+\epsilon\cdot Ric_0\right)  
\end{eqnarray*}
Hence,
\begin{equation}\label{AA}
A_{V_0}''(V_0) \le - \frac{1}{A_{V_0}(V_0)}
\left( \frac12 A_{V_0}'(V_0)^2 + \epsilon \cdot Ric_0 \right)
\end{equation}

\begin{lemma}\label{AB}
Suppose $\Sigma = \partial R$ minimizes area for its volume, 
$R \subset (M^3,g)$, and $Ric(g) \ge \delta>0$.  Then $\Sigma$ has exactly one 
component.
\end{lemma}
{\em Proof.} 
Suppose $\Sigma$ has more than one component.  Consider a flow on $\Sigma$
which is a unit normal flow (flowing out) on the first component 
(parameterized by volume) and a unit normal flow (flowing in) on the second
component (also parameterized by volume).  Then all of the surfaces in this 
family contain the same volume.  However, by equation \ref{AA} the second 
derivative of area is negative with respect to this volume preserving flow
(let $\delta = \epsilon \cdot Ric_0$).
Thus, $\Sigma$ does not minimize area for its volume.  Contradiction. \qed
 
Lemma \ref{AB} will be crucial for getting upper bounds on $A''(V)$ from the
lower bound on scalar curvature, and is one of the reasons we need some kind
of lower bound on Ricci curvature.  

Going back to equation
\ref{AA}, since $A(V_0) = A_{V_0}(V_0)$ and $A(V) \le A_{V_0}(V)$,  
\begin{equation}\label{AC}
A''(V) \le - \frac{1}{A(V)}
\left( \frac12 A'(V)^2 + \epsilon \cdot Ric_0 \right)
\end{equation}
in the sense of comparison functions defined in chapter \ref{Penrose}.

\begin{lemma} \label{AI}
Suppose $(M^3,g)$ satisfies $Ric(g) \ge \delta > 0$.  Then $A(V)$ is strictly
increasing on the interval $[0,\frac12 \mbox{Vol}(M^3)]$.
\end{lemma}
{\em Proof.}
It is always true that $A(V) = A(\mbox{Vol}(M^3) - V)$, since the boundaries of a 
region and its complement are the same.  By equation \ref{AC}, $A''(V)$ is 
strictly negative (again, let $\delta = \epsilon \cdot Ric_0)$.  The lemma
follows.  \qed

Now we want an equation like equation \ref{AC} which follows from the lower
bound on scalar curvature.  From before,
\[A_{V_0}(V_0)^2 A_{V_0}''(V_0) = \int_{\Sigma(V_0)} -||\Pi||^2-Ric(\nu,\nu)\] 
By the Gauss equation,
\[Ric(\nu,\nu) = \frac12 R - K + \frac12 H^2 - \frac12 ||\Pi||^2 \]
where $R$ is the scalar curvature of $M^3$ and $K$ is the Gauss curvature of 
$\Sigma(V_0)$.  Substituting we get,
\[A_{V_0}(V_0)^2 A_{V_0}''(V_0) = \int_{\Sigma(V_0)} 
-\frac12 R + K - \frac12 H^2 - \frac12 ||\Pi||^2 \]
By Lemma \ref{AB}, $\Sigma(V_0)$ has only one component, so by the 
Gauss-Bonnet theorem, $\int_{\Sigma(V_0)} K = 2\pi X(\Sigma(V_0)) \le 4\pi$.  
Since $R \ge R_0$ and $||\Pi||^2 \ge \frac12 H^2$, we have
\begin{eqnarray*}
A_{V_0}(V_0)^2 A_{V_0}''(V_0) &\le& 4\pi - \int_{\Sigma(V_0)} 
\frac12 R_0 + \frac34 H^2 \\
     &=& 4\pi - A_{V_0}(V_0)\left(\frac12 R_0 + \frac34 H^2 \right)
\end{eqnarray*}
Hence,
\[A_{V_0}''(V_0) \le \frac{4\pi}{A_{V_0}(V_0)^2} - \frac{1}{A_{V_0}(V_0)}
                     \left(\frac34 A_{V_0}'(V_0)^2 + \frac12 R_0\right)\]
As before, since $A(V_0) = A_{V_0}(V_0)$ and $A(V) \le A_{V_0}(V)$, 

\begin{equation}\label{AD2}
A''(V) \le \frac{4\pi}{A(V)^2} - \frac{1}{A(V)}
                     \left(\frac34 A'(V)^2 + \frac12 R_0\right)
\end{equation}
in the sense of comparison functions defined in chapter \ref{Penrose}.

Notice that we had to use the Gauss-Bonnet theorem to get equation \ref{AD2}.  
If we tried to generalize equation \ref{AD2} for higher dimensions, we would 
need to get an upper bound for $\int_{\Sigma(V_0)} R^\Sigma$, where 
$R^\Sigma$ is the scalar curvature of $\Sigma(V_0)$.  Since we don't have 
such a bound in general, the argument, as presented here, only works when $M$ 
is a 3-manifold.

However, equation \ref{AC} does generalize for all dimensions.  This allows 
us to give a new proof of Bishop's theorem, which we present in section
\ref{S3}.


\section{Ricci and Scalar Curvature Mass}                                                   \label{S2}

We define

\begin{equation}\label{AF}   F(V) = A(V)^\frac32   \end{equation} 
and choose to deal with $F(V)$ instead of $A(V)$.  Since $F(V)$ and $V$ have
the same units and $F(V)$ is roughly a linear function of $V$ for 
small $V$, the equations for $F(V)$ turn out to be simpler than the equations
for $A(V)$.  Of course, $F(V)$ and $A(V)$ will have the same roots, 
$0$ and $\mbox{Vol}(M^3)$,
and we will want to use upper bounds on $F''(V)$ to prove that the roots of 
$F(V)$ are close together, thereby getting an upper bound on $\mbox{Vol}(M^3)$.

Plugging equation \ref{AF} into equations \ref{AC} and \ref{AD2} and 
simplifying, we get

\begin{equation}\label{AG}
F''(V) \le - \frac{3 \epsilon \cdot Ric_0}{2} F(V)^{-\frac13}
\end{equation}
         and 
\begin{equation}\label{AH}
F''(V) \le \frac{36\pi - F'(V)^2}{6F(V)} - \frac{3 R_0}{4} F(V)^{-\frac13}
\end{equation}
in the sense of comparison functions defined in chapter \ref{Penrose}.  
We comment that it follows that these inequalities are also true distributionally.  
Given inequalities like equations \ref{AG} and \ref{AH}, it is natural to 
want to integrate them. 

\begin{definition} Let
\[m_{Ric}(V) = \left(36\pi - F'(V)^2 \right)-
               \frac{9 \epsilon \cdot Ric_0}{2} F(V)^{\frac23}\]
\[m_R(V) = F(V)^\frac13 \left(36\pi - F'(V)^2 \right) - \frac{3 R_0}{2} F(V)\]
and we call these two quantities ``Ricci curvature mass'' and ``scalar curvature mass''
respectively.
\end{definition}

$F(V)$ is continuous, but $F'(V)$ does not necessarily exist
for all $V$, although it does exist almost everywhere
since $F(V)$ is monotone increasing on 
$[0,\frac12\mbox{Vol}(M^3)]$ and monotone decreasing on
$[\frac12\mbox{Vol}(M^3), \mbox{Vol}(M^3)]$.  
The left and right hand derivatives,
$F'_+(V)$ and $F'_-(V)$, do always exist though.  This follows
from the fact that $F(V)$ has comparison functions
$F_{V_0}(V) = A_{V_0}(V)^\frac32$,
for all $V_0 \in (0,\mbox{Vol}(M^3))$ with uniformly bounded
second derivatives.  Hence, we can add a quadratic to $F(V)$
to get a concave function, from which it follows that the 
left and right hand derivatives exist and are equal except at
a countable number of points.
We define $F'_-(0) = \sqrt{36\pi}$ and 
$F'_+(\mbox{Vol}(M^3)) = -\sqrt{36\pi}$, which is natural
for smooth manifolds.

Furthermore, $F'_+(V) \le F'_-(V)$ using the comparison
function argument again since 
$F'_+(V_0) \le F'_{V_0}(V_0) \le F'_-(V_0)$.  If $F'(V)$ 
does not exist, then it is natural to define $F'(V)$ to be a 
multivalued function taking on every value in the interval
$(F'_+(V), F'_-(V))$.  This is consistent, since if $F'(V)$
does exist, then $F'_+(V) = F'_-(V)$.  Hence, $m_{Ric}(V)$
and $m_R(V)$ are 
multivalued for some $V$, which can be interpreted as the mass
``jumping up'' at these $V$, and the set of $V$ for which 
$m(V)$ and $F(V)$ are multivalued is a countable set.  
Alternatively, one
could replace $F'(V)$ with $F'_+(V)$ (or $F'_-(V)$) in the 
formulas for $m_{Ric}(V)$ and $m_R(V)$ 
so that they would always be single 
valued.

\begin{lemma} \label{AJ2}
The quantities $m_{Ric}(V)$ and $m_R(V)$ are nonnegative, nondecreasing 
functions of $V$ on the interval $[0,\frac12 \mbox{Vol}(M^3)]$ and 
$m_{Ric}(0) = m_R(0) = 0$.  
\end{lemma}
{\em Proof.}  
Since $M^3$ is a smooth manifold, $F(V) \approx \sqrt{36\pi} V$ for small $V$ and  
$F'(0) = \sqrt{36\pi}$.  Since $F(0) = 0$, it follows that
$m_{Ric}(0) = m_R(0) = 0$.  In addition, we observe that
if $F(V)$ were smooth, then 

\[m_{Ric}'(V) = 2 F'(V) \left(-F''(V) 
              - \frac{3 \epsilon \cdot Ric_0}{2} F^{-\frac13} \right)\]
and

\[m_R'(V) = 2 F^\frac13 F'(V) \left(-F''(V) + \frac{36\pi - F'(V)^2}{6F(V)}
          - \frac{3 R_0}{4} F^{-\frac13} \right)\]
Then by lemma \ref{AI}, we would have $F'(V) \ge 0$, 
so that by equations \ref{AG} and \ref{AH},
$m_{Ric}'(V) \ge 0$ and $m_R'(V) \ge 0$
proving that $m_{Ric}(V)$ and $m_R(V)$ are nondecreasing,
and hence nonnegative,
on the interval $[0,\frac12 \mbox{Vol}(M^3)]$.  

More generally, we need to prove that $m_{Ric}'(V) \ge 0$ 
and $m_R'(V) \ge 0$ as distributions, which follows as
before in the proof of lemma \ref{AJ}.\qed

The reason the we call the two quantities $m_{Ric}$ and $m_R$ ``mass'' is motivated by the fact that if we set 
$R_0 = 0$, $m_R(V) = m(V)$, where $m(V)$ is the mass function 
from chapter \ref{Penrose}.  However, beyond being nonnegative,
nondecreasing functions which are very similar
to $m(V)$ in form,
the author is not currently aware of any
physical interpretations of $m_{Ric}(V)$ and $m_R(V)$ in the
context of general relativity, although that is
an interesting possibility.


\section{A New Proof of Bishop's Theorem}                                 \label{S3}

In this section, we will give a new proof of Bishop's theorem using an argument which is 
very similar to the one we will use to prove theorem \ref{T2}.  
Whereas the rest of this chapter deals specifically with $3$-manifolds, 
in the next section we will study $n$-manifolds.  Because of this, we will need to
generalize a few definitions and equations just for this
section.

\begin{theorem} \label{AK}
{\bf (Bishop)}  Let $(S^n,g_0)$ be the standard metric (with any scaling) on $S^n$
with constant Ricci curvature $Ric_0 \cdot g_0$.  If
$(M^n,g)$ is a complete Riemannian manifold ($n\ge 2$) with $Ric(g) \ge Ric_0 \cdot g$,
then $\mbox{Vol}(M^n) \le \mbox{Vol}(S^n)$. 
\end{theorem}
{\em Proof.}
Modifying equation \ref{AC}, since $||\Pi||^2 \ge \frac1{n-1} \mbox{trace}(\Pi)^2 = 
\frac1{n-1} H^2$, ($\epsilon = 1$)
\begin{equation}\label{AL}
A''(V) \le - \frac{1}{A(V)}
\left( \frac1{n-1} A'(V)^2 + Ric_0 \right)
\end{equation}
Now we let $F(V) = A(V)^\frac{n}{n-1}$, from which it follows that
\begin{equation}\label{AM}
F''(V) \le - \frac{n \cdot Ric_0}{n-1} F(V)^{-\frac{n-2}{n}}
\end{equation} 
The correct definition for $m_{Ric}(V)$ then becomes
\begin{equation}\label{AN}
m_{Ric}(V) = \left(n^2 (\omega_{n-1})^\frac2{n-1} - F'(V)^2 \right)-
               \frac{n^2 \cdot Ric_0}{n-1} F(V)^{\frac2n}
\end{equation}
where $\omega_{n-1}$ is the surface area of the sphere $S^{n-1}$ of radius $1$ in $\real^n$.
As before, on the interval $[0,\frac12 \mbox{Vol}(M^n)]$, $F'(V) \ge 0$ and 
$m_{Ric}(V)$ is a nonnegative, nondecreasing function of $V$.
The proof is the same as before.

Now consider phase space which we will view as the $x$-$y$ plane where $x = F(V)$ and
$y = F'(V)$.  Let $\gamma$ be the path in phase space of F(V) for $V$ between $0$ and 
$\frac12 \mbox{Vol}(M^n)$.  Then we note that
\begin{equation} \label{AO}
\frac12 \mbox{Vol}(M^n) = \int_\gamma dV = \int_\gamma \frac{dx}{y}
\end{equation}   
We also observe that since $F(0) = 0$ and $F'(\frac12 \mbox{Vol}(M^n)) = 0$ (by the 
symmetry of $F(V)$), $\gamma$ is a path from the $y$ axis to the $x$ axis. Since
$F(V)$ is strictly increasing and $F'(V)$ is strictly decreasing (by inequality \ref{AM}), 
the $x$ position of $\gamma$ is nondecreasing and the $y$ position of $\gamma$ is strictly
decreasing. 
Since $F'(V)$ is sometimes multivalued, taking on the values
of an interval, $\gamma$ is sometimes vertical.

Now we want to find the $\gamma$ which maximizes equation \ref{AO}, with the 
constraint that $m_{Ric}(V)$ stays nonnegative and nondecreasing, which is
equivalent to satisfying equation \ref{AM}.  Consider all the 
possible paths which terminate at a given point on the $x$ axis, $(x_0,0)$, and 
think of these paths as beginning at this point and then follow the paths backwards.  
The path which maximizes equation \ref{AO} will be 
the one which has the smallest $y$ values.  Since $F''(V) = y \frac{dy}{dx}$, 
we can rewrite inequality \ref{AM} as

\begin{equation}\label{AP}
\frac{dy}{dx} \le - \frac{n \cdot Ric_0}{n-1} x^{-\frac13} y^{-1}
\end{equation} 
The $\gamma$ terminating at $(x_0, 0)$ with the smallest $y$ values will be the path 
which has equality in inequality \ref{AP}.  

Hence, this path is given by the $F(V)$ which has equality in inequality \ref{AM}.  
But equality for inequality \ref{AM} is 
equivalent to $m_{Ric}'(V) = 0$, which implies that $m_{Ric}(V) = m_0$, where 
$m_0$ is some positive constant.  By equation \ref{AN}, this path can be computed
explicitly and is given by
\[m_0 = \left(n^2 (c_{n-1})^\frac2{n-1} - y^2 \right)-
               \frac{n^2 \cdot Ric_0}{n-1} x^{\frac2n} \]
which can be rewritten as
\begin{equation} \label{AQ}
y = \left[\left(n^2 (c_{n-1})^\frac2{n-1} - m_0 \right)-
               \frac{n^2 \cdot Ric_0}{n-1} x^{\frac2n}\right]^\frac12
\end{equation}
Different values of $m_0$ correspond to curves terminating at different points on 
the $x$ axis.  Hence, the $\gamma$ which maximizes equation \ref{AO} is a curve which
is the graph of equation \ref{AQ} for some $m_0$.  By a simple change of variables, it
is easy to compute that
\begin{eqnarray*} 
\frac12 \mbox{Vol}(M^n) &=& \int_\gamma \frac{dx}{y} \le \sup_\gamma\int_\gamma \frac{dx}{y}  \\
                 &=& \sup_{m_0}\left(n^2 (c_{n-1})^\frac2{n-1} - m_0 \right)^\frac{n-1}2
                     \left(\frac{n-1}{n^2 \cdot Ric_0}\right)^\frac{n}{2}
                     \int_0^1 \left[1-z^\frac2n\right]^{-\frac12} dz \\
\end{eqnarray*}  
Now we recall that $m(V)$ is nonnegative, so $m_0$ must also be nonnegative.  Hence, the above
expression is maximized when $m_0 = 0$. 
But the standard sphere $(S^n,g_0)$
with constant Ricci curvature $Ric_0 \cdot g_0$ has 
$m_{Ric}(V) = 0$.  This can be verified by direct computation
using the fact that the isoperimetric spheres of $(S^n,g_0)$
are the spherically symmetric $(n-1)$-spheres, or from noticing
that since we get equality in inequality \ref{AL}
when $(M^n,g) = (S^n,g_0)$, $m'(0) \equiv 0$, so $m(0) \equiv 0$.
Let $\gamma_0$ be the path in phase space
corresponding to this standard sphere with zero mass.  Then
\begin{equation}\label{AR} 
\frac12 \mbox{Vol}(M^n) = \int_\gamma \frac{dx}{y} \le \sup_\gamma\int_\gamma \frac{dx}{y} 
                 = \int_{\gamma_0} \frac{dx}{y} = \frac12 \mbox{Vol}(S^n) 
\end{equation}  
proving the theorem.  \qed

\section{Proof of the Volume Comparison Theorems \\
involving Scalar Curvature}   \label{S6}

{\em Proof.}
The approach we take here is the same as we used to prove Bishop's theorem in section
\ref{S3}.  Going back to section \ref{S2} and combining equations \ref{AG}
and \ref{AH}, we get    
\begin{equation}\label{AW}
F''(V) \le \mbox{min} \left\{\frac{36\pi - F'(V)^2}{6F(V)} - \frac{3 R_0}{4} F(V)^{-\frac13},
                      - \frac{3 \epsilon \cdot Ric_0}{2} F(V)^{-\frac13} \right\}
\end{equation}
Since  
\begin{equation}\label{BA}
m_{Ric}(V) = \left(36\pi - F'(V)^2 \right)-
               \frac{9 \epsilon \cdot Ric_0}{2} F(V)^{\frac23}
\end{equation}
and
\begin{equation}\label{BB}
m_R(V) = F(V)^\frac13 \left(36\pi - F'(V)^2 \right) - \frac{3 R_0}{2} F(V)
\end{equation}
we can rewrite inequality \ref{AW} as
\begin{equation}\label{AX}
F''(V)  \le  -\frac12 F^{-\frac13} \cdot \mbox{max} \left\{ L(V), 3 \epsilon \cdot Ric_0 \right\} 
\end{equation}
where
\begin{equation}\label{BE}
L(V) =  R_0-\frac{m_R(V)}{3 F} = \frac32 (R_0-\epsilon\cdot Ric_0)-\frac{m_{Ric}(V)}{3 F^\frac23} 
\end{equation}

As before, we consider phase space which we will view as the $x$-$y$ plane where $x = F(V)$ and 
$y = F'(V)$.  Let $\gamma$ be the path in phase space of $F(V)$ for $V$ between $0$ and 
$\frac12 \mbox{Vol}(M^3)$.  Then we recall that
\begin{equation}\label{AZ} 
\frac12 \mbox{Vol}(M^3) = \int_\gamma dV = \int_\gamma \frac{dx}{y}
\end{equation}   
Since $F(0) = 0$ and $F'(\frac12 \mbox{Vol}(M^3)) = 0$ (by the 
symmetry of $F(V)$), $\gamma$ is a path from the $y$ axis to the $x$ axis.  
Since
$F(V)$ is strictly increasing and $F'(V)$ is strictly decreasing (by inequality \ref{AX}), 
the $x$ position of $\gamma$ is nondecreasing and the $y$ position of $\gamma$ is strictly
decreasing.  Again, since $F'(V)$ is sometimes multivalued,
taking on the values of an interval, $\gamma$ is sometimes
vertical.

We want to find the $\gamma$ which maximizes equation \ref{AZ}, while still satisfying
inequality \ref{AX}.  
Consider all the 
possible paths which terminate at a given point on the $x$ axis, $(x_0,0)$, and 
think of these paths as beginning at this point and then follow the paths backwards.  
The path which maximizes equation \ref{AZ} will be 
the one which has the smallest $y$ values.  Since $F''(V) = y \frac{dy}{dx}$, 
we can rewrite inequality \ref{AX} as
\begin{equation}\label{BC}
\frac{dy}{dx} \le -\frac12 x^{-\frac13} y^{-1} \cdot 
                   \mbox{max} \left\{L(V), 3 \epsilon \cdot Ric_0 \right\}
\end{equation}
where we think of $L(V)$, $m_R(V)$, and $m_{Ric}(V)$ as functions of $x$ and $y$ instead
of $F(V)$ and $F'(V)$.
The $\gamma$ terminating at $(x_0, 0)$ with the smallest $y$ values will be the path 
which has equality in inequality \ref{BC}, and thus has equality in 
inequality \ref{AW}.  Let's call this path $\gamma(x_0)$.  Then we see that 
$\gamma(x_0)$ maximizes equation \ref{AZ}
among all paths which terminate at $(x_0, 0)$.  

By the computations in the proof of lemma \ref{AJ2}, equality in inequality \ref{AW} 
is equivalent to either $m_R'(V) = 0$ or $m_{Ric}'(V) = 0$ for each $V$.  This, combined
with equation \ref{BE} and the fact that $F(V)$ is a strictly increasing function of $V$, 
gives us that 
$L(V)$ is strictly increasing as a function of $V$ for the path $\gamma(x_0)$.  
Furthermore, 
\begin{eqnarray*}
L(V) \ge 3 \epsilon \cdot Ric_0  & \Rightarrow &  m_R'(V) = 0 \\
L(V) \le 3 \epsilon \cdot Ric_0  & \Rightarrow &  m_{Ric}'(V) = 0
\end{eqnarray*}
Hence, we see that there are three cases.  

Case 1:  If $L(V)$ is always less than $3 \epsilon \cdot Ric_0$, then
$\gamma(x_0)$ is the curve given by $m_{Ric}(V) = c$ for some constant $c \ge 0$ which
depends on $x_0$.  

Case 2:  If $L(V)$ is initially smaller than $3 \epsilon \cdot Ric_0$ but becomes larger than
$3 \epsilon \cdot Ric_0$ for $V \ge \tilde{V}$,
then $\gamma(x_0)$ will be the union of two segments of curves,
one given by $m_{Ric}(V) = c_2$ for $V \le \tilde{V}$ and the other given by $m_R(V) = c_1$
for $V \ge \tilde{V}$, for three constants $c_1, c_2 \ge 0$ and $\tilde{V}$ which depend on $x_0$.

Case 3:  If $L(V)$ is always greater than $3 \epsilon \cdot Ric_0$, then $m_R(V) = c$ 
for some constant $c \ge 0$.  By equation \ref{BE}, though, we see that if this constant
were positive, then $L(V)$ would approach $- \infty$ for small $F$, thus violating our
assumption that $L(V) \ge 3 \epsilon \cdot Ric_0$. 
Hence, $c = 0$, so $m_R(V) = 0$ for all $V$. 

From these observations, we explicitly compute $\gamma(x_0)$. 
We spare the reader some
of the routine details and summarize the results.  For convenience, we define
two special values for $x_0$.  Let 
\[x_S = \left(\frac{24\pi}{R_0}\right)^\frac32  \,\,\,\,\mbox{ and }\,\,\,\,
x_{FB} = \left(\frac{8\pi}{R_0 - 2\epsilon\cdot Ric_0}\right)^\frac32 \]
The subscripts
stand for ``sphere'' and ``football,'' since the standard 3-sphere
produces the curve $\gamma(x_S)$ and the metric (which turns out
to have two singularities) which produces
$\gamma(x_{FB})$ looks like an axially symmetric football
(with two pointy ends) when embedded in $\real^4$.

For $0 \le x_0 \le x_{FB}$,
$\gamma(x_0)$ is the graph of the function
\begin{equation}\label{BH}
y   = \left[  \frac{9 \epsilon \cdot Ric_0}{2} (x_0^\frac23   - x^{\frac23}) \right] ^\frac12
\end{equation}
These are the curves from case 1. 

For $x_{FB} \le x_0 \le x_S$, $\gamma(x_0)$ is the graph of the function 
\begin{equation}\label{BM}
y = \left\{ \begin{array}{cl} \left( 36\pi - c_2 - \frac{9 \epsilon \cdot Ric_0}{2} 
    x^\frac23 \right)^\frac12  &  , 0 \le x \le x_1  \\    
\left( 36\pi - c_1 x^{-\frac13} - \frac{3 R_0}{2} x^\frac23 \right)^\frac12  &
    , x_1 \le x \le x_0  \end{array} \right.
\end{equation}
where 
\[ c_1 = x_0^\frac13 \left( 36\pi - \frac{3 R_0}{2} x_0^\frac23 \right)\]
\[ x_1 = \frac{c_1}{3(R_0 - 3 \epsilon \cdot Ric_0)}\]
\[ c_2 = \frac32 \left[ 3 (R_0 - 3 \epsilon \cdot Ric_0) c_1^2 \right]^\frac13 \]
These are the curves from case 2, and the constants $c_1$ and $c_2$ are the same constants
that are mentioned in case 2.  Case 3 is also included here, and occurs when
$x_0 = x_S$, which implies that $c_1 = x_1 = 0$.    

There are no paths which terminate at $(x_0,0)$ for 
$x_0 > x_S$.  This follows from the definition of
$m_R(V)$ in equation \ref{BB} and the fact that $m_R(V) \ge 0$.

Now let's define
\begin{equation}\label{BJ} 
W(x_0) = \int_{\gamma(x_0)} dV = \int_{\gamma(x_0)} \frac{dx}{y}
\end{equation}   
Then we have that

\begin{equation}\label{BK} 
\frac12 V = \frac12 \mbox{Vol}(M^3) = \int_\gamma \frac{dx}{y} \le \sup_\gamma \int_\gamma \frac{dx}{y} 
                 = \sup_{x_0} W(x_0) 
\end{equation}  
where $M^3$ is any arbitrary $3$-manifold satisfying the curvature conditions of 
theorem \ref{T2}.  Using equations \ref{BH} and \ref{BM}, we can compute $W(x_0)$ explicitly.
\begin{equation}\label{BN}
W(x_0) = \left\{ \begin{array}{cl} \left( \frac{9 \epsilon \cdot Ric_0}{2} \right)^{-\frac12}
                 x_0^\frac23 \int_0^1 \left(1-z^\frac23\right)^{-\frac12} dz,
                 &  0 \le x_0 \le x_{FB}  \\ & \\

                 \int_0^{x_1} \left( 36\pi - c_2 - \frac{9 \epsilon \cdot Ric_0}{2} 
                 x^\frac23 \right)^{-\frac12} dx  & \\

                 + \int_{x_1}^{x_0} \left( 36\pi - c_1 x^{-\frac13}
                 - \frac{3 R_0}{2} x^\frac23 \right)^{-\frac12} dx, & x_{FB} < x_0 \le x_S
\end{array} \right. \end{equation}
where we've simplified the top integral using a change of variables.  Unfortunately, it is
not so easy to simplify the bottom integral.  However, we can simplify the bottom integral when
$x_0$ equals $x_{FB}$ or $x_S$ (because the values of $c_1$ and $x_1$ work out nicely), 
and we find that 
\begin{equation}
W(x_S) =  \frac{36\pi}{\left( \frac{9 Ric_0}{2} \right)^\frac32}        
          \int_0^1 \left(1-z^\frac23\right)^{-\frac12} dz
\end{equation}
(using $R_0 = 3 Ric_0$) and
\begin{equation}
W(x_{FB}) =  \frac{x_{FB}^\frac23}{\left( \frac{9 \epsilon \cdot Ric_0}{2} \right)^\frac12}        
             \int_0^1 \left(1-z^\frac23\right)^{-\frac12} dz
\end{equation}
Now
\begin{equation}
\int_0^1 \left(1-z^\frac23\right)^{-\frac12} dz = 3\pi/4
\end{equation}
so it is easy to check that 
\[W(x_S) = \frac12 \mbox{Vol}(S^3,g_0) = \frac12 V_0.\]
Furthermore, by the definition of $x_{FB}$, 
\begin{equation}
W(x_{FB}) = W(x_S) \frac{1}{\epsilon^\frac12 (3-2\epsilon)}
\end{equation}

We can simplify things further if we recognize the fact that everything
scales as it should.  Using the values from the 3-sphere of radius one embedded
in $\real^4$, we use $R_0 = 6$, $Ric_0 = 2$, and $V_0 = 2\pi^2$ to get
\[x_S = (4\pi)^\frac32\]
\[x_{FB} = \left(\frac{4\pi}{3-2\epsilon}\right)^\frac32\] 
\[c_1 = x_0^\frac13 (36\pi - 9 x_0^\frac23) \]
\[x_1 = \frac{c_1}{18(1-\epsilon)} \]
\[c_2 = \frac32[18(1-\epsilon)c_1^2]^\frac13  \]
so that plugging in these values for $W$ and scaling appropriately we get
\[ V \le  \alpha(\epsilon) \cdot V_0 \]
where
\[ \alpha(\epsilon) = \sup_{0 \le x_0 \le (4\pi)^\frac32} w_\epsilon(x_0) \]
where
\[ w_\epsilon(x_0) = \frac{1}{\pi^2}\left\{ \begin{array}{cl}  
\frac{\pi}{4}\cdot\epsilon^{-\frac12}\cdot x_0^\frac23, & 
0 \le x_0 \le \left(\frac{4\pi}{3-2\epsilon}\right)^\frac32 \\ 
\int_0^{\frac{c_1}{18(1-\epsilon)}} 
\left(36\pi - \frac32[18(1-\epsilon)c_1^2]^\frac13 - 9\epsilon\cdot
x^\frac23  \right)^{-\frac12} dx & \\ 
+ \int_{\frac{c_1}{18(1-\epsilon)}}^{x_0} 
(36\pi - c_1 x^{-\frac13} - 9 x^{\frac23})^{-\frac12} dx, &
\left(\frac{4\pi}{3-2\epsilon}\right)^\frac32 < x_0 \le  (4\pi)^\frac32 
\end{array}\right.\]  
where again
\[c_1 = x_0^\frac13 (36\pi - 9 x_0^\frac23). \]

We can simplify the notation a bit by changing variables.  Let 
$x_0 = z^\frac32$ and $c_1 = 18(1-\epsilon)y$.  Then we have
\[ \alpha(\epsilon) = \sup_{0 \le z \le 4\pi} w_\epsilon(z) \]
where
\[ w_\epsilon(z) = \frac{1}{\pi^2}\left\{ \begin{array}{cl}  
\frac{\pi}{4}\cdot\epsilon^{-\frac12}\cdot z, & 
0 \le z \le \frac{4\pi}{3-2\epsilon} \\ 
\int_0^{y(z)} 
\left(36\pi - 27 (1-\epsilon) y(z)^\frac23 - 9\epsilon\cdot
x^\frac23  \right)^{-\frac12} dx & \\ 
+ \int_{y(z)}^{z^\frac32} 
(36\pi - 18(1-\epsilon)y(z) x^{-\frac13} - 9 x^{\frac23})^{-\frac12} dx, &
\frac{4\pi}{3-2\epsilon} < z \le  4\pi 
\end{array}\right.\]  
where 
\[y(z) = \frac{z^\frac12 (4\pi - z)}{2(1-\epsilon)}. \]

Since $w_\epsilon$ is continuous, the maximum value must occur for 
$z \in [\frac{4\pi}{3-2\epsilon}, 4\pi]$.  Hence, we have the following 
theorem.

\begin{theorem} \label{scalar}
Let $(S^3,g_0)$ be the constant curvature metric on $S^3$ with scalar 
curvature $R_0$, Ricci curvature $Ric_0 \cdot g_0$, and volume $V_0$.
If $\epsilon \in (0,1]$   
and $(M^3,g)$ is any complete smooth Riemannian manifold of volume
$V$ satisfying
\[R(g) \ge R_0 \] 
\[Ric(g) \ge \epsilon \cdot Ric_0 \cdot g\] 
then 
\[ V \le \alpha(\epsilon) V_0 \]
where 
\[\alpha(\epsilon)=\sup_{\frac{4\pi}{3-2\epsilon} \le z \le 4\pi}
\frac{1}{\pi^2}\left(\begin{array}{c}  
\int_0^{y(z)} 
\left(36\pi - 27 (1-\epsilon) y(z)^\frac23 - 9\epsilon\cdot
x^\frac23  \right)^{-\frac12} dx \\
+ \int_{y(z)}^{z^\frac32} 
(36\pi - 18(1-\epsilon)y(z) x^{-\frac13} - 9 x^{\frac23})^{-\frac12} dx
\end{array}  \right)\]
where
\[y(z) = \frac{z^\frac12 (4\pi - z)}{2(1-\epsilon)}. \]
Furthermore, this expression for 
$\alpha(\epsilon)$ is sharp.
\end{theorem}

We note that the reason that this expression for 
$\alpha(\epsilon)$ is sharp is that it is possible to 
construct spherically symmetric manifolds which satisfy 
the curvature conditions of theorem \ref{scalar} and have volumes
as close to $\alpha(\epsilon)V_0$ as desired, and equal to  
$\alpha(\epsilon)V_0$ if we allow the manifolds to have 
singularities.  These manifolds look like long and skinny
axially symmetric footballs when embedded in $\real^4$ with 
two pointy ends where the manifold is not smooth.  The smaller
$\epsilon$ is, the longer and skinnier 
these manifolds become, and as 
$\epsilon$ goes to zero, these ``case of equality'' manifolds
converge to the standard cylinder $S^2 \times \real$ 
which has constant
scalar curvature $R_0$ 
(and zero Ricci curvature in the directions along the length
of the cylinder).  These manifolds can be constructed by looking
at the function $W(x_0)$ from equation \ref{BN} 
(for each value of $\epsilon$) and defining
$\bar{x}(\epsilon)$ to be the value of $x_0$ which maximizes $W$.
Then the curve $\gamma(\bar{x}(\epsilon))$ in phase space as 
described before corresponds to an $F(V)$ function, which
yields an $A(V)$ function using $F(V) = A(V)^{3/2}$.  Given an
$A(V)$ function, we can then construct a spherically symmetric
manifold such that the spherically symmetric spheres which 
contain a volume $V$ have surface area $A(V)$, and it is easy to
verify that these are in fact ``case of equality'' manifolds. 

Direct computation shows that $w_\epsilon(z)$ is a $C^1$ function on 
$[0,4\pi]$ and that $w_\epsilon(4\pi) = 1$,
so $\alpha(\epsilon) \ge 1$ for all $\epsilon \in (0,1]$.  Furthermore,
since direct calculation also shows that $w_\epsilon(z)$ is a nonincreasing
function of $\epsilon$ when $z$ is held fixed, it follows that 
$\alpha(\epsilon)$ is nonincreasing.  Hence, if $\alpha$ equals one at one
value of $\epsilon$, then $\alpha$ equals one for all larger values of 
$\epsilon$ in the interval $(0,1]$.  Let
\[\epsilon_0 = \inf \{ \epsilon \in (0,1] \,\,|\,\, \alpha(\epsilon) = 1 \}\]
Then we have the following theorem.

\begin{theorem} \label{vcthm}
Let $(S^3,g_0)$ be the constant curvature metric on $S^3$ with scalar 
curvature $R_0$, Ricci curvature $Ric_0 \cdot g_0$, and volume $V_0$.
If $(M^3,g)$ is any complete smooth Riemannian manifold of volume
$V$ satisfying
\[R(g) \ge R_0 \] 
\[Ric(g) \ge \epsilon_0 \cdot Ric_0 \cdot g\] 
then 
\[ V \le V_0. \]
\end{theorem}

Naturally it would be desirable to estimate the actual value of 
$\epsilon_0$.  It is straightforward (although messy)
to show that 
$\epsilon_0 < 1$.  However, getting an accurate estimate for $\epsilon_0$
definitely seems to be a job for a computer, and it seems reasonable to conjecture
that $\epsilon_0$ is
transcendental.  From preliminary computer calculations,
it looks like $.134 < \epsilon_0 < .135$, although these bounds 
are not rigorous.

\section{Estimates for $\epsilon_0$}

\begin{figure}
\vspace{6in}
\includegraphics{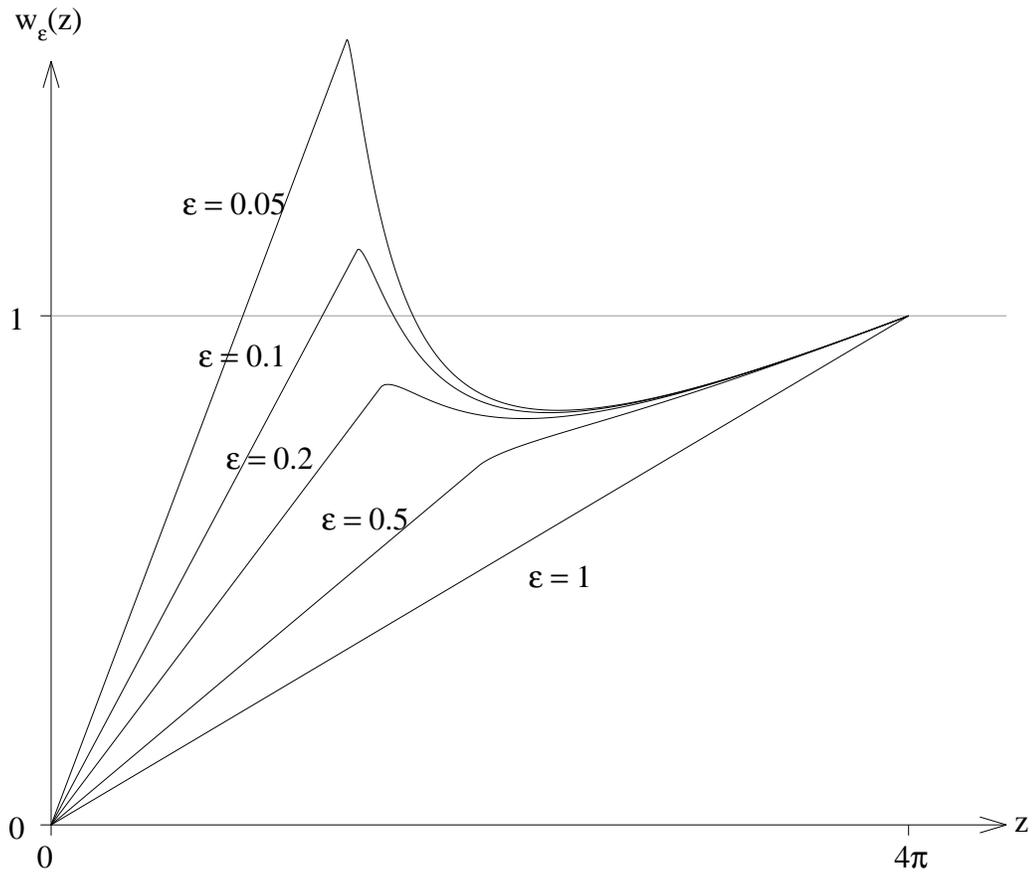}
\caption{Graphs of $w_\epsilon(z)$ for 
         $\epsilon = 0.05, 0.1, 0.2, 0.5$, and $1$.
         \label{epsilon1}}
\end{figure}

\begin{figure}
\vspace{6in}
\includegraphics{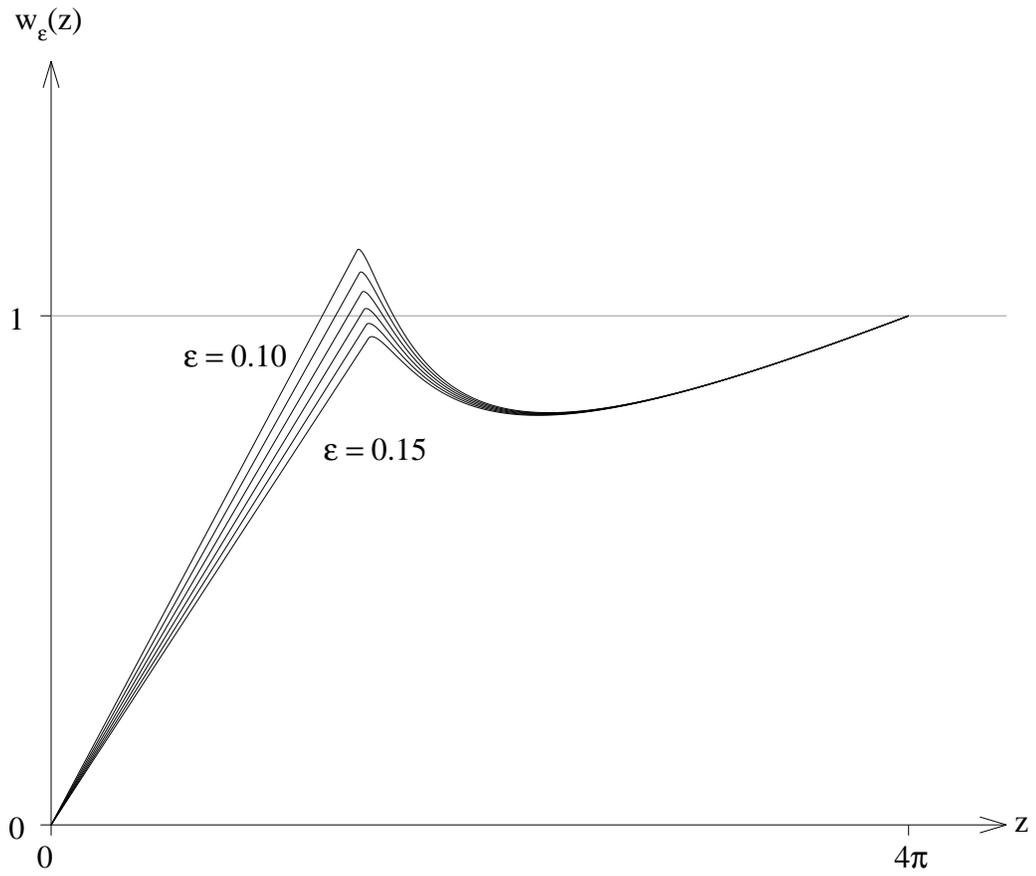}
\caption{\label{epsilon2}
Graphs of $w_\epsilon(z)$ for 
$\epsilon = 0.10, 0.11, 0.12, 0.13, 0.14$, and $0.15$.}
\end{figure}

\begin{figure}
\vspace{6in}
\includegraphics{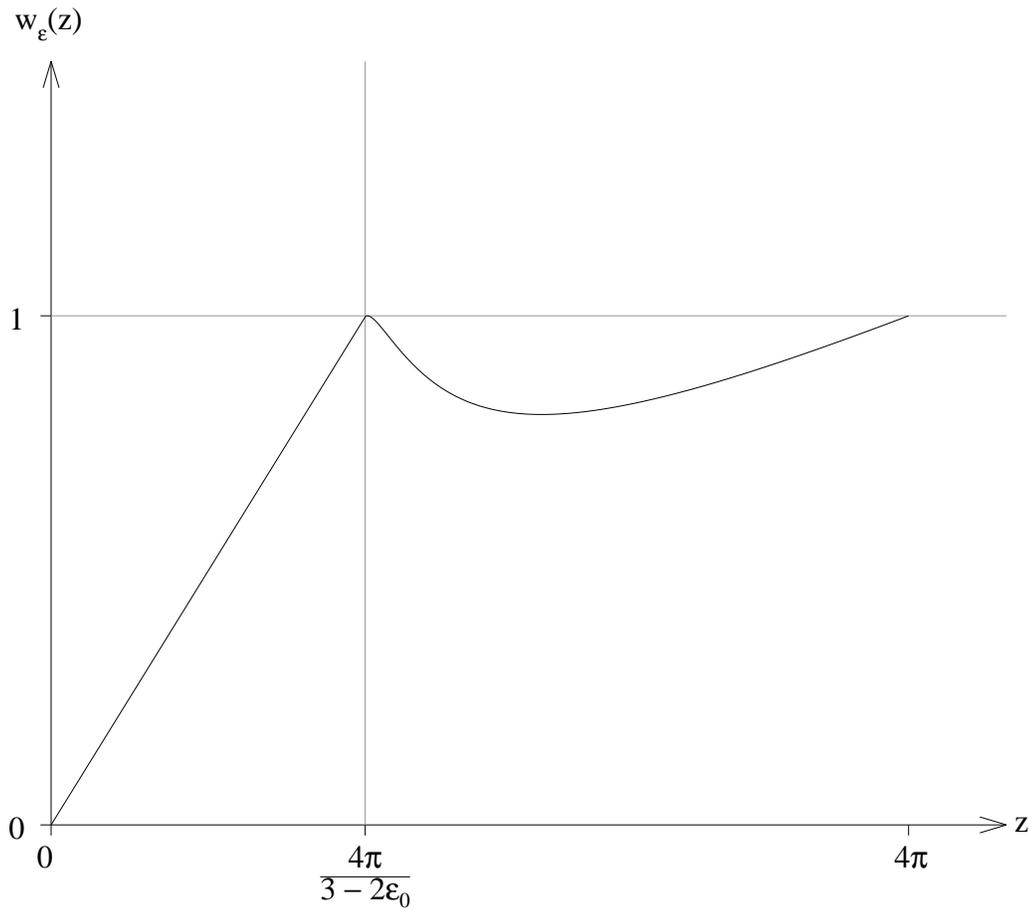}
\caption{\label{epsilon3}
Graph of $w_\epsilon(z)$ for $\epsilon = 0.134727$.}
\end{figure}

\begin{figure}
\vspace{6in}
\includegraphics{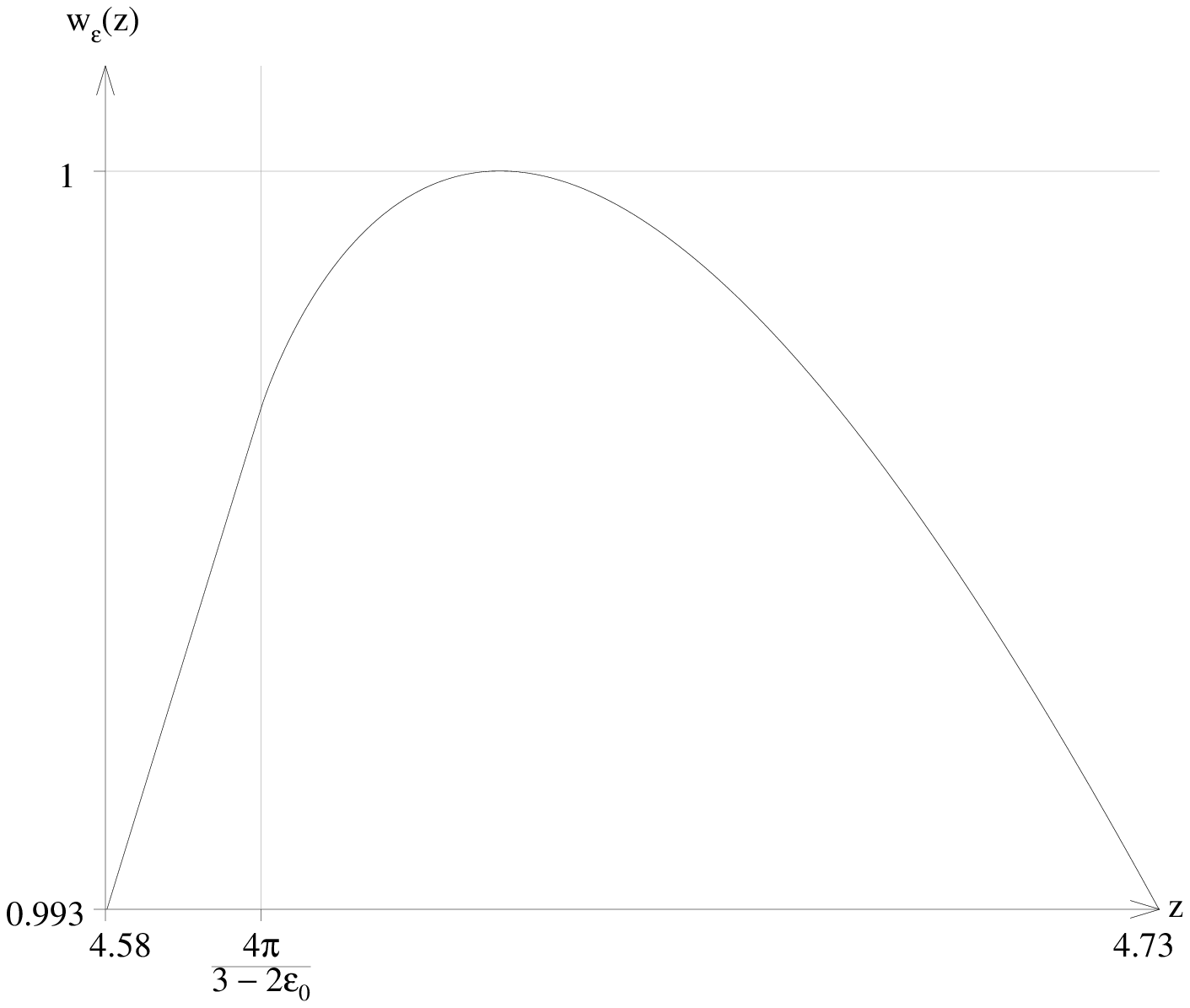}
\caption{\label{epsilon4}
Graph of $w_\epsilon(z)$ for $\epsilon = 0.134727$,
centered on the interior maximum point. }
\end{figure}

The results of this section are due primarily to Kevin Iga
of Stanford University, who wrote several computer programs
using the C programming language on a Sun SPARC station 20
computer to estimate the value
of $\epsilon_0$ from theorem \ref{vcthm}. 
We found that $\epsilon_0 \approx 0.134727$.  However, 
the only 
rigorous bounds that we have are  
\begin{equation}\label{rigorous}
 0.133974 < 1-\frac{\sqrt{3}}{2} < \epsilon_0 < 1,
\end{equation}
but we are reasonably 
confident that $0.134 < \epsilon_0 < 0.135$.  
We leave it those
with greater expertise with computational methods 
to find better rigorous upper and lower 
bounds for $\epsilon_0$. 

We recall that
\[\epsilon_0 = \inf \{ \epsilon \in (0,1] \,\,|\,\, \alpha(\epsilon) = 1 \}\]
where
\[ \alpha(\epsilon) = \sup_{0 \le z \le 4\pi} w_\epsilon(z) \]
where
\[ w_\epsilon(z) = \frac{1}{\pi^2}\left\{ \begin{array}{cl}  
\frac{\pi}{4}\cdot\epsilon^{-\frac12}\cdot z, & 
0 \le z \le \frac{4\pi}{3-2\epsilon} \\ 
\int_0^{y(z)} 
\left(36\pi - 27 (1-\epsilon) y(z)^\frac23 - 9\epsilon\cdot
x^\frac23  \right)^{-\frac12} dx & \\ 
+ \int_{y(z)}^{z^\frac32} 
(36\pi - 18(1-\epsilon)y(z) x^{-\frac13} - 9 x^{\frac23})^{-\frac12} dx, &
\frac{4\pi}{3-2\epsilon} < z \le  4\pi 
\end{array}\right.\]  
where 
\[y(z) = \frac{z^\frac12 (4\pi - z)}{2(1-\epsilon)}. \]
We recall that 
since $w_\epsilon(z)$ is continuous, the maximum value must occur for 
$z \in [\frac{4\pi}{3-2\epsilon}, 4\pi]$.  In fact, since
$w_\epsilon(z)$ is $C^1$, 
the maximum value can not occur on the 
left end point of this interval, $z=\frac{4\pi}{3-2\epsilon}$,
although there is a local maximum very close to this point 
when $\epsilon$ is less than about $0.2$.  In fact, we find
that the maximum value of $w_\epsilon(z)$ either occurs
at $z=1$ or at a $z$ value only slightly greater than
$\frac{4\pi}{3-2\epsilon}$.  This phenomenon can be seen 
in figure \ref{epsilon1} where we can recognize the 
location of $z=\frac{4\pi}{3-2\epsilon}$ on the graph 
using the fact that
$w_\epsilon(z)$ is linear for 
$0 \le z \le \frac{4\pi}{3-2\epsilon}$. 

As previously mentioned, it is easily shown that
$w_\epsilon(z)$ is a decreasing function of $\epsilon$ when
$z$ is held fixed.  Thus, from figure \ref{epsilon1} we 
see that the maximum value of $w_\epsilon(z)$ is greater
than one when $\epsilon = 0.1$, 
so this must be the case for all $\epsilon < 0.1$ as well.
Hence, $\epsilon_0 > 0.1$.  Using this idea again we 
conclude from figure \ref{epsilon2} that $0.13 < \epsilon_0
< 0.14$, and continuing this procedure is
how we estimated that $\epsilon_0 \approx 0.134727$.  

In figure \ref{epsilon3} we see the graph of $w_\epsilon(z)$
when $\epsilon = 0.134727$, so that the maximum value of  
$w_\epsilon(z)$ is roughly one and occurs (to the accuracy
of the computer) at two $z$ values, $z = 1$ and $z$ slightly
greater than $\frac{4\pi}{3-2\epsilon}$.  Figure 
\ref{epsilon4} is an enlargement of figure \ref{epsilon3}
around this second maximum.

Notice that from the form of the formula for $w_\epsilon(z)$
that solving for the explicit values of the critical 
points using $w_\epsilon'(z) = 0$
seems very difficult, and this is why it seems necessary
to resort to numerical computations.

The first integral in the formula for $w_\epsilon(z)$ can
be computed in closed form.  However, we 
used Simpson's rule to estimate the integral in the 
formula for $w_\epsilon(z)$.  To use Simpson's rule, we 
need the function we are integrating to be bounded, so we 
subtract the function 
$k (1 - \frac{x}{z^{3/2}})^{-1/2}$ for some $k$ 
from the second integrand to make it a bounded function.   
We then use Simpson's rule with $N$, the number of intervals,
equal to one thousand.  We have not attempted any rigorous 
error estimates, although we have observed that the value
of $\epsilon_0$ which we compute is the same to six digits
for $N = 100$, which is a good sign.

To get the rigorous bounds in inequality \ref{rigorous}, 
we note that
$w_\epsilon(\frac{4\pi}{3-2\epsilon}) = \frac{1}
{\epsilon^\frac12 (3-2\epsilon)} = 1$ when $\epsilon = 
1-\sqrt{3}/2$.  Since this endpoint is never the maximum value
of $w_\epsilon(z)$ since $w_\epsilon(z)$ is $C^1$ and has positive
slope at $z =\frac{4\pi}{3-2\epsilon}$, $\alpha(1-\sqrt{3}/2) > 1$.  
Hence, $\epsilon_0 > 1 - \sqrt{3}/2$.  Finally, to show that
$\epsilon_0 < 1$, it is sufficient to prove that for some 
$\epsilon < 1$, $w_\epsilon'(z) \ge 0$ which implies that the 
maximum value of $w_\epsilon(z)$ occurs at $z=4\pi$ and equals
$1$.  Choosing $\epsilon$ very close to $1$ we find that this 
is true, although the computations are not trivial.  Thus, 
$\alpha(\epsilon) = 1$ for some $\epsilon < 1$, so 
$\epsilon_0 < 1$.

\section{Conjectures}

The most natural generalization of theorem \ref{vcthm} is to
propose that it is true in higher dimensions. 

\begin{conjecture} \label{higher_dim}
Let $(S^n,g_0)$ be the constant curvature metric on $S^n$ with scalar 
curvature $R_0$, Ricci curvature $Ric_0 \cdot g_0$, and volume $V_0$.
Then for each $n \ge 3$, there exists a positive 
$\epsilon_0(n) < 1$ such that
if $(M^n,g)$ is any complete smooth Riemannian manifold with volume
$V$ satisfying
\[R(g) \ge R_0 \] 
\[Ric(g) \ge \epsilon_0(n) \cdot Ric_0 \cdot g\] 
then 
\[ V \le V_0. \]
\end{conjecture}

Other problems which relate to scalar
curvature include questions connected to the Yamabe problem and 
Einstein metrics \cite{S}.  Given a manifold $M^n$, 
consider metrics of volume one and define
the energy to be the integral of scalar curvature.  Einstein metrics are
critical points of this functional.  One approach to finding critical points
of energy is to define $I(g)$ to be the infimum of the energy of all metrics
conformal to $g$, and then to define $\sigma(M)$ to be the supremum of
$I(g)$ over all conformal classes of metrics.  If $\sigma(M) \le 0$, then
it is known that $I(g)$ is always realized by a unique metric, so that 
$\sigma(M)$ is realized by a metric which is a critical point of the 
energy functional and hence is Einstein.  
However, for $\sigma(M) > 0$, it is not known under what 
circumstances this procedure
yields a critical point of the energy functional.

Also, if $M$ is a manifold
which admits a constant curvature metric, then it is conjectured 
by Schoen \cite{S} that the above procedure 
produces the constant curvature metric and that $\sigma(M)$ equals the
energy of this metric.  Schoen's conjecture splits
naturally into two cases, depending on whether the constant curvature
metric is negatively curved or positively curved.  Considering these two cases
separately motivates the following two conjectures. 
 
\begin{conjecture}{\bf (Schoen)}   \label{R10}
Suppose $M^n$, $n \ge 2$, admits a hyperbolic metric $g_0$ 
with constant negative 
scalar curvature $R_0$.  If $g$ is any other metric on $M^n$ 
with $R(g) \ge R_0$, then $Vol(g) \ge Vol(g_0)$.       
\end{conjecture}

\begin{conjecture}     \label{R11}
Let $(S^n,g_0)$, $n \ge 2$, 
be the standard constant curvature metric on $S^n$ with 
first nonzero eigenvalue
of the Laplacian operator $\lambda_0$.  
Let $G$ be any finite isometric group action on $(S^n,g_0)$ 
without fixed points, so that
$(M^n,g_0) = (S^n,g_0) / G$ is a constant curvature metric on $M^n$ 
with scalar curvature $R_0$ and volume $V_0$.
If $g$ is a metric on $M^n$ with $R(g) \ge R_0$ and first eigenvalue
$\lambda(g) \ge \lambda_0$, then $Vol(g) \le V_0$.
\end{conjecture}

Conjectures \ref{R10} and \ref{R11} imply Schoen's conjecture
respectively in the negatively and positively curved cases.  (In the case
that $M^n$ admits a flat metric the conjecture is already known to be true.)
Furthermore, if either conjecture \ref{R10} or \ref{R11} turns out to
be false, then there would be a good chance that a counterexample
to Schoen's conjecture could be found.

Conjectures \ref{higher_dim} and \ref{R11} have the similarity
that both attempt to use a lower
bound on scalar curvature to achieve an upper bound on the 
total volume.  However, both conjectures are false without
additional assumptions.  For conjecture \ref{higher_dim}, we
need a lower bound on the Ricci curvature, and for 
conjecture \ref{R11} 
we need a lower bound on the first nonzero eigenvalue. 
These last two inequalities
are weak in the sense that they are not equalities for the 
constant curvature metrics (unless $G$ is trivial in conjecture
\ref{R11}).
Hence, both conjectures say that
for metrics close to the constant curvature metric (on $S^n$
in conjecture \ref{higher_dim} and on $S^n/G$ in conjecture
\ref{R11}, for nontrivial $G$) 
that $R \ge R_0$ implies $V \le V_0$.
Conjecture \ref{R10}, on the other hand, is a volume comparison
conjecture for scalar curvature for hyperbolic metrics, and is 
particularly compelling because of its simplicity.

\appendix
\chapter{Some Geometric Calculations}

Let $\Sigma^2$ be a smooth compact surface without boundary in $(M^3,g)$.  
In this appendix we compute the rate of change of the mean 
curvature and the area form of $\Sigma^2$ given a smooth
variation of $\Sigma^2$.  
We
define a variation of $\Sigma^2$ as follows.  For $-\epsilon<t<\epsilon$ and
$x\in \Sigma^2$, suppose $\Sigma^2(x,t)$ takes values in $M^3$, is smooth,
$\Sigma^2(t)=\{\Sigma^2(x,t)|x\in\Sigma^2\}$ is a smooth family of surfaces
around $\Sigma^2$, and the vector $\frac{\partial\Sigma^2(x,t)}{\partial t}$
is perpendicular to $\Sigma^2(t)$ at $\Sigma^2(x,t)$.  Let $\vec{\mu}(x,t)$
be the outward-pointing unit normal to $\Sigma^2(t)$ at $\Sigma^2(x,t)$, so
that we must have
\begin{equation}
\frac{\partial \Sigma^2(x,t)}{\partial t} = \eta(x,t)\vec{\mu}(x,t)
\label{eqn:flow}
\end{equation}
for some real-valued function $\eta(x,t)$.  Then we see that the surfaces
$\Sigma^2(t_0)$ can be thought of as the surface created by starting
at $\Sigma^2$ and flowing in the outward unit normal (to
$\Sigma^2(t)$) direction at speed $\eta(x,t)$ for $t$ between $0$ and
$t_0$.  We call $\eta(x,t)$ the flow rate.  In fact, given any smooth flow rate
$\eta(x,t)$, for $x\in\Sigma^2$ and $t\in(-\delta,\delta)$, we can always
find a smooth mapping $\Sigma^2(x,t)$ as above satisfying equation
(\ref{eqn:flow}) such that $ \Sigma^2(t)$ is a smooth family of surfaces
around $\Sigma^2$, for $t\in(-\epsilon,\epsilon)$, for some $\epsilon>0$.

Let $du(x)$ be the area form on $\Sigma^2$, $\pi(x)$ be the second fundamental
form of $\Sigma^2$ in $M^3$ at $x$, and $H(x)=\mbox{trace}(\pi(x))$ be the
mean curvature of $\Sigma^2$ at $x$.  Let $d\mu(x,t)$ be the area form on
$\Sigma^2(t)$, $\pi(x,t)$ be the second fundamental form of $\Sigma^2(t)$
in $M^3$ at $\Sigma^2(x,t)$, and $H(x,t)=\mbox{trace}(\pi(x,t))$ be the
mean curvature of $\Sigma^2(t)$ at the point $\Sigma^2(x,t)$, for
$t\in(-\epsilon,\epsilon)$.  In this section we will verify the formulas
\begin{equation}
\ddt d\mu(x,t) = H(x,t)\eta(x,t) d\mu(x,t)
\label{eqn:darea}
\end{equation}
and
\begin{equation}
\ddt H(x,t) = -\Delta_{\Sigma(t)}\eta(x,t)
-\eta(x,t)\|\pi(x,t)\|^2_{M^3}-\eta(x,t)Ric_{M^3}(\vec{\mu}(x,t),\vec{\mu}(x,t))
\label{eqn:dcurve}
\end{equation}
which we will use for important calculations in chapters \ref{Penrose} and
\ref{volume}.

Let $\alpha:U\to\Sigma^2$ for some $U\subset \real^2$ be a local coordinate
chart for $\Sigma^2$.  Then we can define $\Sigma^2(x,t)$ equivalently locally
on $U\times[-\epsilon,\epsilon]\subset \real^3$ with coordinates $(x_1, x_2,
t)$.  Let $\partial_i$ be the vector $\frac{\partial}{\partial x_i}$,
and define the $2\times 2$ matrix
\[g_{ij}(x_1,x_2,t) = \langle\partial_i,\partial_j\rangle_{M^3}\,\,\, ,
\,\,\, 1\le i,j\le 2\]
where $\langle\cdot,\cdot\rangle_{M^3}$ is the pull-back of the metric of $M^3$ using
the mapping $\Sigma^2(x,t):U\times[-\epsilon,\epsilon]\to M^3$.  Then
$g_{ij}(x_1,x_2,t)$ is the metric for some neighborhood of $\Sigma^2(t)$ so
that
\[d\mu(x_1,x_2,t)=\sqrt{|g(x_1,x_2,t)|}\,dx_1\wedge dx_2\]
where $|g(x_1,x_2,t)|=\det(\{g_{ij}(x_1,x_2,t)\})$.  Computing, we get
\[\ddt\sqrt{|g|}=
\frac{1}{2}|g|^{-1/2}\ddt|g|
=\frac{1}{2}|g|^{1/2}\mbox{trace}(g^{-1}\ddt g)\]
where we have used the formula $\ddt (\det A)
= (\det A)\mbox{trace}(A^{-1}\ddt A)$.
Also,
\begin{eqnarray*}
\ddt g_{ij}&=&\ddt \langle
\partial_i,\partial_j\rangle_{M^3}\\
&=&\langle D_{\partial_t}\partial_i,\partial_j\rangle_{M^3}+
\langle\partial_i,D_{\partial_t}\partial_j\rangle_{M^3}\\
&=&\langle D_{\partial_i}\partial_t,\partial_j\rangle_{M^3}+
\langle\partial_i,D_{\partial_j}\partial_t\rangle_{M^3}
\end{eqnarray*}
since $D_{\partial_i}\partial_t-D_{\partial_t}\partial_i=
[\partial_i,\partial_t]=0$ by the torsion-free property of the connection in
$M^3$ and since $\partial_i$ and $\partial_t$ are coordinate vectors.  Thus,
since by equation (\ref{eqn:flow}) $\partial_t=\eta\vec{\mu}$,
\begin{eqnarray*}
\ddt g_{ij}=\langle D_{\partial_i}(\eta\vec{\mu}),
\partial_j\rangle_{M^3}+\langle\partial_i,D_{\partial_j}(\eta\vec{\mu})
\rangle_{M^3}\\
= \eta\langle D_{\partial_i}\vec{\mu},\partial_j\rangle_{M^3}+
\eta \langle\partial_i,D_{\partial_j}\vec{\mu}\rangle_{M^3}
\end{eqnarray*}
since $\langle\vec{\mu},\partial_i\rangle_{M^3}=0$.  Furthermore, since the
second fundamental form is given by $\pi_{ij}=\langle D_{\partial_i}\vec{\mu},
\partial_j\rangle$ and is symmetric, we have
\begin{equation}
\ddt g_{ij}=2\eta\pi_{ij}.
\label{eqn:dg}
\end{equation}
Thus, putting it all together, we have
\begin{eqnarray*}
\ddt d\mu(x_1,x_2,t)&=&
\ddt \sqrt{|g|}\,dx_1\wedge dx_2\\
&=&\frac{1}{2}|g|^{1/2}\mbox{trace}(g^{-1}\ddt g)\,dx_1
\wedge dx_2\\
&=&\frac{1}{2}|g|^{1/2}\mbox{trace}(g^{ij}2\eta\pi_{jk})\,dx_1\wedge dx_2\\
&=&\mbox{trace}(g^{ij}\pi_{jk})\eta\sqrt{|g|}\,dx_1\wedge dx_2\\
&=&H(x_1,x_2,t)\eta(x_1,x_2,t)\,d\mu(x_1,x_2,t).
\end{eqnarray*}
Thus, equation (\ref{eqn:darea}) is true.

Now we verify equation (\ref{eqn:dcurve}).  Since $H=g^{ij}\pi_{ij}$,
\[\ddt  H = (\ddt  g^{ij})
\pi_{ij} + g^{ij}(\ddt  \pi_{ij}).\]
But since $\ddt (A\cdot A^{-1})=0$, by the product
rule it follows that $\ddt (A^{-1})=
-A^{-1}(\ddt  A)A^{-1}$ so that by equation
(\ref{eqn:dg}),
\begin{eqnarray*}
\ddt H &=& g^{ij}(\ddt \pi_{ij})
-g^{ij}\cdot 2\eta\pi_{jk}\cdot g^{kl}\cdot\pi_{li}\\
&=&g^{ij}(\ddt \pi_{ij})-2\eta\cdot\pi^i_k\pi^k_i\\
&=&g^{ij}(\ddt \pi_{ij})-2\eta\|\pi\|^2_{M^3}.
\end{eqnarray*}
Furthermore,
\begin{eqnarray*}
\ddt \pi_{ij} &=& \ddt \langle D_{\partial_i}\vec{\mu},\partial_j\rangle\\
&=& \langle D_{\partial_t} D_{\partial_i}\vec{\mu},\partial_j\rangle +
 \langle D_{\partial_i}\vec{\mu},D_{\partial_t}\partial_j\rangle\\
&=& -R(\partial_t,\partial_i,\vec{\mu},\partial_j) +
\langle D_{\partial_i} D_{\partial_t} \vec{\mu},\partial_j\rangle +
 \langle D_{\partial_i}\vec{\mu},D_{\partial_t}\partial_j\rangle
\end{eqnarray*}
by the definition of the Riemann curvature tensor.  Then since
$D_{\partial_t}\partial_j=D_{\partial_j}\partial_t$ and $\partial_t
=\eta\vec{\mu}$,
\[\ddt \pi_{ij}=-\eta\cdot R(\vec{\mu},\partial_i,\vec{\mu},\partial_j) +
\langle D_{\partial_i} D_{\partial_t} \vec{\mu},\partial_j\rangle +
 \langle D_{\partial_i}\vec{\mu},D_{\partial_j}\partial_t\rangle.\]
We leave it to the reader to check that $D_{\partial_t}\vec{\mu}
= -\vec{\nabla}_{\Sigma(t)}\eta$.  Furthermore, since $\partial_t
=\eta \vec{\mu}$, and $\langle D_{\partial_j}\vec{\mu},\vec{\mu}\rangle=0$,
\[\ddt \pi_{ij} = -\eta\cdot R(\vec{\mu},\partial_i,\vec{\mu},\partial_j) +
\langle D_{\partial_i}(-\vec{\nabla}_{\Sigma(t)}\eta),\partial_j\rangle +
 \eta \langle D_{\partial_i}\vec{\mu},D_{\partial_j}\vec{\mu}\rangle\]
so that
\[g^{ij}\ddt \pi_{ij} = -\eta Ric(\vec{\mu},\vec{\mu})-\Delta_{\Sigma(t)}\eta
+\eta\|\pi\|^2_{M^3}\]
where $Ric(\cdot,\cdot)$ is the Ricci curvature tensor.

Hence, from before, we have
\[\ddt H = -\Delta_{\Sigma(t)}\eta - \eta\|\pi\|^2_{M^3}-\eta Ric(\mu,\mu)\]
proving equation (\ref{eqn:dcurve}).

One immediate consequence to equation (\ref{eqn:darea}) is that smooth
surfaces which minimize area with a volume constraint must have constant
mean curvature.  Otherwise, we consider a flow on the surface $\Sigma$ with
a flow rate $\eta$ defined on $\Sigma$.  Then since the area of $\Sigma(t)$
is
\[A(t)=\int_{\Sigma(t)} d\mu(x,t)\]
we have that
\[A'(0)=\int_{\Sigma}\ddt d\mu(x,0)=\int_{\Sigma} H(x,0)\eta(x,0)\,d\mu(x,0).\]
Furthermore, since
\[V'(0)=\int_\Sigma \eta(x,0)\]
we can find an $\eta(x,0)$ such that $A'(0)<0$ and $V'(0)=0$ unless $H(x,0)$
equals a constant.  Hence, any smooth surface which even locally minimizes
area among surfaces containing the same volume must have constant mean
curvature.


\begin{thebibliography}{99}
\bibliographystyle{plain}
\bibliography{mybib}


\bibitem{ADM} R. Arnowitt, S. Deser and C. Misner, ``Coordinate
Invariance and Energy Expressions in General Relativity,'' {\em
Phys. Rev.} {\bf 122} (1961) 997-1006.

\bibitem{Ba2} R. Bartnik, ``The Mass of an Asymptotically Flat
Manifold,'' {\em Comm. Pure Appl. Math.} {\bf 39} (1986) 661-693.

\bibitem{Ba4} R. Bartnik, ``New Definition of Quasi-Local 
Mass,'' {\em
Phys. Rev. Lett.} {\bf 62} (1989) 2346.

\bibitem{Ba3} R. Bartnik, ``Quasi-Spherical Metrics and Prescribed
Scalar Curvature,'' {\em J. Diff. Geom.} {\bf 37} (1993) 31-71.

\bibitem{C} D. Christodoulou, ``Examples of Naked Singularity
Formation in the Gravitational Collapse of a Scalar Field,'' {\em
Ann. of Math.} {\bf 140} (1994) 607-653.

\bibitem{CY} D. Christodoulou and S.-T. Yau, ``Some 
Remarks on the Quasi-Local Mass,'' {\em Contemporary 
Mathematics} {\bf 71} (1988) 9-14.


\bibitem{G} R. Geroch, ``Energy Extraction,'' {\em Ann. New York
Acad. Sci.} {\bf 224} (1973) 108-17.

\bibitem{Gi} G. Gibbons, ``Collapsing Shells and the Isoperimetric
Inequality for Black Holes,''  {\em Univ. of Cambridge}, preprint,
1997.

\bibitem{H1} S. W. Hawking, ``Gravitational Radiation in an 
Expanding
Universe,'' {\em J. Math. Phys.} {\bf 9} (1968) 598-604.

\bibitem{H2JW} S. Hawking, {\em Phys. Rev. Lett.} {\bf 26}, 
1344 (1971).

\bibitem{H2} S. W. Hawking, ``Black Holes in General 
Relativity,'' {\em
Comm. Math. Phys.}, {\bf 25} (1972) 152-166.

\bibitem{HE} S. W. Hawking and G. F. R. Ellis, {\em The Large-Scale
Structure of Space-Time}, Cambridge University Press, Cambridge, 
1973.

\bibitem{HP} S. W. Hawking and R. Penrose, ``The Singularities of
Gravitational Collapse and Cosmology,'' {\em Proc. Roy. Soc. A} {\bf
314} (1970) 529-548.

\bibitem{He} M. Herzlich, ``A Penrose-like Inequality for the 
Mass of
Riemannian Asymptotically Flat Manifolds,'' to appear in
{\em Comm. Math. Phys.}

\bibitem{HI} G. Huisken and T. Ilmanen, ``Proof of the Penrose
Inequality (Announcement).''

\bibitem{HY} G. Huisken and S.T. Yau, ``Definition of Center of Mass
for Isolated Physical Systems 
and Unique Foliations by Stable Spheres
with Constant Mean Curvature,'' {\em Invent. Math.} {\bf 124} (1996)
281-311.

\bibitem{IJW} W. Israel, {\em Phys. Rev.} {\bf 164}, 1776 (1967); 
{\em Comm. Math.}
Phys. {\bf 8}, 245 (1968).

\bibitem{J1} P.S. Jang, ``On the Positive Energy Conjecture,''
{\em J. Math. Phys.} {\bf 17} (1976) 141-145.

\bibitem{J2} P.S. Jang, ``On the Positivity of Energy in General
Relativity,'' {\em J. Math. Phys.} {\bf 19} (1978) 1152-1155.

\bibitem{J3} P.S. Jang, ``Note on Cosmic Censorship,'' {\em
Phys. Rev. Lett. D} {\bf 20} (1979) 834-838.

\bibitem{J4} P.S. Jang, ``On the Positivity of the Mass for 
Black Hole
Space-Times,'' {\em Comm. Math. Phys.} {\bf 69} (1979) 257-266.

\bibitem{JW} P.S. Jang and R. M. Wald, ``The Positive Energy
Conjecture and the Cosmic Censor Hypothesis,'' {\em J. Math. Phys.}
{\bf 18} (1977) 41-44.

\bibitem{P} R. Penrose, ``Naked Singularities,'' {\em Ann. New York
Acad. Sci.} {\bf 224} (1973) 125-134.

\bibitem{RJW} D. Robinson, {\em Phys. Rev. Lett.} {\bf 34}, 
905 (1975).

\bibitem{S} R. Schoen, ``Variational Theory for the Total 
Scalar Curvature Functional for Riemannian Metrics and Related
Topics,''{\em Topics in Calculus of Variations (M. Giaquinta, ed.)
Lecture Notes in Math.}, {\bf 1365}, 120-154, 
Springer, Berlin, 1987.

\bibitem{SY1} R. Schoen and S.-T. Yau, ``Incompressible Minimal
Surfaces, Three-Dimensional Manifolds with Nonnegative Scalar
Curvature, and the Positive Mass Conjecture in General Relativity,''
{\em Proc. Nat. Acad. Sci.} {\bf 75} (1978), no. 6, 2567.

\bibitem{SY3} R. Schoen and S.-T. Yau, ``On the Proof of the 
Positive
Mass Conjecture in General Relativity,'' {\em Comm. Math. Phys.} 
{\bf
65} (1979) 45-76. 

\bibitem{SY4} R. Schoen and S.-T. Yau, ``Positivity of the 
Total Mass
of a General Space-Time,'' {\em Phys. Rev. Lett.} {\bf 43} (1979) 1457-1459.

\bibitem{SY5} R. Schoen and S.-T. Yau, ``Proof of the Positive Mass
Theorem II,'' {\em Comm. Math. Phys.} {\bf 79} (1981) 231-260.

\bibitem{SY6} R. Schoen and S.-T. Yau, ``The Energy and the Linear
Momentum of Space-Times in General Relativity,'' {\em
Comm. Math. Phys.} {\bf 79} (1981) 47-51.

\bibitem{SY7} R. Schoen and S.-T. Yau, ``The Existence of a 
Black Hole
due to Condensation of Matter,'' {\em Comm. Math. Phys.} {\bf 90}
(1983) 575-579.

\bibitem{T} K. P. Tod. 
{\em Class. Quant. Grav.} {\bf 9} (1992)
1581-1591.

\bibitem{Wi} E. Witten, ``A New Proof of the Positive Energy
Theorem,'' {\em Comm. Math. Phys.} {\bf 80} (1981) 381-402.

\end{thebibliography}
\end{document}